\documentclass[a4paper,12pt]{amsart}
\usepackage{subfigure}
\usepackage{graphicx}
\usepackage{rotating}
\usepackage{morefloats}
\usepackage[utf8]{inputenc}
\usepackage[T1]{fontenc}
\usepackage{amssymb}
\usepackage{amsmath}    
\def\R{I\kern -0,37 em R}
\def\P{I\kern -0,37 em P}
\def\Z{I\kern -0,37 em Z}

\begin{document}
\title[Lie and Cartan Theory, II]{On the Lie and Cartan Theory of Invariant Differential Systems, II}
\author{A. Kumpera}
\address[Antonio Kumpera]{Campinas State University, Campinas, SP, Brazil}
\email{antoniokumpera@gmail.com}

\date{April, 2016}


\maketitle

\begin{abstract}
We start discussing basic properties of Lie groupoids and Lie pseudo-groups in view of applying these techniques to the analysis of Jordan-Hölder resolutions and the subsequent integration of partial differential equations which is the summit of Lie and Cartan's work. Next, we discuss the integration problem for systems of partial differential equations in one unknown function and special attention is given to the first order systems. The Grassmannian contact structures are the basic setting for our discussion and the major part of our considerations inquires on the nature of the Cauchy characteristics in view of obtaining the necessary criteria that assure the existence of solutions. In all the practical applications of partial differential equations, what is mostly needed and what is in fact hardest to obtains are the solutions of the system or, occasionally, some specific solutions. We continue our discussion by examining the local equivalence problem for partial differential equations, illustrating it with some examples, since  almost any integration process or method is actually a local equivalence problem involving a suitable model. We terminate the discussion by inquiring on non-integrable Pfaffian systems and their integral manifolds of maximal dimension. This work is based on four most enlightening \textit{Mémoires} written by Élie Cartan in the beginning of the last century.
\end{abstract}

\section{Introduction}
We would initially like to ask the reader to look at the first paragraph of the \textit{Introduction} in Élie Cartan's \textit{Thesis} (\cite{Cartan1894}) and, in case he should be satisfied with it having no further curiosity, there is no sense at all in reading this manuscript.

\vspace{2 mm}
Differential equations have been studied for a long time since Newton, Leibniz and Monge. Nevertheless we owe, first to Marius Sophus Lie and later to Élie Joseph Cartan, the structuring of this problem with the introduction of the \textit{finite and infinite continuous groups of transformations} leaving invariant these equations (\cite{Lie1880,Cartan1894,Cartan1904,Cartan1905}). It is however worthwhile to observe (\cite{Kumpera1980,Kumpera1999,Kumpera2015}) that whereas Lie was most incredibly successful in dealing with equations of finite type (the solutions depending only upon a finite number of parameters, \cite{Lie1873,Lie1878,Lie1879}) for he had a deep knowledge and understanding of finite continuous groups, in the general context he was unable to go any further beyond examining a few specific situations involving the simple infinite and transitive continuous groups. He showed, for example, that Jacobi’s last multiplier method was the best possible integration method - shortest, most accurate and of lowest degree - due to the fact that the infinite group of all volume preserving transformations is simple. He was aware of the four classes of (complex) transitive simple pseudo-groups of transformations but rather uneasy whether these were the only ones, thus preventing him to take any benefits stemming from a systematic use of Jordan-Hölder resolutions (\cite{Kumpera1980}). Consequently, it remained to Cartan to develop the infinite dimensional theory with a touch only accessible to the most illuminated. In his doctoral dissertation  (\cite{Cartan1894}), Cartan classified the transitive complex simple Lie pseudo-groups showing that these were precisely the four classes announced by Sophus Lie (see also \cite{Cartan1909} and \cite{Kuranishi1962}). Much later (\cite{Cartan1913}), he obtained the classification of the real forms and could start definitely working on what really mattered to him namely, the devising of integration processes for differential systems (partial differential equations) with the help of Jordan-Hölder resolutions for the Lie pseudo-groups of all the local equivalences, defined on the base 
spaces, of the given systems (\cite{Cartan1904, Cartan1905, Cartan1908, Cartan1909, Cartan1921, Cartan1930, Cartan1938/39})\nocite{Goursat1922}. We also 
find it appropriate to mention here Cartan’s covariant approach to Lie pseudo-groups theory. A Cartan pseudo-group, in the terminology 
of Kuranishi and Matsushima (\cite{Kuranishi1961,Matsushima1953}), is a pseudo-group of local transformations defined on an 
analytic manifold \textit{M} of dimension \textit{n} and satisfying the following conditions (\cite{Kumpera1963}):
 
\vspace{2 mm}
\noindent
1. There exist $m+p=n$ (point-wise) independent analytic linear differential forms $\{\omega^i,\varpi^{\lambda}\},$ $1\leq i\leq m~,$ defined on \textit{M} and that satisfy the structure equations
\begin{equation*}
d\omega^i=1/2~c^i_{jk}\omega^j\wedge\omega^k + a^i_{j\lambda}\omega^j\wedge\varpi^{\lambda}~,
\end{equation*}
where the coefficients are analytic functions, the $c^i_{jk}$ are skew-symmetric in the lower indices, some of the 
forms $\varpi^{\lambda}$ are a basis for the first integrals of the invariance fibering of the pseudo-group action and 
the matrices $a_{\lambda}=[a^i_{j\lambda}]$ are linearly independent at each point.

\vspace{2 mm}
\noindent
2. A local transformation $\varphi$ belongs to the pseudo-group if and only if it preserves each $\omega^i$ \textit{i.e.}, $\varphi^*\omega^i=\omega^i$ for all indices \textit{i}. 

\vspace{3 mm}
A Cartan pseudo-group is a Lie pseudo-group of order one and we observe that the structure equations provide all the information required for the analysis of the so-called \textit{Lie co-algebroid}. Furthermore, we also observe that the $c^i_{jk}$ provide the structure of the corresponding first order Lie algebroid, which can also be defined in the obvious way with the help of the forms $\omega^i~,$ and the matrices $a_{\lambda}$ describe the isotropy algebras at the various points. It should finally be mentioned that any Lie pseudo-group of order \textit{k} can be replaced by an equivalent Cartan pseudo-group since it suffices, for this, to prolong the initially given pseudo-group to the manifold $\Pi_kM~,$ the groupoid of all the invertible $k-$jets, restrict this prolongation to the manifold of the differential equation that defines this $k-th$ order pseudo-group and finally choose an appropriate free set of contact 1-forms in restriction to the equation. Since, for a Cartan pseudo-group, all the finite as well as the infinitesimal data are \textit{contained in the memory of the forms $\omega^i$ coupled with the structure equations}, it is not surprising at all to learn that Lie's Third Theorem is true for Cartan pseudo-groups since we can, actually, define the Cartan algebroid directly with the help of the same $1-$forms $\omega^i~,$ where invariance will now mean the vanishing of these forms under the Lie derivatives of the candidates to membership in the Cartan algebroid Club (\cite{Kumpera1963}). Quite to the contrary, it is most surprising to learn that this feature does not carry over to the Lie pseudo-groups and pseudo-algebras, the problem that arises being the appearance of an undesirable non-trivial \textit{holonomy groupoid} that becomes present when we try to project the Cartan finite and infinitesimal transformations down to the Lie finite and infinitesimal counterpart\footnote{The appropriate setting for this discussion should be held on the level of Lie groupoids and Lie algebroids defined in Jet spaces.} (\cite{Almeida1985}).

\vspace{2 mm}
As for the notion of \textit{invariant} differential 1-form with respect to an integrable Pfaffian system $\mathcal{P}$ (\textit{invariant intégral absolu} in the sense of Cartan, \cite{Cartan1921}), we refer the reader to \cite{Kumpera1999}. These forms play a fundamental role in the integration of \textit{S} by quadratures, are invariant under the finite and infinitesimal automorphisms of $\mathcal{P}$ and, moreover, are the invariant forms, in the sense above, of the pseudo-group (resp. pseudo-algebra) of all the local automorphisms, finite (resp. infinitesimal) of the Pfaffian system $\mathcal{P}~.$ Consequently, these pseudo-groups and pseudo-algebras are both \textit{Cartan}.

\vspace{2 mm}
Many years ago, Jordan-Hölder sequences for finite dimensional Lie groups were examined in \cite{Kumpera1980} and rather curious and surprising situations were unveiled, the results having been later applied to the local equivalence problem for contact stuctures (\cite{Kumpera1974}) and the integration of equations of finite type (\cite{Kumpera1999,Kumpera2002}). In particular, \textit{multi-flag} and \textit{truncated multi-flag} systems are discussed in the latter reference showing their relevance in the integration problem for under-determined systems of ordinary differential equations. As mentioned earlier, differential equations (transitively) invariant under the action of a finite dimensional Lie group are a very restricted subject and, in order to examine arbitrary equations of any order, we are forced to consider much broader structures. In our present context, we shall consider the differentiable categories and groupoids introduced by Ehresmann (\cite{Ehresmann1958}) and study, in a first step, the Jordan-Hölder resolutions for these objects though what really matters is the corresponding knowledge for Lie pseudo-groups. We shall also make an attempt to apply these resolutions towards the integration problem for arbitrary differential systems and this will lead us to take a glance at one of the most beautiful, fascinating and surprising articles ever written by Élie Cartan (\cite{Cartan1910}, see also \cite{Cartan1911} and \cite{Cartan1914}). It stands as a gift, to the reader, of the most authentic $Parisian~hors-d'o\!euvre.$     

\vspace{2 mm}
The differential properties of Lie groupoids are subordinate to an algebraic structure namely, a groupoid structure. However, the restrictions thus imposed are much weaker than in the case of Lie groups since, in general, there is a small part of the groupoid, namely the subset of units, which does not carry any algebraic structure at all, except a trivial action of the groupoid itself, and the composition in the groupoid is only partially defined. As in the case of Lie groups, many algebraic and differential properties of the Lie groupoids can be described as well as derived in terms of the corresponding algebraic properties of their Lie algebras - called Lie algebroids in the present case - thus giving rise to functors from an algebraic to a differential category. The local results are analogous to those for Lie groups, their proofs just requiring a little more patience and effort. Quite to the contrary, the global results are sensibly more sophisticated and certainly more difficult to obtain. \nocite{Kumpera1999} \nocite{Cartan1894}

\vspace{2 mm}
The importance of Lie groupoids relies more in their usefulness as a technique to study geometrical problems than as mathematical objects \textit{per se}. Moreover, they constitute at present the most appropriate setting for the development of certain geometric theories \textit{e.g.}, infinitesimal connexions and Lie pseudo-groups. They find inasmuch their place in the local and global equivalence problem for geometrical structures (\cite{Kumpera2015},\cite{Kumpera2016}) and, most important, in the study of differential systems and partial differential equations. We finally observe that the study of G-structures can be considered as a special case of the study of transitive Lie groupoids and, in the bibliography, we mention several articles that apply Lie groupoids in various differential geometric situations.

\vspace{2 mm}
Most of the definitions and concepts adopted here are standard and can be found in \cite{Kumpera1972} and \cite{Kumpera1980} where, in the second reference, we shall replace the term \textit{group} by that of \textit{groupoid}. It should however be observed that the definition of a \textit{differentiable groupoid} given in \cite{Kumpera1972} relies on a slightly weaker regularity assumption than the original definition given by Ehresmann in, for instance, \cite{Ehresmann1967}. The former just requires the transversality of the product map $source\times target$ with respect to the diagonal in the base space whereas the later requires the submersivity of any one of the above two maps (the same holding, consequently, for the other). The difference being so insignificant, the reader should not bother about it. In the present paper, we limit ourselves to discuss, in what concerns the basic definitions and properties related to Lie groupoids and algebroids, just the sub-structures and quotients since much can be found in the appendix of the first reference, in \cite{Kumpera2017} as well as in \cite{Kumpera2016} and finally in \cite{Kumpera1971}. As for the terminology, we do not have any privileged choices and find that the best option results in looking at the not always too explicit origins. We finally remark that all the spaces and structures considered here are \textit{small} i.e., defined on underlying sets and never on \textit{classes} as considered, in a much broader setting, by Ch. Ehresmann (\cite{Ehresmann1958}) \nocite{Pradines2006} \nocite{Ehresmann1966}.

\vspace{2 mm}
The paper is essentially divided in four parts The first one (sect. $1-5$) deals with Lie groupoids, Jordan-Hölder sequences and their application to the integrability problem for partial differential equations. It should be noted that, in this first part, we just provide the strictly necessary proofs. The section 2 has a full proof showing that the quotient groupoid is \textit{Lie}. As for the section 3, there is no sense in providing proofs for statements already proven by Élie Cartan. Nevertheless, we do provide many hints that exhibit the main facts underlying these proofs. The section 4 provides important motivational examples since Lie's most ingenious integration method is essentially based on reducing the problem to the integration of equations invariant under simple Lie pseudo-groups or, inasmuch, simple Lie groupoids, once the choice of a Jordan-Hölder sequence has been achieved. Concerning the section 5, a detailed description is given showing the steps to be taken in order to implement the Lie integration method. As for the subsequent sections, proof are provided whenever necessary though similar proofs are not repeated. The second part (sect. $6-17$) studies the integrability of partial differential equations, in one unknown function, in the context of contact structures. The third part (sect. $18-20$) puts forward the local equivalence problem for differential equations and finally the fourth part (sect. $21-29$) exposes Cartan's general methods enabling to investigate non-integrable Pfaffian systems, that do appear in the second part, and the obtainment of maximal dimensional integral sub-manifolds of such Pfaffian systems via the Jordan-Hölder setup, as described in the first part.
\vspace{6 mm}
\begin{equation*}
CONTENTS
\end{equation*}

\vspace{3 mm}
2. Sub-groupoids and quotients.

\vspace{2 mm}
3. Pfaffian systems in five variables and the Monge problem.

\vspace{2 mm}
4. On the integration of some differential systems transitively invariant under simple groupoids.

\vspace{2 mm}
5. The Jordan-Hölder resolution and the integration of partial differential equations.

\vspace{2 mm}
6. The Lie and the Cartan layouts.

\vspace{2 mm}
7. Single second order involutive equations.

\vspace{2 mm}
8. Two and more second order involutive equations.

\vspace{2 mm}
9. The Darboux layout.

\vspace{2 mm}
10. Pfaffian systems, definitions and notations.

\vspace{2 mm}
11. Darboux characteristics.

\vspace{2 mm}
12. Cartan characteristics.

\vspace{2 mm}
13. Classical contact structures.

\vspace{2 mm}
14. Cauchy characteristics.

\vspace{2 mm}
15. Contact transformations.

\vspace{2 mm}
16. Brackets.

\vspace{2 mm}
17. Higher order differential equations.

\vspace{2 mm}
18. The local equivalence problem.

\vspace{2 mm}
19. The local equivalence of differential systems.

\vspace{2 mm}
20. The local equivalence of Lie groupoids.

\vspace{2 mm}
21. Integral contact elements.

\vspace{2 mm}
22. Integral sub-manifolds.

\vspace{2 mm}
23. The derived system.

\vspace{2 mm}
24. The pseudo-group of local automorphisms.

\vspace{2 mm}
25. Pfaffian systems whose characters are equal to one.

\vspace{2 mm}
26. The Cartan co-variant approach.

\vspace{2 mm}
27. Integration of the Pfaffian systems with character equal to one.

\vspace{2 mm}
28. Pfaffian systems whose characters are larger than one.

\vspace{2 mm}
29. Singular Pfaffian systems.

\vspace{2 mm}
30. Examples.

\vspace{2 mm}
31. Equations on Jet bundles.

\section{Normal sub-groupoids and quotients}
In this section we discuss very standard facts though some additional care must be taken in the case of groupoids. We further point out that the subsequent argumentation has already been carried out, to a certain extent and many years ago, by \textit{R. Almeida}, in his doctoral dissertation (\cite{Almeida1980}). A Lie sub-groupoid of $\Gamma$ is, of course, a subset $\Sigma$ stable under the product operation and under inverses. Moreover, it is assumed that $\Sigma$ is a (not necessarily regularly embedded) sub-manifold of $\Gamma$ and that its units space $\Sigma_0$ is a sub-manifold of $\Sigma~.$ A \textit{morphism} $\varphi:\Gamma~\longrightarrow~\Gamma_1$ is a mapping that preserves the composition hence also preserves the inverses since forcibly it must preserve the units. Curiously enough, the image of a morphism is not always a Lie sub-groupoid of the image. On the one hand, the images of two non-composable elements can become composable and one cannot claim anything about their composite. Inasmuch, the image needs not be a sub-manifold. \textit{A priori}, the units space of a sub-groupoid needs not coincide with that of the ambient groupoid and we shall therefore assume right from the beginning the following two

\vspace{3 mm}
\noindent
\textbf{Hypotheses:}

\vspace{2 mm}
1. The units space of a sub-groupoid coincides with the ambient units space.

\vspace{2 mm}
2. Any morphism becomes a diffeomorphism upon restriction to the units.

\vspace{3 mm}
\noindent
Under these conditions, the image of a morphism becomes of course a sub-groupoid. Let $\Gamma$ be a given Lie groupoid and let us take a Lie sub-groupoid $\Sigma~.$ Then clearly, every $\alpha-$fibre of $\Sigma$ is contained in the corresponding $\alpha-$fibre of $\Gamma~.$ Next, for any given unit \textit{e}, we denote by $\Gamma_{\alpha e}$  (resp. $\Sigma_{\alpha e}$) and $\Gamma_e$ (resp. $\Sigma_e$) the $\alpha-$fibres and the isotropy groups at this unit element. The union of all the isotropy groups needs not be a differentiable sub-manifold of $\Gamma~.$ However, these isotropies will become a locally trivial differentiable sub-bundle of, say, the $\alpha-$fibres or the $\beta-$fibres bundle, if and only if the fibre product mapping 
\begin{equation*}
\alpha\vee\beta:g\in\Gamma~\longmapsto~(\alpha(g),\beta(g))\in\Gamma_0\times\Gamma_0
\end{equation*}
is transverse to the diagonal in $\Gamma_0\times\Gamma_0~,$ where $\Gamma_0$ is the units space of $\Gamma$ and where by "transverse" we mean that the sum of the tangent space to the diagonal with the image of the tangent space to $\Gamma$ by the map  $(\alpha\vee\beta)_*$ is equal to the tangent space of the product. We shall say that $\Sigma$ is \textit{normal} or \textit{invariant} in $\Gamma$ whenever its isotropy bundle is \textit{normal} in $\Gamma~.$ A little explanation seems to be required. What we are seeking for is the usual condition of normality namely, $\gamma\cdot\Sigma\cdot\gamma^{-1}\subset\Sigma$ for an arbitrary element $\gamma\in\Gamma~.$ Since $\alpha(\gamma)=\beta(\gamma^{-1})$ and $\beta(\gamma)=\alpha(\gamma^{-1})~,$ we infer that the previous condition is only applicable to isotropic elements of $\Sigma$ and, moreover, that the resulting subset $\gamma\cdot\Sigma_e\cdot\gamma^{-1},~e=\alpha(\gamma)~,$ is the isotropy $\Sigma_{e'}$ of $\Sigma$ at the point $e'=\beta(\gamma)~.$ The above normality requirement means consequently that each subgroup $\Sigma_e$ is normal in $\Gamma_e$ and, furthermore, that the isotropy sub-bundle of $\Sigma$ is preserved by the (adjoint) action of arbitrary elements $\gamma\in\Gamma~.$ By inversibility, the previous inclusion becomes in fact an equality.

\vspace{2 mm}
We are now all set to define the quotient $\Gamma/\Sigma$ in the category of Lie groupoids. However, instead of defining the equivalence relation that will factor $\Gamma$ and ultimately produce the quotient, we shall first introduce the equivalence classes. Of course, we can consider either the \textit{right} classes $\Sigma\cdot\gamma=\Sigma_{\beta(\gamma)}\cdot\gamma$ or the \textit{left} classes $\gamma\cdot\Sigma=\gamma\cdot\Sigma_{\alpha(\gamma)}~.$ On account of the inversion $\gamma\longmapsto\gamma^{-1}~,$ both procedures are equivalent and, on account of normality, both procedures will turn out to be equal. Since every $\alpha-$fibre of $\Sigma$ is contained in the corresponding $\alpha-$fibre of $\Gamma~,$ we would like that each $\alpha-$fibre $\Sigma_{\alpha e}$ be an equivalence class that factors onto (identify with) the corresponding unit in the \textit{eventual} quotient groupoid. However, this cannot be so since the elements belonging to an $\alpha-$fibre can have many targets. The same problem appears with respect to the $\beta-$fibres. We shall therefore have to argue with the isotropy groups.

\vspace{2 mm}
Let us, at present, just consider the right classes (right co-sets) and, for any element $\gamma\in\Gamma~,$ we set $\Sigma_e\cdot\gamma~,$ where $e=\beta(\gamma)$ ($\Sigma_e=\Sigma_e\cdot e$). We first show that two such classes either coincide or are disjoint. In fact, assuming that $Z$ belongs to both $\Sigma_e\cdot Y$ ($e=\beta(Y)$) and $\Sigma_{e'}\cdot Y'~,$ then $Z=X\cdot Y=X'\cdot Y'~,$ $\alpha(Z)=\alpha(Y)=\alpha(Y')$ and  $\Sigma_e\cdot Y\cdot(Y')^{-1}=\Sigma_{e'}$ hence ultimately $\Sigma_e \cdot Y=\Sigma_{e'}\cdot Y'~.$ Next, we define the \textit{equivalence relation} on $\Gamma$ whose equivalence classes are precisely the above co-sets. If $Z,~Z'$ both belong to the co-set $\Sigma_e\cdot Y~,$ then $Z=X\cdot Y$ and $Z'=X'\cdot Y$ with $X,X'\in\Sigma_e$ whereupon $Z'\cdot Z^{-1}=X'\cdot X^{-1}\in\Sigma_e~.$ We infer that $Z\sim Z'$ if and only if $Z'\cdot Z^{-1}\in\Sigma~,$ as in group theory. 

\vspace{2 mm}
We next define the groupoid structure on the quotient $\Gamma/\Sigma~,$ composed by the above co-sets, and assume that $\Sigma$ is topologically closed in $\Gamma$ otherwise the quotient topology would be inadequate to underlie a differentiable manifold structure. It follows thereafter that the Lie sub-group $\Sigma_e$ is also closed in $\Gamma_e~.$ Setting $\alpha(\Sigma_e\cdot\gamma)=\alpha(\gamma)~,$ $\beta(\Sigma_e\cdot\gamma)=e$ and taking any two composable co-sets $\Sigma_e\cdot\gamma$ and $\Sigma_{\epsilon}\cdot\delta~,$ $e=\beta(\gamma),~\epsilon=\beta(\delta),~\alpha(\gamma)=\beta(\delta),$ we define the product by 
\begin{equation*}
(\Sigma_e\cdot\gamma)\cdot(\Sigma_{\epsilon}\cdot\delta)=\Sigma_e\cdot\gamma\cdot\delta~,
\end{equation*}
and proceed to show that this operation does not depend upon the representatives. Assuming that $\Sigma_e\cdot\gamma=\Sigma_{e'}\cdot\gamma'~,$ then $\gamma'=X\cdot\gamma~,$ with $X\in\Sigma_e~,$ and consequently $\Sigma_{e'}\cdot\gamma'=\Sigma_e\cdot\gamma'$ hence,
\begin{equation*}
(\Sigma_e\cdot\gamma')\cdot(\Sigma_{\epsilon}\cdot\delta)=\Sigma_e\cdot\gamma'\cdot\delta=\Sigma_e\cdot X\cdot\gamma\cdot\delta=\Sigma_e\cdot\gamma\cdot\delta~,
\end{equation*}
since $\Sigma_e\cdot X=\Sigma_e~.$ It should be observed that this first step does not require the normality condition on $\Sigma~.$ Let us next assume that $\Sigma_{\epsilon}\cdot\delta=\Sigma_{\epsilon'}\cdot\delta'~.$ Then $\delta'=X\cdot\delta,~X\in\Sigma_{\epsilon}$ and consequently
\begin{equation*}
(\Sigma_e\cdot\gamma)\cdot(\Sigma_{\epsilon'}\cdot\delta')=\Sigma_e\cdot\gamma\cdot\delta'=\Sigma_e\cdot\gamma\cdot X\cdot\delta=\Sigma_e\cdot\gamma\cdot X\cdot\gamma^{-1}\cdot\gamma\cdot\delta=\Sigma_e\cdot Y\cdot\gamma\cdot\delta~,
\end{equation*}
where $Y=\gamma\cdot X\cdot\gamma^{-1}\in\Sigma_e~,$ hence the above 
product is equal to $\Sigma_e\cdot\gamma\cdot\delta$ since $\Sigma_e\cdot Y=\Sigma_e~.$ Needless to say that the inverse of $\Sigma_e\cdot\gamma$ is equal to $\Sigma_{\epsilon}\cdot\gamma^{-1}$ where $\epsilon=\alpha(\gamma)~.$ 

\vspace{2 mm}
The fact that $\beta\vee\alpha$ is transverse to the diagonal in $M\times M$ guarantees that the quotient $\Gamma/\Sigma$ has a differentiable manifold structure compatible with its groupoid structure. Moreover, it can be shown that the differentiable structure of the sub-groupoid $\Sigma$ is regularly embedded in $\Gamma$ inasmuch as what is true for Lie groups. We shall however skip the details and only remark that the above definition of the quotient does not agree with that outlined in \cite{Pradines1966} and described, more detailed, in \cite{Krill1987} and \cite{Pradines1986},  though both are closely related since, apparently, one cannot escape from operating with isotropies. The advantage of our approach relies on the simplicity in expressing the equivalence classes of given elements $\gamma\in\Gamma~,$ identical to that found in group theory, as well as in  attributing coordinates to the quotient spaces that arise from coordinates defined in the initially given groupoid. However, the main reason for not adopting here the definition found in the above references lies in the fact that such a definition, for the quotient groupoids, enters in conflict with the standard prolongation procedure applied to object defined on finite order jet spaces like, for instance, differential systems and higher order jet groupoids.

\vspace{2 mm}
As for the infinitesimal aspects of the above discussion, there is not much more to be added. It suffices to define the isotropy sub-algebra of $\mathcal{L}$ as being the sub-sheaf composed by those germs whose representatives vanish at the unit element involved. We shall say that $\mathcal{S}$ is a sub-algebroid 
in ideals of $\mathcal{L}$ when each stalk of the isotropy of $\mathcal{S}$ is an ideal in the corresponding stalk of $\mathcal{L}$ at that point. Integrating this condition for each germ of vector field belonging to $\mathcal{L}~,$ we infer that $\mathcal{S}$ is invariant under the local actions associated to arbitrary sections of $\mathcal{L}~.$ The process of defining now a quotient Lie algebroid becomes obvious. It is also self-evident that each isotropy algebra is the Lie algebra of the corresponding isotropy group and that the quotient algebroid is the Lie algebroid associated to the quotient Lie groupoid whenever such an association exists.

\section{Pfaffian systems in five variables and the Monge problem}
In this section we recall the main topics discussed by Cartan in his article\footnote{In \cite{Stormark2000}, the author dedicates the whole chapter 17, as well as some previous sections, to this same \textit{Mémoire}. Essentially all of Cartan's calculations and a few more are indicated there with the sole purpose of devising integration procedure instructions for a large number of differential equations. Curiously enough, the methods employed for solving the equations that appear in \cite{Cartan1910} were already known and well established since Monge (\textit{Gaspard Monge, Comte de Péluse (1746-1818)}) and Cartan wrote this \textit{Mémoire} with the sole purpose of using it as an experimental test-ground where he could examine Lie's ingenious ideas, in all their details, on specific and non-trivial known situations.} (\cite{Cartan1910}) and the first step consists in examining a little closer the Darboux system $\mathcal{D}=\mathcal{D}_3$ on the space $\textbf{R}^3.$ It is the system of rank one generated by the form $\omega=dx^3+x^2~dx^1$ of Darboux class three. The Cartan class of the system is also maximum and equal to three hence its space of characteristics is equal to $\textbf{R}^3$ \textit{i.e.}, its characteristic system is equal to $T^*\textbf{R}^3$ (generated by $\{dx^1,dx^2,dx^3\}$) which means that the system $\mathcal{D}$ cannot be factored locally to a space of lower dimension. Moreover, the first derived system $\mathcal{D}_1$ is null, a property that characterizes $\mathcal{D}~.$ The Darboux system is, up to (local)  equivalence, the only system of rank one and of maximum class three on the space $\textbf{R}^3.$ Its infinitesimal automorphisms are defined by
\begin{equation}
Aut(\mathcal{D})=\{\xi\in\chi(\textbf{R}^3):\vartheta(\xi)\omega\equiv 0~mod~\omega\}~,
\end{equation}
where $\vartheta(\xi)$ denotes the Lie derivative by the vector field $\xi$ and $\chi(\textbf{R}^3)$ the Lie algebra of all the vector fields on $\textbf{R}^3.$ There is a well known Hamiltonian characterization of $Aut(\mathcal{D})$ (\cite{Kumpera2014}) by means of which the Lie bracket of $Aut(\mathcal{D})$ identifies with the Lagrange bracket on $\mathcal{F}~,$ the ring of $C^{\infty}$ functions on $\textbf{R}^3,$ defined by $[f,g]=\xi g-g\zeta f~,$ where $\zeta=H^{-1}(1)$ and \textit{H} is a chosen contact Hamiltonian (up to a multiplicative factor, \textit{cf.} \cite{Kumpera2014}). Since the function \textit{f} is arbitrary, we infer the

\vspace{2 mm}
\newtheorem{lagrange}[LemmaCounter]{Lemma}
\begin{lagrange}
The algebra $Aut(\mathcal{D})$ is simple and operates transitively on $\textbf{R}^3~.$
\end{lagrange}

\vspace{2 mm}
We shall now examine the Pfaffian systems on $\textbf{R}^5$ as discussed in \cite{Cartan1910} and, many years later, in \cite{Kumpera1982}, Sect.11 and \cite{Kumpera1998}. A Pfaffian system is said to be totally regular when all its derived systems are regular. The definitions of the \textit{class} of a system, of the \textit{derived} system and of the \textit{co-variant} system associated to a given system can be found in the last reference and are, of course, borrowed from Élie Cartan. In what follows, we shall only be interested in the systems of rank two and three. Nevertheless, let us just mention that the systems of rank equal to four are all integrable and have the model $\{dx^1,dx^2,dx^3,dx^4\}$ in the standard coordinates. As for the systems of rank equal to one, they are either integrable and have the model $\{dx^1\}$ or else, are the Darboux systems $\mathcal{D}_1=\mathcal{D}$ (a system of $rank~1$ and $class~3=2\times 1+1$ defined on a 5-space) and $\mathcal{D}_2$ (a system of $rank~1$ and $class~5=2\times 2+1$ defined on a 5-space), the latter admitting the model $dx^5+x^4~dx^3+x^2~dx^1$ (the Cartan class of a rank 1 system is always odd).

\vspace{2 mm}
The determination of all the Pfaffian systems of ranks two and three is based on the following results (\cite{Giaro1978},\cite{Kumpera1982}).

\vspace{2 mm}
\newtheorem{co-rank}[PropositionCounter]{Proposition}
\begin{co-rank}
Let S be a Pfaffian system of co-rank equal to 2 with a regular derived system. Then, either S is integrable and its class is equal to the co-dimension of its leaves or else its class is equal to the dimension of M.
\end{co-rank}

\vspace{2 mm}
\newtheorem{rank2}[PropositionCounter]{Proposition}
\begin{rank2}
On a manifold of dimension 5, to each Pfaffian system S of rank 2 with null derived system $S_1=0$ corresponds, in a unique way, a Pfaffian system $\overline{S}$ of rank 3, the co-variant system, determined by the condition $S\subset\overline{S}_1~.$ Furthermore, when $\overline{S}_1$ is regular, then either $\overline{S}_1$ is integrable or else $\overline{S}_1=S.$
\end{rank2}

\vspace{2 mm}
\noindent
\textbf{Remark.} The condition $S\subset\overline{S}_1$ implies that $\overline{S}_1$ is integrable if and only if $\overline{S}$ also is, in which case $\overline{S}=\overline{S}_1~.$

\vspace{2 mm}
\newtheorem{rank3}[PropositionCounter]{Proposition}
\begin{rank3}
On a manifold of dimension 5, the correspondence that to each totally regular Pfaffian system S, of rank 3 and verifying $S_2=0~,$ associates $S_1$ is bijective.
\end{rank3}

\vspace{3 mm}
\noindent
\textbf{Remark.} Following Cartan, we shall, in the above situation, also call $S_1$ the \textit{co-variant} system of \textit{S}.

\vspace{2 mm}
On account of the second proposition, we can say that the Pfaffian systems of rank 2 in a five dimensional space are firstly the integrable systems with the model $\{dx^1,dx^2\}~,$  when $S=S_1~,$ secondly those for which $S_1$ has rank one and finally the case where $S_1=0~.$ The second mentioned systems are, when $S_2=0~,$ the flag systems of length two and are all equivalent to (the pullback of) the Engel flag in four space admitting the model $\{dx^2+x^3dx^1,dx^3+x^4dx^1\}~.$ When $S_1=S_2$ \textit{i.e.}, when $S_1$ is integrable, the system admits the model $\{dx^1,dx^2+x^3dx^4\}~.$ In the third case, when $\overline{S}_1$ is integrable \textit{i.e.}, when $\overline{S}_2=\overline{S}_1~,$ then \textit{S} can be brought to the local form $\{dx^1+x^4dx^3,dx^2+x^5dx^3\}$ and, in the non-integrable case \textit{i.e.}, when $\overline{S}_1=S~,$ the system \textit{S} can be brought to the normal form $\{dx^1+(x^3+x^4x^5)dx^4,dx^2+x^3dx^5\}$ as soon as the Cartan quartic form  (\textit{forme biquadratique binaire}, \cite{Cartan1910}, p.152, \textit{l.}-6) vanishes. The remaining cases can also be examined by making appeal to the Cartan \textit{quartic form} though in these cases we shall find continuous moduli of rank 2 Pfaffian systems that we can describe via their \textit{pseudo-normal forms} (\cite{Kumpera1982},\cite{Kumpera2014}) that will carry parameters (the constants that appear in the pseudo-normal forms). The discussion being rather extensive, we shall limit ourselves to mention just a few results but, firstly, we examine those Pfaffian systems of rank 3 that have discrete models.

\vspace{2 mm}
Let us therefore consider the Pfaffian systems \textit{S} of rank 3 
defined on a 5-dimensional space. There are always either the integrable systems with $S=S_1$ admitting the model $\{dx^1,dx^2,dx^3\}~,$ or those systems for which $S_1=S_2$ \textit{i.e.}, when $S_1$ is integrable, admitting 
the models $\{dx^1,dx^2,dx^3+x^4dx^5\}~,$ that has Cartan class equal to 5, and $\{dx^1,dx^3,dx^2+x^3dx^4\}$ that has class 4, then those systems for 
which $S_2=S_3$ \textit{i.e.}, when $S_2$ is integrable, that admit the models $\{dx^1,dx^2+x^3dx^5,dx^3+x^4dx^5\}~,$ that has class 5, and $\{dx^1,dx^2+x^3dx^4,dx^3+x^4dx^4\}~,$ that has class 4. Next the Cartan flags, when $rank~S_1=2~,~rank~S_2=1$ and $S_3=0~,$ that are all equivalent either to the  homogeneous or the inhomogeneous models namely, the models  $\{dx^2+x^3dx^1,dx^3+x^4dx^1,dx^4+x^5dx^1\}$ and  $\{dx^2+x^3dx^1,dx^3+x^4dx^1,dx^1+x^5dx^4\}~,$ and those systems for which $rank~S_1=2$ and $rank~S_2=0$ \textit{i.e.}, those for which $S_2=0,$ that admit the discrete model $\{dx^1+(x^3+x^4x^5)dx^4,dx^2+x^3dx^5,dx^3+x^4dx^5\}~.$ Finally there is the case $S_1=0$ where the model of class 5 is $\{dx^2,dx^4,dx^5+x^4dx^3+x^2dx^1\}$ and that of class 4 $\{dx^1,dx^3,dx^2+x^3dx^4\}~.$ Apart from these, there are the systems for which $S_2=0~,$ where $rank~S_1$ can be either 2 or 1 and that belong to continuous moduli. These shall be examined below.

\vspace{2 mm}
We firstly examine the systems \textit{S} of rank 3 for which $S_2=0$ and $rank~S_1=2~.$ These always admit representatives $\{\omega^1,\omega^2,\omega^3\}$ that satisfy the following structure equations:
\begin{equation*}
d\omega^1\equiv\omega^3\wedge\omega^4\hspace{5 mm}mod(\omega^1,\omega^2)~,
\end{equation*}
\begin{equation}
d\omega^2\equiv\omega^3\wedge\omega^5\hspace{5 mm}mod(\omega^1,\omega^2)~,
\end{equation}
\begin{equation*}
\hspace{7 mm}d\omega^3\equiv\omega^4\wedge\omega^5\hspace{5 mm}mod(\omega^1,\omega^2,\omega^3)~,
\end{equation*}
where the forms $\{\omega^4,\omega^5\}~,$ chosen adequately, complete the previous forms to a local basis in 5-space. We observe that the first two forms generate the co-variant (or first derived) system $S_1$ associated to \textit{S}. To determine those systems that are locally equivalent, one procedure is to determine the differential invariants of the pseudo-group that realizes such equivalences and, according to Cartan, the first step consists in exhibiting this pseudo-group, composed by all the local or infinitesimal automorphisms of such systems. Ultimately, the author puts forward a fundamental set of 27 higher order differential invariants that establishes the criterion for the equivalence. An outstanding situation occurs when the Cartan quartic form vanishes and, in this case, Cartan proves the following result.

\vspace{2 mm}
\newtheorem{rank5}[TheoremCounter]{Theorem}
\begin{rank5}
On a manifold of dimension 5, all the systems S of rank three with a vanishing Cartan quartic form are locally equivalent. The restriction, to a neighborhoods of the origin, of the infinitesimal automorphisms of such a system (or, equivalently, the set of all the germs, at the origin, of the infinitesimal automorphisms) is a simple algebra of dimension 14 isomorphic to the exceptional complex simple Lie algebra $g_2$ with real form $g_{2(2)}.$
\end{rank5}

\vspace{2 mm}
\noindent
The argumentation is extremely simple and entirely based on what Cartan calls \textit{the fundamental identity} namely, $d^2=0~$. Quite to the contrary, the calculations are extremely long to the point that the author himself does not even exhibit a single line of them and just mentions the conclusions. We should also point out that Cartan does not try to determine directly the infinitesimal transformations that carry one system into an equivalent one but examines initially the linear group, with variable (function type) entries (a matrix group), that transforms, via linear combinations, a given system represented by the forms $\{\omega^i\}$ into another system $\{\overline{\omega}^i\}~,$ both satisfying the structure equations (4). This last group has 7 parameters (dimension 7). There are other situations to be considered namely those when the Cartan quartic form does not vanish identically. Nevertheless, we shall not enter here into the details since the context becomes rather complicated and somewhat inconclusive.

\vspace{2 mm}
Let us now take a look at what Cartan does with the systems of two second order partial differential equations in one unknown function \textit{z} of two independent variables \textit{x,y} namely, a system of the form
\begin{equation}
r=\textbf{R}(x,y,z,p,q,t)~,\hspace{10 mm}s=\textbf{S}(x,y,z,p,q,t)~,
\end{equation}
where \textit{s} is the second order mixed derivative. This is the choice made by Cartan but we could as well interchange \textit{r} with \textit{t}. However Cartan's argument will break down if we interchange \textit{s} with \textit{t}. Such systems are also called second order Monge equations. In our present days tedious formalism, this can be rephrased as follows (\cite{Kumpera1972}).

\vspace{2 mm}
\noindent
We start by taking the fibration $\pi:\textbf{R}^3~\longrightarrow~\textbf{R}^2~,$ $(x,y,z)~\longmapsto~~(x,y)~,$ and consider a locally trivial sub-bundle $\mathcal{E}_2\subset J_2\pi~,$ where the latter is the bundle of 2-jets of local sections of $\pi~.$ The \textit{symbol} of the equation $\mathcal{E}_2$  is the kernel of the restricted projection $\rho_{1,2}:\mathcal{E}_2~\longrightarrow~J_1\pi$ (\textit{loc.cit.}). An equation is said to be \textit{involutive} when its symbol is involutive \textit{i.e.}, for any positive integer \textit{k}, the symbol of the $(k+1)-$st prolongation of the equation is the \textit{derived} space of the symbol of the $k-$th prolongation and if furthermore the initial symbol is 2-acyclic. An equation is said to be linear when its ambient space is a vector sub-bundle of a jet space relative to an initially given vector bundle. The Pfaffian system associated to a $k-$th order differential equation is the restriction of the $k-$th order contact Pfaffian system, defined on the space of $k-$jets, to the sub-manifold representing the ambient space of the equation. It is a Pfaffian system generated by the restrictions of the $k-$th order contact forms. In order to reproduce the system (2) considered by Cartan, it suffices to take the above considered sub-bundle $\mathcal{E}_2$ and require that the following two conditions be satisfied:

\vspace{3 mm}
\textit{(a)} The restricted projection $\rho_{1,2}:\mathcal{E}~\longrightarrow~J_1\pi$ is a submersion with a 1-dimensional fibre and

\vspace{2 mm}
\textit{(b)} The differential \textit{dt} never vanishes in restriction to the above 1-dimensional fibres.

\vspace{3 mm}
\noindent
We next write down the contact Pfaffian system associated to the equation $\mathcal{E}_2$ namely,
\begin{equation*}
\omega^1=dz-pdx-qdy~,
\end{equation*}
\begin{equation*}
\hspace{2 mm}\omega^2=dp-\textbf{R}dx-\textbf{S}dy~,
\end{equation*}
\begin{equation*}
\omega^3=dq-\textbf{S}dx-tdy~,
\end{equation*}

\vspace{3 mm}
\noindent
and are faced\footnote{The reader should observe the misprint on pg.122, \textit{l.}6, where, in the expression of $\omega_1,$ the term $ydy$ should read $qdy.$} with a rank 3 Pfaffian system on a six dimensional manifold $\mathcal{E}_2~.$ Let us now look for conditions under which the above Pfaffian system can be \textit{reduced} to a five dimensional manifold, the purpose being two-fold. Not only Cartan wants to take advantage of the knowledge of his rank 3 models on 5-spaces but, most important, this is the basic key of all the argumentation. In other terms, he wants to find conditions under which the differentials $d\omega^1,~d\omega^2$ and $d\omega^3$ can be expressed in terms of $\omega^1,~\omega^2,~\omega^3$ and two additional independent Pfaffian forms $\omega^4,~\omega^5$ that can be written as linear combinations of the forms $dx,dy,dp,dq$ and $dt~.$ A simple computation shows that

\vspace{2 mm}
\begin{equation}
d\omega^1\equiv 0 \hspace{5 mm}mod(\omega^1,\omega^2,\omega^3)~,\hspace{59 mm}
\end{equation}
\begin{equation*}
d\omega^2\equiv dx\wedge d\textbf{R}+dy\wedge d\textbf{S}\equiv dx\wedge(\frac{\partial\textbf{R}}{\partial t}dt+\frac{d\textbf{R}}{dt}dy)+dy\wedge(\frac{\partial\textbf{S}}{\partial t}dt+\frac{d\textbf{S}}{dx}dx)~,
\end{equation*}
\begin{equation*}
\hspace{5 mm}d\omega^3\equiv dx\wedge d\textbf{S}+dy\wedge dt\equiv dx\wedge(\frac{\partial\textbf{S}}{\partial t}dt+\frac{d\textbf{S}}{dy}dy)+dy\wedge dt~,
\end{equation*}

\vspace{2 mm}
\noindent
where $\frac{df}{dx}$ and $\frac{df}{dy}$ denote the \textit{total derivatives} in jet spaces (\textit{loc.cit.}). Furthermore, we can simplify the expressions (4) to

\vspace{2 mm}
\begin{equation}
d\omega^1\equiv 0,\hspace{61 mm}
\end{equation}
\begin{equation*}
d\omega^2\equiv\frac{\partial\textbf{R}}{\partial t}dx\wedge dt+\frac{\partial\textbf{S}}{\partial t}dy\wedge dt+(\frac{d\textbf{R}}{dy}-\frac{d\textbf{S}}{dx})dx\wedge dy\hspace{5 mm}mod(\omega^1,\omega^2,\omega^3)~, 
\end{equation*}
\begin{equation*}
d\omega^3\equiv\frac{\partial\textbf{S}}{\partial t}dx\wedge dt+dy\wedge dt+\frac{d\textbf{S}}{dy}dx\wedge dy~.
\end{equation*}

\vspace{2 mm}
\noindent
In order that the differentials of the generators $\omega^i$ of the contact Pfaffian system associated to the initial equations be expressible under the form (2), where the two forms $\omega^4$ and $\omega^5$ are linear combinations of $\{dx,dy,dt\}~,$ it is necessary and sufficient that the appropriate proportionality relations hold between the coefficients of the second members in (5) (\textit{cf.}, \cite{Cartan1910}, pg.123, formulas (3)) and we find precisely the conditions under which the system (3) is involutive. We can now set:
\begin{equation*}
\omega^4=dy+\frac{\partial\textbf{S}}{\partial t}dx~,
\end{equation*}
\begin{equation*}
\omega^5=dt+\frac{\partial\textbf{S}}{\partial y}dx~,
\end{equation*}
and obtain the congruences
\begin{equation*}
d\omega^1\equiv 0~,\hspace{12 mm}
\end{equation*}
\begin{equation*}
d\omega^2\equiv\frac{\partial\textbf{S}}{\partial t}\omega^4\wedge\omega^5\hspace{5 mm}mod(\omega^1,\omega^2,\omega^3)~,
\end{equation*}
\begin{equation*}
d\omega^3\equiv\omega^4\wedge\omega^5~,
\end{equation*}
the characteristic system of (3) being generated by $\{\omega^1,\omega^2,\omega^3,\omega^4,\omega^5\}~,$ hence the system (2) (or the equation $\mathcal{E}_2$) as well as the associated contact Pfaffian system \textit{S} both factor to a space in five variables (the space of the characteristic variables). Two cases are to be distinguished:

\vspace{3 mm}
\textit{(a)} The second order partial derivative $\frac{\partial^2\textbf{S}}{\partial t^2}$ vanishes identically.

\vspace{2 mm}
\noindent
In this case, we find the local (discrete) models for the second order contact Pfaffian system \textit{S} to be $\{dx^1,dx^2+x^3dx^4,dx^3+x^5dx^4\}$ when $S_2$ is integrable and the homogeneous flag system $\{dx^1+x^2dx^4,dx^2+x^3dx^4,dx^3+x^5dx^4\}$ otherwise. These two cases present themselves when the initially given system (1) is linear \textit{i.e.}, the two equations are linear functions of their arguments:
\begin{equation*}
s=\alpha t+\beta~,\hspace{10 mm} r=\alpha s+\gamma~,
\end{equation*}
where $\alpha,\beta$ and $\gamma$ are functions of the independent arguments $x,y,z,p,q~.$

\vspace{3 mm}
\textit{(b)} The second order partial derivative $\frac{\partial^2\textbf{S}}{\partial t^2}$ does not vanish identically.

\vspace{2 mm}
\noindent
In this case, the system of two second order involutive partial differential equations falls into the generic context where the congruences relative to the associated contact Pfaffian system are equal to (2). No local discrete models are available since continuous moduli parametrize such systems.

\vspace{2 mm}
We finally consider second order single partial differential equations given by
\begin{equation*}
F(x,y,z,p,q,r,s,t)=0
\end{equation*}
and assume that both partial derivatives $\frac{\partial F}{\partial x}$ and $\frac{\partial F}{\partial y}~,$ with respect to the independent variables, do not vanish. This situation can be transcribed in jet spaces by a locally trivial sub-bundle of $J_2\pi$ for which the restriction of $\rho_{1,2}$ is a submersion with a $2-$dimensional fibre. The discussion here is rather extensive but the main conclusions are the following. The associated second order contact Pfaffian system is of rank 2 and factors to a Pfaffian system in 5-space when the double Monge characteristics of the equation are integrable. In this situation, we can single out the well known Monge-Ampère equation (\cite{Goursat1891,Goursat1898,Goursat1905})
\begin{equation*}
(r+A)(t+C)-(s+B)^2=0~.
\end{equation*} 
The associated Pfaffian system is a flag system of rank and length two, hence an Engel flag, its \textit{covariant} system, of rank three, is therefore integrable (\cite{Kumpera1999,Kumpera2002}) and its pseudo-group of local automorphisms is the natural prolongation of $Aut(\mathcal{D})~.$ In other terms, this peudo-group of local automorphisms can be obtained by firet lifting and subsequently prolonging the algebra of all the infinitesimal transformations leaving invariant the Pfaffian equation $dx^2+x^3dx^1=0$ \textit{i.e.}, the Darboux system in the space of three variables. Cartan further shows that all the equations with such an associated Pfaffian system are locally equivalent, the above Monge-Ampère equation being the standard model, and that the algebra of all the infinitesimal second order contact automorphisms of these equations is (locally) isomorphic to the algebra of all the vector fields on $\textbf{R}^3~.$ 

Another noteworthy example is given by the second order Goursat equations (\cite{Goursat1898,Kumpera1999}) whose algebra of all the infinitesimal second order contact automorphisms is isomorphic to the exceptional complex simple Lie algebra $g_2$ with real form $g_{2(2)}.$ Cartan finally shows that the Goursat systems are equivalent to the previously discussed second order involutive systems of two equations. As for the first order Monge problem, we refer the reader to \cite{Kumpera1976}.

\vspace{2 mm}
One last word is due. In all the above discussion, we never mentioned neither Jordan nor Hölder and the reason is the following. Firstly and concerning the above discussed differential equations, the only possible Jordan-Hölder resolutions are the shortest ones namely, those consisting of the initially given defining groupoid of the pseudo-group and a maximal normal sub-groupoid. However these do not appear explicitly for the following reason. Cartan is forced to make extensive calculations, which unfortunately he never published, and was convinced that the best place to achieve these calculations was, in fact, the initially given space of all the variables and not in some quotient space. Consequently, he produced, with the help of adequately chosen additional linear forms $\varpi,$ \textit{clones} of the quotient objects and, in particular, he produced a projectable clone pseudo-group, defined on the whole space, that projects isomorphically onto the pseudo-group defined in the quotient (Sections IX-XIII, see in particular p.171, \textit{l.}12-14, p.180, \textit{l.}10-14, p.185, \textit{l.}1-5, p.181, \textit{l.}9-14). We should mention as well that, Cartan working with differential $1-$forms, can use the nice trick: \textit{The dual of a quotient is canonically isomorphic to the annihilator of the kernel}, which enables him to consider the differential forms, defining the quotient system, as also living on the total space.

\section{On the integration of some differential systems transitively invariant under simple groupoids}
Perhaps the simplest (though not \textit{simple}) example already considered by Lie concerns, using a slightly different terminology, those Pfaffian systems of rank 1 that are transitively invariant under 
a \textit{solvable} Lie group, \textit{integrable} in the terminology 
of Lie and Cartan (\cite{Lie1874},\cite{Lie1885},\cite{Kumpera1973}). 
In this case, the system $\mathcal{S}$  is integrable (its class being equal to 1), the integration being achieved by successive quadratures since the corresponding automorphic systems are invariant under 1-dimensional abel-ian Lie groups (\textit{cf.} also \cite{Cartan1921}).

\vspace{2 mm}
Let us next examine the integration of a differential system invariant under the action of a first order \textit{unimodular} groupoid (\cite{Lie1876},\cite{Lie1877},\cite{Cartan1938/39},\cite{Kumpera1965}). We consider a manifold \textit{M} together with a volume form $\Omega,$ two arbitrary unimodular vector fields $\xi$ and $\zeta$ (\textit{i.e.}, that preserve $\Omega$) and inquire, together with Lie, on the set of all local transformations $\varphi$ that send $\zeta$ into $\xi.$ To avoid at present unnecessary complications, we shall assume that these vector fields never vanish. Such transformations are the local solutions of the first order differential system
\begin{equation*}
\mathcal{S}=\{X\in\Pi_1M~:~X_*\zeta_{\alpha(X)}=\xi_{\beta(X)},~X^*\Omega=\Omega\}
\end{equation*}
that is invariant under the \textit{left} action (on the target) of the Lie pseudo-group $\Delta$ of all local unimodular transformations preserving $\xi~,$ this action being transitive. Equivalently, the equation is invariant under the first order groupoid $J_1\Delta$ equal to the set of all invertible $1-$jets \textit{X} such that $X\cdot Y\in J_1\Omega$ for every element $Y\in J_1\Omega~.$ Next, we denote by $\Phi$ the foliation defined by the integral curves of $\xi$ and by $N=M/\Phi$ the corresponding quotient space. We shall not worry whether this quotient has a compatible differentiable structure since we can always argue locally and we set $p:M\longrightarrow N$ for the associated quotient projection. Writing $\Pi_1(M,N)$ for the set of all $1-$jets of local surmersions, the projection $p$ prolongs to a projection $p_1:\Pi_1M\longrightarrow\Pi_1(M,N)$ and we define the differential system, of order 1 , $\mathcal{S}_1=p_1(\mathcal{S})~.$ Denoting by $\Delta_0$ the sub-pseudogroup of $\Delta$ composed by those transformations that preserve each leaf, we obtain a normal sub-pseudogroup and $\Delta$ factors to a quotient pseudogroup $\Delta_1$ that can be considered as the quotient $\Delta_1=\Delta/\Delta_0~.$

\vspace{2 mm}
We now consider the $(n-1)-$form, $n=dim~\!M~,$ $\tilde{\omega}=i(\xi)\Omega$ and, since $i(\xi)\tilde{\omega}=0~,$ this $(n-1)-$form factors to a volume form $\omega$ defined on the quotient manifold \textit{N}, it being clear that $Aut(\omega)=\Delta_1~.$
The differential system $\mathcal{S}_1$ is then transitively invariant under the action of the unimodular pseudo-group $\Delta_1$ and consequently can be integrated via the Jacobi last multiplier method. Choosing arbitrarily $n-2$ independent functions, the remaining one is obtained by a quadrature.

\vspace{2 mm}
To proceed with the integration of $\mathcal{S}~,$ we continue as follows. For any local solution $\sigma$ of $\mathcal{S}_1,$ we consider the inverse image of $im~\!\sigma$ in $\mathcal{S}~,$ with respect to the projection $p_1~,$ denote it by $\mathcal{S}_{\sigma}$ and thus obtain a differential system invariant under the transitive action of $\Delta_0~.$
The integration of this system is quite simple and elementary. In fact, it suffices to determine a function \textit{f} such that $\zeta f=1~,$ which only requires a quadrature along the integral curves of $\zeta~.$ As for the Jordan-Hölder resolution involved, we are just considering the two terms series $\Delta\supset\Delta_0.$ The system $\mathcal{S}_1$ is called, by Lie and Cartan, the \textit{automorphic} system associated to $\mathcal{S}$ and to the given Jordan-Hölder resolution. The system $\mathcal{S}_{\sigma}$ is called the \textit{resolvent} system associated to the solution $\sigma~.$

\vspace{2 mm}
Another transitive simple pseudo-group is the \textit{symplectic} pseudo-group consisting of all the local transformations that leave invariant an exterior differential $2-$form of maximal rank (equal to $dim~\!M$). We shall not dwell here into the symplectic domain but just refer the reader to \cite{Lie1876} and \cite{Lie1877} for some classical applications. Nowadays, several very interesting situations do occur in connection with \textit{Poisson Structures}.

\section{The Jordan-Hölder resolutions and the integration of partial differential equations}
In this section we simply transcribe, in the context of Lie groupoids, some of the concepts that can be found in \cite{Kumpera1980} and \cite{Kumpera1999} for Lie groups.

\vspace{2 mm}
A Lie algebroid is said to be \textit{simple} when its only sub-sheave in ideals is the trivial one. As is customary, an algebroid reduced only to its unit elements \textit{is not} considered to be simple. A Lie algebroid which is a direct sum of simple and non-abelian (at the isotropy level) sub-algebroids is called \textit{semi-simple}. In the  case that abelian factors are present, the algebroid is called \textit{reductive}. A \textit{Jordan-Hölder resolution or series}, for a Lie algebroid $\mathcal{S}~~ ,$ is a composition series \textit{i.e.}, a descending chain of Lie sub-algebroids, each term being an ideal in the preceding one, and for which the successive quotients are simple. Jordan-Hölder resolutions are precisely those composition series that cannot be properly \textit{refined} (\cite{Kumpera1980}). We could also define \textit{principal resolutions} (all the terms are ideals in the initially given algebroid) but shall not do so since it is rather exceptional to find differential systems that place in evidence such resolutions. It should however be remarked that any Jordan-Hölder resolution is always a finite chain and that any composition series can always be refined to a Jordan-Hölder series, omitting of course the repeated terms. Moreover, any two Jordan-Hölder series have the same length and are equivalent.

\vspace{2 mm}
A Lie groupoid is said to be \textit{simple} when the only normal sub-groupoids are the trivial ones. A groupoid reduced to its unit elements \textit{is not} considered to be simple. A Lie groupoid that is a direct product of simple and non-abelian (at the isotropy level) sub-groupoids is called \textit{semi-simple}. Otherwise, it is called \textit{reductive}. A \textit{Jordan-Hölder resolution or series} for a Lie groupoid $\Gamma$ is a composition series \textit{i.e.}, a descending chain of (topologically) closed Lie sub-groupoids, each term being normal in the preceding one, and for which the successive quotients are simple. The present definition does not agree with that given in \cite{Kumpera1980} for Lie groups where we could establish a one-to-one correspondence between finite and infinitesimal Jordan-Hölder resolutions as soon as we assumed that the sub-groups were connected. In the present context we could as well give a similar definition by considering, for example, a \textit{weak} Jordan-Hölder resolutions as being a composition series that cannot be properly \textit{refined} (\cite{Kumpera1980}). Nevertheless, it is safer not to fondle  within the devil's workshop since Lie's second fundamental theorem does not hold, in general, for groupoids and algebroids (\cite{Almeida1985}). In fact, we cannot establish a one-to-one correspondence between the finite and the infinitesimal resolutions even if we were to adopt the \textit{weak} context. Inasmuch, we shall not discuss \textit{principal resolutions} for Lie groupoids. It should however be observed that any Jordan-Hölder resolution is always a finite chain and that any composition series can always be refined to a Jordan-Hölder series, omitting of course the repeated terms. Moreover, any two Jordan-Hölder series have the same length and are locally equivalent in a neighborhood of the units.

\vspace{2 mm}
There do not exist general methods or techniques for constructing Jordan-Hölder resolutions but only specific particular procedures that can be devised for certain types of groupoids. The same situation is already present in Lie groups theory and, in this section, we shall try to extend to the present context some of the results and techniques that are described in \cite{Kumpera1980}. For that matter, we assume from now on that all the manifolds considered, including the systems of partial differential equations envisaged as sub-manifolds contained in some $k-th$ order jet space, are \textit{connected} and \textit{second countable} and we first mention the following almost self evident statement

\vspace{2 mm}
\newtheorem{jhe}[PropositionCounter]{Proposition}
\begin{jhe}
Let $\Gamma_e$ be the isotropy group of the Lie groupoid $\Gamma$ at a fixed unit element $e~.$ Then, there is a one-to-one correspondence between the normal sub-groups of $\Gamma_e$ and the isotropy bundles of normal Lie sub-groupoids $\Sigma$ of $\Gamma~.$
\end{jhe}

\vspace{2 mm}
A similar statement holds for Lie algebroids. We next observe that there are many possibilities for the normal sub-groupoids (resp. normal algebroids) that admit the above mentioned isotropy bundles. Nevertheless, we can make a privileged choice that shows to be the appropriate one for $\Sigma$ and consists in adding to the isotropy sub-bundle $\Sigma_M$ the set of all the elements that belong to $\Gamma$ and whose sources and targets are distinct. We shall call these sub-groupoids \textit{saturated}.

\vspace{2 mm}
\newtheorem{jhi}[PropositionCounter]{Proposition}
\begin{jhi}
There is a one-to-one correspondence between the Jordan-Hölder series in $\Gamma_e$ and those in $\Gamma$ with successively saturated normal sub-goupoids i.e., where at each step, $\Sigma_{i+1}$ is a saturated, normal and closed sub-groupoid of $\Sigma_i~.$ 
\end{jhi}

\vspace{2 mm}
\noindent
Observe that the normality and the closeness of the sub-groupoids can be read by the corresponding properties of the sub-groups belonging to the isotropy. 

\vspace{2 mm}
Much can be said about these resolutions but, since we shall be concerned with partial differential equations, we content ourselves by only inquiring in this subject for the Lie groupoids that are sub-groupoids of a jet space $\Pi_kM~,$ one of the first challenges being the determination of those couples of sub-groupoids, one of them being normal and closed in the other, for which the quotients are simple. Similarly, we have a parallel infinitesimal theory for Lie algebroids that represents a first step towards the finite theory and, in this case, we also have at our disposal an additional bonus namely, a linear structure.\nocite{Barros1969} The above mentioned challenge seems a very hard one but, surprisingly enough, it is not so since it has already been essentially solved by Élie Cartan more than a hundred years ago. Let us say a few words about this. By definition, a Lie pseudo-group $\Gamma$ of order \textit{k} (\textit{groupe de transformations fini et continu}, in the sense of Cartan) is the set of all the local transformation, operating on a manifold \textit{M} , that satisfy a $k-th$ order differential system contained in $\Pi_kM~.$ We further require that the elements of $\Gamma$ verify the "group conditions" where the operation is the usual composition of transformations or, eventually, of their restrictions. When this is the case, the defining differential system is a sub-groupoid of $\Pi_kM$ and we usually require that it also be Lie. In the \textit{Mémoire} \cite{Cartan1909}, Cartan classified all the classes of \textit{simple} real (analytic) pseudo-groups, both transitive and intransitive, hence it is now quite easy to guess which should be the simple sub-groupoids contained in any $\Pi_kM~.$ Moreover, in order to do so, Cartan was forced to examine thoroughly the respective defining equations \textit{i.e.}, the defining sub-groupoids. In order to determine the simple quotients, it will take a little additional effort and, of course, will entail a visible factorization procedure to be performed, at least locally, on the manifold \textit{M} and on the pseudo-groups. One last word should be mentioned. The pseudo-groups considered above are in fact assumed to be $k-th$ order transitive which means that the $k-$th order prolongations of their elements, when restricted to the $k-th$ order defining differential systems, operate transitively.

\vspace{2 mm}
We next discuss differential equations whose invariance groupoids display Jordan-Hölder resolutions and assume in advance all the necessary regularity conditions so as everything should run smoothly though, in later sections, we shall also consider singular data that are indispensable in this theory. It has been already shown by B. Malgrange (\cite{Malgrange1970},\cite{Malgrange1972}) that we cannot always avoid singularities. In the above mentioned case, the singularities are not too bad being all contained in a closed \textit{Zariski} subset. More recently, it has also been shown that singularities must be taken into consideration, this having been emphasized in the study of \textit{Flag Systems} (\cite{Kumpera2014}, see also \cite{Kumpera1982} and \cite{Kumpera1991}). \nocite{Brown2010} \nocite{Matsushima1954}

\vspace{2 mm}
Recalling the definitions stated in \cite{Kumpera1972} and given a fibration (surmersion) $\pi:P~\longrightarrow~M~,$ a $k-th$ order differential system or partial differential equation of order \textit{k} is a locally trivial sub-bundle $\mathcal{S}_k$ of the jet space $J_k\pi~.$ The groupoid of all the invertible $k-$jets on a manifold \textit{P} will be denoted by $\Pi_kP~.$ We recall that a local diffeomorphism $\varphi:U~\longrightarrow~U'$ of \textit{P} is said to be an \textit{automorphism} of the differential system $\mathcal{S}_k$ when the $k-th$ prolongation of $\varphi$ preserves the total space of the equation, this prolongation requiring of course some transversality hypotheses. The local automorphisms can be recognised by their \textit{clones}, namely the $k-th$ order (local) contact transformations that preserve the total space of the equation. We finally denote by $\widetilde{\mathcal{C}}_k$ the restriction of the contact Pfaffian system $\mathcal{C}_k~,$ defined on $J_k\pi,$ to the ambient space of the $k-th$ order equation $\mathcal{S}_k~.$ Let us denote by $Aut(\mathcal{S}_k)$ the pseudo-group composed by all the local diffeomorphism operating on \textit{P} whose $k-th$ order  prolongations preserve the equation $\mathcal{S}_k~.$ When this equation is in involution then $Aut(\mathcal{S}_k)$ is a Lie pseudo-group of order \textit{k} and, furthermore, its $k-th$ order prolongation restricted to the underlying manifold of $\mathcal{S}_k$ is a \textit{Cartan pseudo-group}, as described in the introduction, whose defining forms $\omega^i$ can be chosen among the local sections of $\widetilde{\mathcal{C}}_k~.$ No wonder that Élie Cartan never considered other that his legitimately own $(Cartan-)$ pseudo-groups.

\vspace{2 mm}
The process of applying the techniques of the Jordan-Hölder resolutions to the integration of differential systems is an iterative process where, at each step $\{i,i+1\},$ essentially the same argumentation will take place. Therefore, in order to exhibit the method, we need not even worry that the series be Jordan-Hölder since, in fact, it can be any composition series. \nocite{Brown2010} \nocite{Matsushima1954}

\vspace{2 mm}
Let us therefore take a sub-groupoid $\Gamma_k$ of $\Pi_kP~,$ that leaves invariant and operates transitively on a given $k-th$ order differential system $\mathcal{S}_k$ and let us simply take a normal and closed  sub-groupoid $\Sigma_k$ of $\Gamma_k~.$ The first groupoid operates on the total space \textit{P} of the fibration $\pi$ by sending a \textit{source} point onto the corresponding  \textit{target} point and we can start by taking the quotient manifold $P/\widetilde{\Sigma}_k=\widetilde{P}~i.e.,$ the set of intransitivity classes, modulo the action of the sub-groupoid $\widetilde{\Sigma}_k~,$ by all the elements of $\Sigma_k$ that preserve the $\pi-$fibres of \textit{P} \textit{i.e.}, that transform points belonging to an arbitrary fibre into points belonging to the same fibre. The groupoid $\Sigma_k$ operates, by the (conjugate) left action\footnote{care must be taken to assure the projectability, in \textit{M}, of an element belonging to $\Sigma_k~,$ transversality conditions with respect to $\pi$ being consequently  required, where after the action on the $k-$jets is achieved by simple conjugation.}, on $J_k\pi$ and leaves invariant, by assumption, the system $\mathcal{S}_k$ where after we can take the quotient $\mathcal{S}_k/\Sigma_k=\widetilde{\mathcal{S}}_k$ and, as Élie Cartan would candidly claim, \textit{un calcul qui n'offre aucune difficulté montre alors que} $\widetilde{\mathcal{S}}_k\subset J_k\tilde{\pi}~,$ where $\tilde{\pi}:\widetilde{E}~\longrightarrow~M$ is the obvious quotient projection. The differential system $\widetilde{\mathcal{S}}_k~,$ the \textit{automorphic} system, is then invariant by the quotient groupoid $\Gamma_k/\Sigma_k~,~\Sigma_k$ being normal and closed. Given any solution $\tilde{\sigma}$ of the differential system $\widetilde{\mathcal{S}}_k~,$ its inverse image in $\mathcal{S}_k~,$ via the quotient projection, becomes a differential system, the \textit{resolvent system}, invariant under the groupoid $\Sigma_k$ (\cite{Kumpera1965,Kumpera1973,Kumpera1975,
Kumpera1998,Kumpera1985}).

\vspace{2 mm}
We might now ask ourselves why should we consider Jordan-Hölder resolutions? The groupoid $\Pi_kM$ is, for any manifold \textit{M}, always simple since the pseudo-group of all the local transformations defined on a manifold \textit{M} is simple (Lie dixit et Cartan confirmatur). However, there is not much sense in considering differential equations invariant under this groupoid since its solutions would be all the local sections. The $k-th$ order groupoid that leaves invariant, up to order \textit{k}, a given volume form is also simple and the integration method, for a system invariant under such a groupoid (or under the corresponding volume preserving pseudo-group), is the well known (generalized) \textit{Jacobi last multiplier} method (Lie dixit, \cite{Lie1885,Lie1895}). However, we should point out that Lie actually considered only \textit{linear} transformations with constant coefficients that leaved invariant a fixed volume form. As for the $k-th$ order Darboux groupoid associated to the pseudo-group $Aut(\mathcal{D}_q),~k=2q+1~,$ it is a simple groupoid (Cartan dixit) and the integration method corresponding to differential systems invariant under this groupoid is the topic that has been discussed in the previous section concerning the beautiful \textit{mémoire} of Élie Cartan. It deals, in fact, with the specific case of second order involutive systems (Monge equations) where the associated contact Pfaffian system is a \textit{flag systems} of length two or three. In the next sections we shall give more attention to this subject and, in particular, shall also mention the integration method that corresponds to the \textit{symplectic} groupoid that is also simple (Lie quoque Cartan dixerunt). As for the resolvent system, unfortunately anything can pop up. The only benefits that can be expected here are lower dimensions and an eventual \textit{order} reduction. It turns out, as shown in the specific situation examined in \cite{Kumpera1973}, that the resolvent systems are simply very convenient auxiliary equations.\nocite{Lie1871}\nocite{Lie1872}\nocite{Lie1875}
\nocite{Lie1876}\nocite{Lie1877}\nocite{Lie1884}\nocite{Lie1888}

\vspace{2 mm}
We now examine some special situations where the Jordan-Hölder series can be constructed and, inasmuch, the integrability problem solved. Unfortunately, the general context is far from being not only resolved but simply investigated.

\vspace{3 mm}
\textit{(a) The abelian data.}

\vspace{3 mm}
\noindent
A groupoid $\Gamma$ is said to be abelian when the composition of any two elements commutes. On account of the composability requirements, this simply means that all he isotropy groups $\Gamma_e$ are abelian. We can now choose a Jordan-Hölder resolution for a specific isotropy group $\Gamma_{e_0}~,$ this resolution having $1-$dimensional quotient and such a resolution will generate a corresponding Jordan-Hölder resolution for $\Gamma~.$ The integration of the differential system proceeds by quadratures. In fact, one can argue as in \cite{Kumpera1999} by choosing a basis of invariant forms for the contact Pfaffian system associated to the equation $\mathcal{S}~.$

\vspace{3 mm}
\textit{(b) The solvable data.}

\vspace{3 mm}
\noindent
The abelian case extends \textit{ipso fato} to \textit{solvable} groupoids since these groupoids also admit Jordan-Hölder resolutions with $1-$dimensional isotropy quotients.

\vspace{3 mm}
\textit{(c) The nilpotent data.}

\vspace{3 mm}
\noindent
Nilpotency being a particular case of solvability, not much can be added. However, this condition carries with itself some outstanding properties in what concerns principal resolutions. It can be shown that, in this case, Jordan-Hölder and principal resolutions coincide. \nocite{Giaro1978}

\vspace{3 mm}
\textit{(d) The semi-simple and reductive data.}

\vspace{3 mm}
\noindent
Here again we can proceed as in \cite{Kumpera1999} in what concerns the several particular cases discussed there and where we shall replace the finite dimensional objects by their infinite dimensional counterparts (\cite{Cartan1893},\cite{Cartan1938/39} p.49). Needless to say that, in this extended context, we shall also find homotopy obstructions.

\section{The Lie and the Cartan Layouts}
In general, we consider a fibration $\pi:P~\longrightarrow~M$ and denote by $\mathcal{S}_k$ a differential system of order \textit{k} defined over \textit{P} namely, defined in the Jet space $J_k\pi$ and invariant under the action of a Lie pseudo-group $\Gamma$ defined on \textit{P}. As already mentioned earlier, it is unrealistic to try avoiding singularities but presently, and just for the sake of convenience, we shall assume that the differential system $\mathcal{S}_k$ is a locally trivial fibre sub-bundle of the Jet space $J_k\pi$ of all the $k-$jets of local sections of $\pi$ (\textit{cf.} \cite{Kumpera1972}). Further, we shall assume that the pseudo-group $\Gamma$ operates transitively, after $k-th$ order prolongation, on the given differential equation as well as on all its prolongations and say that it operates transitively at all the orders, this latter condition being also paraphrased by saying that all the defining groupoids of orders $k+h,~h\geq 0,$ of the above pseudo-group, operate transitively (to the right) on the given equation and on all its prolongations. The first step, in view of integrating an equation, consists in determining the characteristics of $\mathcal{S}_k$ and, thereafter, reduce this differential system to its characteristic variables, the pseudo-group $\Gamma$ inducing as well an action on the reduced system. The second step, thereafter, consist in associating to each equation the descending chain (composition series) formed by the iterated derived equations or derived Pfaffian systems (\textit{systèmes déduits} in the terminology of Cartan, \textit{cf.} \cite{Kumperaa1982}) and, with the help of this descending chains, determine the nature of the integration problem involving the given equation.

\vspace{2 mm}
It also seems useful to recall, in this section, a few standard facts concerning Jet spaces. Given a fibration $\pi:P~\longrightarrow~M$ as above and denoting by $J_k\pi$ the space of all the $k-$jets of local sections of $\pi~,$ we first observe that $dim~J_k\pi$ is odd for any $k\geq 1~.$ A local section $\tau~,$ with respect to the \textit{source} projection $\alpha_k~,$ is said to be holonomic whenever it is of the form $j_k\sigma$ where  $\sigma$ is a local section of $\pi$ \textit{i.e.}, when it is equal, at each point, to the $k-$jet of a given (fixed) section $\sigma~.$ A first order contact element on $J_k\pi$ is also said to be holonomic whenever it is the tangent space, at some point, to the image of a holonomic section. The $k-th$ order canonical contact system $\mathcal{C}_k$ on the Jet space $J_k\pi$ is the Pfaffian system of maximum rank that annihilates all the holonomic linear contact elements. Equivalently, it is the system that vanishes on all the holonomic sections. A set of independent generators being then provided by:
\begin{equation*}
\omega^{\lambda}=dy^{\lambda}-\sum~p^{\lambda}_i dx^i~,
\end{equation*}
\begin{equation*}
~\omega^{\lambda}_i=dp^{\lambda}_i-\sum~p^{\lambda}_{ij} dx^j~,
\end{equation*}

\hspace{40 mm}...................................

\begin{equation*}
\omega^{\lambda}_{i_1~\dots~i_{k-1}}=dp^{\lambda}_{i_1~\dots~i_{k-1}}-\sum~p^{\lambda}_{i_1~\dots~i_{k-1}j} dx^j~ ,
\end{equation*}

\vspace{3 mm}
\noindent
where $(x^1,~\dots~,x^m)$ is a local coordinate system given in \textit{M},
\begin{equation*}
(x^1,~\dots~,x^m, y^1,~\dots~,y^n)
\end{equation*}
a compatible system in \textit{P} , the $p'$s are the remaining coordinates in $J_k\pi$ and where the above differential forms are invariant under any permutation of the lower indices (the order of the successive differentiations being irrelevant). In what follows, we shall only be concerned with partial differential equations whose solutions are scalar functions depending on the variables $x^i$ and consequently the coordinates $y^i$ reduce to a single coordinate still denoted by \textit{y}. In sum, the total space \textit{P} of the fibration $\pi$ admits local coordinates of the form $(x^1,~\dots~,x^m,y).$ Furthermore, in order to comply with the standard indexation notations, we shall replace \textit{m} by \textit{n} hence, in definite, the local coordinates of \textit{P} do write $(x^1,~\dots~,x^n,y).$ Consequently, we shall rewrite the above mentioned differential forms as adapted to the present situation namely,
\begin{equation*}
\omega=dy-\sum~p_i dx^i~,\hspace{6 mm} 1\leq i\leq n~,
\end{equation*}
\begin{equation*}
\omega_i=dp_i-\sum~p_{ij} dx^j~,\hspace{22 mm}
\end{equation*}

\hspace{28 mm}...................................

\begin{equation*}
\omega_{i_1~\dots~i_{k-1}}=dp_{i_1~\dots~i_{k-1}}-\sum~p_{i_1~\dots~i_{k-1}j} dx^j~,\hspace{8 mm}
\end{equation*}

\vspace{3 mm}
\noindent
No doubt, there is much interest in studying differential equations where the unknowns are several scalar functions of the given \textit{n} variables but this would force us to write in excess. It is no coincidence that Cartan had exactly the same feeling and, in \cite{Cartan1911}, he developed a full theory, in all its generality, by only discussing it in the special case where $n=3~.$ Believe it or not, he wrote 91 exquisite pages on the subject. We shall refer to the above limited context as being \textit{the Cartan realm}\footnote{In view of helping the reader with Cartan's unorthodox notations, we explain the meaning of the expression $\partial\omega'_3/\partial(\omega_1\omega_{23})$ (\textit{e.g.}, pg.35,  \textit{l.}15,-2 and the following pages). Inasmuch as the \textit{total derivative} of a function, we can as well calculate the \textit{total derivative} of a differential form $\omega$ (\textit{cf.}\cite{Kumpera1999} and \cite{Kumpera2002}), the quotient meaning actually $\partial\omega'_3\hspace{2 mm}mod~\partial\omega_1\omega_{23}~,$ $\omega'$ meaning $d\omega$ and finally, $\omega_1\omega_2$ meaning $\omega_1\wedge\omega_2~.$}.

\vspace{2 mm}
We finally observe that the integration of the equation $\mathcal{S}_k$ is entirely equivalent to the search of the $n-$dimensional integral sub-manifolds of the restriction $\widetilde{\mathcal{C}}_k$ of $\mathcal{C}_k$ to the equation $\mathcal{S}_k$ that, moreover, are transversal to the $\alpha-$fibres. The above described constructions are compatible with the prolongation procedures. It should be remarked that $\mathcal{C}_k$ is not integrable (in the sense of Pfaffian systems) though it is \textit{in involution} in the sense of Cartan. Let us indicate the references \cite{Kuranishi1967} and \cite{Kumpera1972} concerning the (standard) notion of involutiveness for systems of partial differential equations and Pfaffian system. Nevertheless, perhaps the best references, in the  sense of pleasantness and usefulness providing its practical calculating aspects, the reader will find in \cite{Cartan1945} and \cite{Cartan1911}, pg.6, \textit{l.}-6, where Cartan provides the explicit values, necessarily maximal, for the dimensions of the successive algebraic symbols of the equation together with its prolongations. We also refer the reader to \cite{Kumperaa1982} for the notions of \textit{characteristics} and \textit{derived systems} together with their calculating algorithm.
\noindent
As a first step, we begin by examining systems that consist of a single equation which, nevertheless, provides all the flavor of what should be expected in the general case.

\section{Single second order involutive equations}
We begin by considering the $k-th$ order equation
\begin{equation*}
    F(x^i,y,p_i,p_{ij},~\dots)=0
\end{equation*}
and do not specify, for the time being, any class requirements, just the everywhere non-vanishing of the partial derivatives $\partial F/\partial x^i~.$ In a more formal setting, we can take a fibration $\pi:P~\longrightarrow~M~,$ with the previously mentioned dimensional limitations, and a sub-manifold $\mathcal{S}_k$ of $J_k\pi$ satisfying the following two properties:

\vspace{5 mm}
(a) $codim~\mathcal{S}_k=1$ and

\vspace{3 mm}
(b) $\alpha_k:\mathcal{S}_k~\longrightarrow~M$ is surmersive while $\beta_k:\mathcal{S}_k~\longrightarrow~P$ is just submersive, $\alpha_k$ and $\beta_k$ being the source and the target maps respectively. 

\vspace{5 mm}
\noindent
Loosely speaking, the characteristics of a system of partial differential equations are a finite family of independent functions such that the given system can be re-written in terms of these functions, the number of these functions being minimum. The formal aspects of the previous hint run as follows. The equation $\mathcal{S}_k$ being given, we consider its associated contact Pfaffian system $\widetilde{\mathcal{C}}_k~,$ restriction of $\mathcal{C}_k$ to $\mathcal{S}_k$ and, via \textit{inner products} (\textit{cf.}\cite{Kumpera1982}), define subsequently its \textit{characteristic system} $\overline{\mathcal{C}}_k~.$ This latter system is not necessarily a regular Pfaffian system though the following result, due to Cartan, holds (its proof being far from trivial).

\vspace{3 mm}
\newtheorem{cartan}[LemmaCounter]{Lemma}
\begin{cartan}
The characteristic system associated to an involutive system of partial differential equations is always regular.
\end{cartan}

\vspace{3 mm}
A simple though fancy example shows up when the rank of the characteristic system is maximum \textit{i.e.}, when $\overline{\mathcal{C}}_k = T^*\mathcal{S}_k~,$ this being the case for the \textit{Laplacian} or, more precisely, the \textit{Laplace equation}. We invite the reader to perform here \textit{all} the necessary calculations that are far from being few. Curiously enough, a simple glimpse at the Laplace equation will convince anybody that it \textit{cannot} be re-written with a lesser number of variables.

\vspace{2 mm}
Continuing with our discussion on characteristics, we claim that the characteristic system is always integrable (\textit{cf.}\cite{Kumpera1982}) and therefore determines an \textit{integral foliation}. When this system has maximum rank, the leaves reduce to points whereas, when the rank is null, the leaves become the connected components of the underlying manifold. Most often, a foliation does not factor differentiably, modulo its leaves, since non-trivial holonomy groups at some points might represent an obstruction. A criterion for the quotientability of a foliation can be found via Pradine's holonomy groupoid (\cite{Pradines2007}). Nevertheless, we can always proceed with a local factorization in an open neighborhood of any point. Proceeding thus, locally, not only we obtain a quotient manifold \textit{W} but, far more interesting, the contact system $\widetilde{\mathcal{C}}_k$ factors to a Pfaffian system $\widetilde{\mathcal{C}}_{k,W}$ that has null characteristics which means that it cannot be further reduced. Any local coordinate system in \textit{W} becomes then a characteristic system for the initially given equation $\mathcal{S}_k$ which, on its turn, can be reduced and re-written in terms of a lesser number of variables \textit{i.e.}, can be re-defined by a $k-th$ order equation in a lower dimensional Jet space. A typical example is shown in the section 3.

\vspace{2 mm}
Once the above reduction to the sole characteristic variables has been performed, the next step consists in trying to detect the appropriate integration method or procedure which shall be unveiled by the nature of the corresponding \textit{derived composition series} of the contact Pfaffian system, at $k-th$ order, associated to the equation together with those corresponding to lower orders. To each given configuration of these composition series some integration procedure should be devised and, of course, the invariance pseudo-group as well as the invariance groupoids for the given equation will become crucial, the Jordan-Hölder decompositions for these finite or infinitesimal objects providing considerable help in this task as was already shown and discussed in the section 3 and in \cite{Kumpera1999}. At present, we just mention that Cartan, in \cite{Cartan1911}, gave a full account for the differential equations in three independent variables. Not only he determined, in this realm, all the reductions via the characteristics as well as placed in evidence the different integration methods by exhibiting all the corresponding \textit{derived series} but went much further by indicating the explicit integration procedures. Though veiled in behind of his calculations mainly, as already mentioned before, due to the fact that it is rather inconvenient to compute in the quotients, one can in fact discern, in this work, some traces of Jordan-Hölder series. Let us just add a few more words concerning this remarkable \textit{Mémoire}.

\vspace{2 mm}
In the Cartan realm, the first order full contact system is given by the form $\omega=dy-p_i dx^i$ and, in second order, by the forms $\{\omega,\omega_i\},$ where $\omega_i=dp_i-p_{ij}dx^j~,$ $1\leq i,j\leq 3~.$ Furthermore, three possibilities for the derived sequence do arise. The first one shows up, exemplified by the Laplace equation, when the class of the system is null \textit{i.e.}, when the corresponding characteristic system is equal to the whole co-tangent bundle. It then results directly from the definitions that this property is equivalent to the vanishing of the first derived system of the associated second order contact system, thus the derived sequence reduces to $\widetilde{\mathcal{C}}_2\supset 0~.$ At the other end, exemplified for instance by the truncated Laplace equation $\partial^2 f/\partial x_1^2=0~,$ we find the descending derived chain $\widetilde{\mathcal{C}}_2\supset (\widetilde{\mathcal{C}}_2)_1 =(\widetilde{\mathcal{C}}_2)_2$ with $rank(\widetilde{\mathcal{C}}_2)_1 =3~,$ the second term being therefore integrable. As for the remaining two cases (\textit{cf.}\cite{Cartan1911}, p.18, \textit{l.}(-7)-(-4) and p.21, \textit{l.}2-5), we find for both the same type of derived chain namely, a chain of the form $\widetilde{\mathcal{C}}_2\supset(\widetilde{\mathcal{C}}_2)_1\supset 0$ where the rank of $(\widetilde{\mathcal{C}}_2)_1$ is equal to one. The difference between the two cases relies in the dimensions of the characteristics, two in the first case and one in the second though, of course, mainly relies in the different structures of the initial contact systems, these structures being most relevant in devising the corresponding integration methods. The above rather succinct discussion is very carefully examined by Cartan and we should point out that, when integration is actually under consideration, we must assure ourselves that all the data is \textit{real analytic} since there are examples, in the $C^{\infty}$ realm, of involutive partial differential equations that do not admit any solution at all. We should also mention that Cartan examines partial differential equations defined by 2, 3, 4 and 5 equations, more than 5 equations becoming redundant, and the richest situation is by far that involving two equations. In this case, a second family of characteristics - the double Monge characteristics - can be defined, the confrontation of the two rendering much useful additional information. The reader might well guess that the most interesting situation arises when the two characteristics coincide. In any case, the reading of Cartan's 91 pages is pure delight.

\section{Two and more second order involutive equations}
We start the discussion with systems formed by two equation or, in other words, with differential systems defined by sub-manifolds of co-dimension equal to 2 and shall find considerable resemblance with the corresponding situation in the space of five variables (section 3). In fact, two distinct situations do arise though the derived sequences are quite similar.

\vspace{2 mm}
Examining the congruences explicited in the Cartan realm (\cite{Cartan1911}, pg.29, \textit{l.}13-16 and pg.34, \textit{l.}2-5), we find for both derived sequences the same form namely, $\mathcal{S}\supset(\mathcal{S})_1\supset 0$ with $rank~(\mathcal{S})_1=1~.$ Moreover, $(\mathcal{S})_1$ is simply the system of rank 1 obtained by lifting, to second order jets, the first order contact system hence $(\mathcal{S})_1$ is Darboux since its first derived system is null. The distinction between the two cases has to do, mainly, with their characteristics though the natures of their first derived systems are quite different. In the first case, the standard characteristics are $1-$dimensional whereas, the so called \textit{double Monge characteristics} always do exist and are $2-$dimensional, the tangent spaces to the above distinct  characteristics being transversal. Cartan shows that such systems can always be integrated by the \textit{Monge method} (\textit{cf.} \cite{Cartan1910} and sect.3). The second case, considerably more intricate than the first one, is in fact separated in three distinct sub-cases depending on whether the system possesses double Monge characteristics or not, the standard characteristics being always $1-$dimensional. Cartan exhibits remarkably clarifying examples of these three possibilities and eventually proves that such systems can, notwithstanding, be integrated by means of ordinary differential equations.

\vspace{2 mm}
We finally mention that, for systems with the largest number of equations \textit{i.e.}, those for which the representative sub-manifolds, in Jet spaces, have the least admissible (and significant) dimensions, necessarily equal to $dim~P+1~,$ the associated contact Pfaffian systems are always integrable. 
It should nevertheless be mentioned that Cartan, as is customary, writes many pages of calculations and, worse, omits many more, to such an extent that in this \textit{Mémoire} it seems hopeless to try to unveil, in behind the curtains, any Jordan-Hölder scheme quite to the contrary of what becomes rather apparent in \cite{Cartan1910}.

\section{The Darboux Layout}
Perhaps the most interesting and inspiring integration method devised by Élie Cartan has to do with the Darboux contact structures (\cite{Cartan1910},\cite{Kumpera1999}) though, undoubtedly, the \textit{Mémoires} \cite{Cartan1911} and \cite{Cartan1914} are also extraordinary and bring invaluable information.

\vspace{2 mm}
Let us recall that on a manifold \textit{M} of odd dimension 2q+1, a \textit{Darboux contact structure} $\mathcal{D}_q$ is a Pfaffian system of rank 1 and of maximum Cartan class 2q+1 \footnote{The class of a rank 1 Pfaffian system is always odd independently of the underlying manifold's dimension.}. It is locally generated by a Pfaffian form
\begin{equation*}
\omega= dx^0+x^1~dx^2+...+x^{2q-1}~dx^{2q}~.
\end{equation*}
According to Lie (\cite{Lie1895}), the pseudo-group $\Gamma(\mathcal{D}_q)$ of all the local transformations of \textit{M} that preserve the above Darboux system is simple. We shall now have a closer look at the differential systems transitively invariant under this pseudo-group or, for that matter, under the corresponding defining groupoids that, of course, are Lie and show how such systems do also reduce to spaces with 2n+1 variables. In \cite{Cartan1910}, Cartan examines this phenomenon for $q=1,2$ and, in \cite{Cartan1911}, for $q=3~.$ The argument, in the general case, is quite similar though worthwhile to examine since it involves the characteristics of the associated contact Pfaffian systems on Jet spaces and, more precisely, on their restrictions to the given differential systems and to their prolongations.

\vspace{2 mm}
Let us then consider a fibration $\pi:P~\longrightarrow~M$ and denote by $\mathcal{S}_k$ a differential system of order \textit{k} defined over $P$ and invariant under the action of the simple pseudo-group $\Gamma(\mathcal{D}_q)~.$ We also assume that the fibre dimension of $\pi$ is equal to one and, in order to avoid (unduly) technical complications, that the differential system $\mathcal{S}_k$ is a locally trivial fibre sub-bundle of the Jet space $J_k\pi$ of all the $k-$jets of local sections of $\pi$ (\textit{cf.},\cite{Kumpera1972}). Finally, we assume that the given action is transitive at all the orders as already specified  previously.

\vspace{2 mm}
Since $dim~J_k\pi$ is odd and the \textit{co-dimension} of $\mathcal{S}_k$ is equal to 1 (we are still considering a single equation), it follows that the \textit{dimension} of $\mathcal{S}_k$ is even. The second order full contact Pfaffian system as well as the restricted one have both rank $n+1$. This restricted system factors into a system defined on an odd dimensional manifolds, one dimension less than that underlying the given equation, as soon as the double Monge characteristics, of the given equation, are integrable. Inasmuch as in the case of 5 variables, it can be shown that all such systems are locally equivalent, the Pfaffian systems  corresponding to the reductions becoming a Cartan homogeneous flag of length $n~.$ Furthermore, the pseudo-groups of all their local automorphisms are the appropriate prolongations of $Aut(\mathcal{D}_q),~2q+1=n+1~,$  that, on the infinitesimal level, can be identified with the algebra of all the vector fields on $\textbf{R}^{n+1}~.$ Finally, we claim that all these systems are locally equivalent to the generalized Monge-Ampère equation involving all the second order partial derivatives, this equation becoming a model.

\section{Pfaffian systems, definitions and notations.}
We shall now examine the integrability problem for systems $\mathcal{S}$ of partial differential equations in one unknown function, special emphasis being given to first order systems. The geometrical setting is provided by a Pfaffian system, the so-called  canonical contact structure, defined on a Grassmann bundle of contact elements of a certain order, our approach being twofold. We first examine the Cauchy characteristics of the Pfaffian system $\mathcal{P}=\mathcal{P(S)}$ associated to a given system of partial differential equations $\mathcal{S}$ and obtained by restricting the above mentioned canonical contact structure to the sub-manifold $\mathcal{S}$ that defines the equation. We prove that the integrability of $\mathcal{S}$ is equivalent to the regularity of the Cauchy characteristics allied to the appropriate dimension. Next, we construct characteristic vector fields of $\mathcal{C(S)}$ via the geometry of contact transformations and prove that integrability is equivalent to the following analytical condition: The Lagrange bracket of any two functions vanishing on $\mathcal{S}$ also vanishes on $\mathcal{S}~.$ Joining the two ends, we infer that the regularity of the characteristics together with the appropriate dimension are equivalent to the vanishing condition for the Lagrange bracket. Our results extend, to systems of partial differential equations, the classical theory of the Cauchy characteristics and the Lie construction of these characteristics for a single first order equation, in which case no compatibility conditions are required. In other terms, a single first order equation in one unknown function is always integrable whereas two or more equations are, in most case, incompatible hence non-integrable. Along the way, we show how to integrate an integrable system using a simple extension of the Cauchy method. \nocite{Lie1874} \nocite{Lie1877} \nocite{Lie1888}

\vspace{2 mm}
Let \textit{M} be a finite dimensional connected differentiable ($C^{\infty}$) manitold, $T=TM$ its tangent bundle and $T^*=T^*M$ the cotangent bundle. A regular Pfaffian system defined on \textit{M} is, by definition, a locally trivial vector sub-bundle $\mathcal{P}$ of $T^*.$ Its annihilator $\Sigma=\mathcal{P}^{\perp}$ is a sub-bundle of \textit{T}, $(T/\Sigma)^*\simeq\mathcal{P}$ and $T^*/\mathcal{P}\simeq\Sigma^*.$ Such a sub-bundle of \textit{T} is called a regular distribution, occasionally a field of first order contact elements, on the manifold \textit{M}. Given a vector bundle \textit{E} (or, for that matter, any bundle), we denote by $\Gamma(E)$ the module of all the global sections and by $\Gamma_{\ell}(E)$ the pre-sheaf composed by the local sections. A local automorphism of $\mathcal{P}$ (or of $\Sigma$) is a local diffeomorphism $\varphi$ of \textit{M} satisfying the property $\varphi^*\mathcal{P}=\mathcal{P},$ this property being equivalent, by duality, to $\varphi_*\Sigma=\Sigma.$ The set $Aut(\mathcal{P})$ of all the local automorphisms of $\mathcal{P}$ is a pseudo-group of transformations of order 1 though not always Lie since it might fail to be complete. An infinitesimal automorphism of $\mathcal{P}$ (or $\Sigma$) is a local vector field $\xi$ defined on \textit{M} and generating a local 1-parameter group of local automorphisms. It can be characterized by the condition $\theta(\xi)\Gamma_{\ell}(\mathcal{P})\subset\Gamma_{\ell}(\mathcal{P})$ or, equivalently, by $\theta(\xi)\Gamma_{\ell}(\Sigma)\subset\Gamma_{\ell}(\Sigma),$ where $\theta(\xi)$ is the Lie derivative along the vector field $\xi.$ The set $\mathcal{L}(\mathcal{P})$ of all the infinitesimal automorphisms of $\mathcal{P}$ is a not necessarily complete infinitesimal pseudo-algebra of order 1. The eventually singular distribution \textit{L} induced by $\mathcal{L}(\mathcal{P})$ \textit{i.e.}, defined by $L_x=\mathcal{L}(\mathcal{P})_x,~x\in M,$ satisfies the integrability criterion of Stefan and Sussmann, its maximal integral manifolds being the connected components of the orbits of $Aut(\mathcal{P}).$ We say that $\mathcal{P}$ is homogeneous when $Aut(\mathcal{P})$ operates transitively and infinitesimally homogeneous when $\mathcal{L}(\mathcal{P})$ is transitive \textit{i.e.}, when $L=TM.$ Since \textit{M} is assumed to be connected, infinitesimal homogeneity implies homogeneity, the converse statement being however inexact. Finally, given a second Paffian system $\mathcal{P}'$ on a manifold $M',$ we shall say that $(\mathcal{P},x)$ is locally equivalent to $(\mathcal{P}',x')$ when there exists a local diffeomorphism $\varphi:U\longrightarrow U',$ $U$ and $U'$ being open neighborhoods of \textit{x} and \textit{x'}, such that $\varphi(x)=x'$ and $\varphi^*\mathcal{P}'=\mathcal{P}.$ Two systems $\mathcal{P}$ and $\mathcal{P}'$ are said to be locally equivalent when the above mentioned property holds for any pair of points $(x.x')\in M\times M'.$ Occasionally, we shall consider distributions and Pfaffian systems that are not regular in the sense adopted in the beginning of this section and shall specify, in each case, the appropriate regularity  conditions.

\section{Darboux characteristics.}
In this section we only consider scalar differential forms defined on the manifold \textit{M} though some of the considerations apply as well to vector valued forms.

\vspace{2 mm}
Let $\omega$ be an exterior differential form defined on an open subset \textit{U} of \textit{M} and denote by $\omega_x$ the induced form on $T_xM.$ The \textit{annihilator} of $\omega$ at the point $x\in U$ is, by definition, the vector subspace $V_x\subset T_xM$ defined by
\begin{equation*}
V_x=\{v\in T_xM~|~i(v)\omega_x=0\}.
\end{equation*}
When $\omega$ is a 1-form then $V_x=ker~\omega_x.$ It is easy to prove that $\omega_x$ factors to an exterior form defined on $T_x/V_x$ (\textit{i.e.}, $\omega_x=q^*\tilde{\omega}_x$ where \textit{q} is the quotient map) and that $V_x$ is the largest subspace of $T_xM$ having this property. If the degree of $\omega$ is equal to \textit{r}, then $V_x^{\perp}\subset T^*M$ is generated by the family
\begin{equation*}
\{i(e_{j_1}\wedge~\cdots~\wedge e_{j_{r-1}})\omega\},
\end{equation*}
where $\{e_j\}$ is any set of generators for $T_xM.$ Moreover, $V_x^{\perp}$ is the smallest sub-space \textit{W} of $T_x^*M$ such that $\omega_x$ belongs to the sub-algebra of $\wedge T_x^*M$ generated by $\textbf{R}+W.$ \nocite{Cartan1911} \nocite{Cartan1914}

\vspace{2 mm}
A local automorphism of the form $\omega$ is a local diffeomorphism $\varphi$ of \textit{M} such that $\varphi^*\omega=\omega$ and an infinitesimal automorphism is a local vector field $\xi$ generating a local 1-parameter group of local automorphisms. It is characterized by the condition $\theta(\xi)\omega=0$ in terms of Lie derivatives. The infinitesimal automorphism $\xi$ is said to be \textit{characteristic} when $i(\xi)\omega=0$ which means that $\xi_x\in V_x$ for any point \textit{x} in the domain of $\omega.$ Characteristic infinitesimal automorphisms are therefore defined in terms of the equations
\begin{equation*}
i(\xi)\omega=\theta(\xi)\omega=0
\end{equation*}
or, equivalently, by
\begin{equation}
i(\xi)\omega=i(\xi)d\omega=0
\end{equation}
in view of the relation $\theta(\xi)=[i(\xi),d].$

\vspace{2 mm}
The sum and the constant scalar multiples of infinitesimal automorphisms are infinitesimal automorphisms since $\theta([\xi,\eta])=[\theta(\xi),\theta(\eta)]$ and the bracket of infinitesimal automorphisms is again an infinitesimal automorphism. Furthermore, the formula $i([\xi,\eta])=[i(\xi),\theta(\eta)]$ shows that the bracket of characteristic infinitesimal automorphism remains characteristic and, finally, the formulas $i(f\xi)=fi(\xi)$ and $\theta(f\xi)=f\theta(\xi)+df\wedge i(\xi)$ show that $f\xi$ is characteristic when $\xi$ is, whatever the function \textit{f}. In sum, the set $\mathcal{L}_c(\omega)$ of all the characteristic infinitesimal automorphisms of $\omega$ is closed under addition, scalar multiplication by an arbitrary function and under the Lie bracket.

\vspace{2 mm}
We denote by $\Delta\subset TM$ the eventually singular distribution on the manifold \textit{M} (or, rather, on the open set \textit{U}) induced by the family $\mathcal{L}_c(\omega)$ \textit{i.e.},
\begin{equation*}
\Delta_x=\{\xi_x~|~\xi\in\mathcal{L}_c(\omega)\}
\end{equation*}
and call it the characteristic distribution of $\omega.$ Its annihilator $\Delta^{\perp}\subset T^*M$ is the characteristic system of $\omega.$ The distribution $\Delta$ satisfies the integrability criterion of Stefan and Sussmann and its maximal integral manifolds are the Cauchy characteristics of $\omega.$ In view of the equation (6), the elements of $\mathcal{L}_c(\omega)$ are the local vector fields $\xi$ on \textit{M} taking values in $\Delta.$ When the dimension of $\Delta_x,~x\in U,$ is constant, we say that $\omega$ has regular characteristics in which case $\Delta$ is a regular and integrable distribution on \textit{M} (rather \textit{U}) and $\Delta^{\perp}$ is a regular Pfaffian system. The following result enables to reduce, locally, the form $\omega$ to an expression involving only the characteristic variables \textit{i.e.}, the first integrals of $\Delta.$ The proof is straightforward.  

\vspace{2 mm}
\newtheorem{mlagrange}[PropositionCounter]{Proposition}
\begin{mlagrange}
Let $\omega$ be an exterior differential form of degree r with regular characteristics and let $\{y_1,~\cdots~,y_s\}$ be a fundamental system of independent first integrals of $\Delta$ defined in a neighborhood of a given point $x\in M.$ Then $\omega$ has the local expression
\begin{equation*}
\omega=\sum~a_{i_1~\cdots~i_r}(y)~dy^{i_1}\wedge~\cdots~\wedge dy^{i_r}~,
\end{equation*}
where the sum is extended to all the sequences $1\leq i_1<~\cdots~<i_r\leq s$~.
\end{mlagrange}

\vspace{2 mm}
\newtheorem{lagra}[CorollaryCounter]{Corollary}
\begin{lagra}
Let $\omega$ be an exterior differential form with regular characteristics and let us assume that its degree is equal to the co-dimension of $\Delta.$ Then, for any $x\in M,$ there exist s first integrals $\{z_1,~\cdots~,z_s\}$ of $\Delta,$ defined in a neighborhood of x and such that $\omega=dz_1\wedge~\cdots~\wedge dz_s.$ These forms are therefore closed and locally decomposable.
\end{lagra}

\vspace{2 mm}
Based on the equation (6), we now examine a point-wise approach to the characteristics and define, for any $x\in U,$ the vector sub-space of $T_xM$
\begin{equation*}
\tilde{\Delta}_x=\{v\in T_xM~|~i(v)\omega_x=i(v)d\omega_x=0\}~.
\end{equation*}
In general, $\tilde{\Delta}=\cup~\tilde{\Delta}_x$ is a singular vector sub-bundle of \textit{TM} since the dimension of the fibres can vary and clearly $\Delta\subset\tilde{\Delta}.$ When $\omega$ has degree \textit{r}, then $\tilde{\Delta}$ is defined by the linear equations
\begin{equation*}
[i(w_1\wedge~\cdots~\wedge w_{r-1})\omega](v)=[i(w_1\wedge~\cdots~\wedge w_r)d\omega](v)=0~,
\end{equation*}
where $(w_1,~\cdots~,w_r)$ runs through all $r-$tuples of vectors in $T_xM$ and therefore its annihilator $\tilde{\Delta}^{\perp}\subset T^*M$ is generated by the linear forms
\begin{equation*}
\{i(w_1\wedge~\cdots~\wedge w_{r-1})\omega_x~,~i(w_1\wedge~\cdots~\wedge w_r)d\omega_x\}~.
\end{equation*}
It suffices, of course, to consider $r-$tuples $(e_{i_1},~\cdots~,e_{i_r}),$ where $\{e_i\}$ is a basis of $T_xM,$ such that $1\leq i_1<~\cdots~<i_r\leq dim~M.$ Since $dim~\Delta_x$ is lower semi-continuous and $dim~\tilde{\Delta}_x$ upper semi-continuous, it follows that $dim~\tilde{\Delta}_x,~x\in U,$ is constant if and only if $\Delta=\tilde{\Delta}.$ Furthermore, since the dimensions are integers and, in the present case bounded, they are locally constant on open dense subsets of \textit{U} where after the equality of $\Delta_x$ with $\tilde{\Delta}_x$ also takes place on an open dense subset. The integer $c_x=codim~\tilde{\Delta}_x$ is precisely the \textit{class} of $\omega$ at the point \textit{x}. When $\Delta=\tilde{\Delta},$ then $c_x$ is constant (assuming \textit{U} connected), this integer being the Cauchy class of the form $\omega.$ In what follows, we state the main results involving the notion of class.

\vspace{2 mm}
\newtheorem{range}[PropositionCounter]{Proposition}
\begin{range}
The class of a differential $1-$form $\omega$ is equal to $2p+1,$ at the point x, if and only if $\omega_x\wedge (d\omega_x)^p\neq 0$ and $(d\omega_x)^{p+1}=0.$ Under these conditions,
\begin{equation*}
\tilde{\Delta}_x(\omega\wedge d\omega^{\ell})=V_x(\omega)\cap V_x(d\omega)
\end{equation*}
for any $\ell=0,~\cdots~,p$. 
\end{range}

\vspace{2 mm}
\newtheorem{ranger}[PropositionCounter]{Proposition}
\begin{ranger}
The class of a differential $1-$form $\omega$ is equal to $2p,$ at the point x, if and only if  $(d\omega_x)^p\neq 0$ and $\omega_x\wedge(d\omega_x)^p=(d\omega_x)^{p+1}=0.$ This being the case,
\begin{equation*}
\tilde{\Delta}_x(\omega\wedge d\omega^{\ell})=\tilde{\Delta}_x(d\omega^q)=V_x(d\omega)
\end{equation*}
for any $\ell=0,~\cdots~,p-1$ and $q=1,~\cdots~,p$. If, moreover, $\omega_x\neq 0$ then the last condition $(d\omega_x)^{p+1}=0$ is a consequence of the first two. 
\end{ranger}

\vspace{2 mm}
\newtheorem{ranges}[PropositionCounter]{Proposition}
\begin{ranges}
The class of a closed differential $2-$form $\omega$ is, at every point, an even integer. This class is equal to $2p,$ at the point x, if and only if $\omega_x^p\neq 0$ and $\omega_x^{p+1}=0.$ Under these conditions,
\begin{equation*}
\tilde{\Delta}_x(\omega^{\ell})=V_x(\omega)
\end{equation*}
for any $\ell=1,~\cdots~,p$ and $\omega_x^p$ is decomposable. 
\end{ranges}

\vspace{2 mm}
\newtheorem{angels}[TheoremCounter]{Theorem (Darboux)}
\begin{angels}
Let $\omega$ be a differential $1-$form of constant class defined on the manifold M. When this class is equal to $2p+1$ then, at any point $x\in M,$ there exist local coordinates $(x^i)$ vanishing at x and such that $\omega$ has the local expression
\begin{equation*}
\omega=dx^1+x^2dx^3+~\cdots~+x^{2p}dx^{2p+1}.
\end{equation*}
When this class is equal to $2p$ and, moreover, when $\omega$ is everywhere non-singular, then it admits the local expression
\begin{equation*}
\omega=(1+x^1)dx^2+x^3dx^4+~\cdots~+x^{2p-1}dx^{2p}.
\end{equation*}
\end{angels}

\vspace{2 mm}
\newtheorem{angelical}[TheoremCounter]{Theorem (Darboux)}
\begin{angelical}
Let $\omega$ be a closed differential $2-$form of constant class equal to $2p$ defined on the manifold M. Then, at any point $x\in M,$ there exist local coordinates $(x^i),$ vanishing at x, such that $\omega$ has the local expression
\begin{equation*}
\omega=dx^1\wedge dx^2+dx^3\wedge dx^4+~\cdots~+dx^{2p-1}\wedge dx^{2p}.
\end{equation*}
\end{angelical}

\vspace{2 mm}
We finally observe, as a consequence of the Corollary 1, that non-singular $n-$forms ($n=dim~M$) always have the local expression $\omega=dx^1\wedge~\cdots~\wedge dx^n,$ non-singular closed $(n-1)-$forms the local expression $\omega=dx^1\wedge~\cdots~\wedge dx^{n-1}$ and non-singular $(n-1)-$forms satisfying $d\omega_x\neq 0$ the local expression
\begin{equation*}
\omega=(1+x^1)dx^2\wedge~\cdots~\wedge dx^n
\end{equation*}
in a neighborhood of \textit{x}. In all the above displayed expressions of the local \textit{canonical forms}, the coordinates $x^i$ figuring in these expressions are necessarily first integrals of the characteristic distribution $\Delta$ namely, they are the so-called \textit{characteristic functions}.

\section{Cartan characteristics.}
Let $\mathcal{P}$ be a Pfaffian system defined on the manifold \textit{M}. Inasmuch as above, an infinitesimal automorphism $\xi\in\mathcal{L}(\mathcal{P})$ is said to be characteristic whenever it is tangent to $\Sigma$ \textit{i.e.}, when $\xi\in\Gamma_{\ell}(\Sigma).$ The set $\mathcal{L}_c(\mathcal{P})=\mathcal{L}(\mathcal{P})\cap\Gamma_{\ell}(\Sigma)$ of all the characteristic vector fields of the system is an infinitesimal pseudo-algebra of order 1. Moreover, it is a pre-sheaf of modules with respect to the local $C^{\infty}-$functions on \textit{M} hence $f\xi$ is characteristic whenever $\xi$ is, \textit{f} being any local function on \textit{M}. We remark that Frobenius' integrability condition, in terms of brackets, amounts to say that $\mathcal{P}$ is integrable if and only if $\mathcal{L}_c(\mathcal{P})=\Gamma_{\ell}(\Sigma).$ The module $\mathcal{L}_c(\mathcal{P})$ induces an eventually singular integrable distribution $\Delta:x\in M\longrightarrow\Delta_x=\mathcal{L}_c(\mathcal{P})_x,$ called the characteristic distribution of $\mathcal{P},$ its annihilator $\Delta^{\perp}\subset T^*M$ being the characteristic system of $\mathcal{P}.$ The maximal integral manifolds of $\Delta$ are the Cartan characteristics of $\mathcal{P}.$ Again, $\mathcal{P}$ is integrable if and only if $\Delta=\Sigma,$ in which case the Cartan characteristics coincide with the integral leaves of $\Sigma.$ Since $\theta(\xi)=i(\xi)d+di(\xi),$ it follows that $\mathcal{L}_c(\mathcal{P})$ is the pre-sheaf of all the local vector fields $\xi$ on \textit{M} satisfying the equations
\begin{equation}
<\xi,\omega>=<\xi,i(\eta)d\omega>~,
\end{equation}
for all $\omega\in\Gamma_{\ell}(\mathcal{P})$ and $\eta\in\Gamma_{\ell}(\Sigma),$ and therefore $\mathcal{L}_c(\mathcal{P})$ is the set of all the local vector fields on \textit{M} taking values in $\Delta.$ When the dimension of $\Delta_x,~x\in M,$ is locally constant, we shall say that $\mathcal{P}$ has regular characteristics. The result that follows provides a method enabling us to reduce locally the Pfaffian system $\mathcal{P}$ to the sole characteristic variables \textit{i.e.}, to the first integrals of $\Delta.$ 

\vspace{2 mm}
\newtheorem{antenor}[TheoremCounter]{Theorem (Cartan)}
\begin{antenor}
Let $\mathcal{P}$ be a Pfaffian system with regular characteristics and $\{y^1,~\cdots~,y^s\}$ a fundamental system of independent first integrals of $\Delta$ defined in a neighborhood of a given point $x\in M.$ There exists then a local basis $\{\omega_1,~\cdots~,\omega_r\}$ of $\mathcal{P}$ defined in a neighborhood of x such that
\begin{equation*}
\omega^i=\sum~a^i_j(y^1,~\cdots~,y^s)dy^j\hspace{2 mm},\hspace{2 mm}1\leq i\leq r~.
\end{equation*}
\end{antenor}

\vspace{2 mm}
\noindent
Given a regular integrable distribution $\Sigma$ on the manifold \textit{M}, we say that the open set \textit{U} is \textit{simple} with respect to $\Sigma$ when the foliation associated to the restricted distribution $\Sigma|U$ admits a quotient, this meaning precisely that the quotient space $U/\Sigma,$ modulo the leaves of $\Sigma|U,$ admits a manifold structure for which the quotient map $\rho:U\longrightarrow U/\Sigma$ is a surmersion. The manifold structure of $U/\Sigma$ is then uniquely determined by the structure of \textit{U}.

\vspace{2 mm}
Expressed in geometrical terms, the previous theorem can be transcribed as follows enhancing the role of the finite automorphisms.

\vspace{2 mm}
\newtheorem{antena}[CorollaryCounter]{Corollary}
\begin{antena}
Let $\mathcal{P}$ be a Pfaffian system with regular characteristics, U an open set of M simple with respect to the characteristic distribution $\Delta$ and $\rho:U\longrightarrow U/\Delta$ the quotient surmersion modulo the leaves of $\Delta.$ Then there exists a unique Pfaffian system $\widetilde{\mathcal{P}},$ defined on $U/\Delta,$ such that $\rho^*\widetilde{\mathcal{P}}=\mathcal{P}|U.$ Moreover, $\mathcal{L}_c(\widetilde{\mathcal{P}})=0,$ the finite and infinitesimal automorphisms of $\mathcal{P}|U$ are locally $\rho-$projectable onto the corresponding automorphisms of $\widetilde{\mathcal{P}}$ (they are globally projectable when their domains have connected $\rho-$fibres) and, for any $x\in U,$ there exists an open neighborhood $U_x$ such that the projections $\mathcal{L}(\mathcal{P}|U_x)\longrightarrow\mathcal{L}(\widetilde{\mathcal{P}}|\rho(U_x))$ and $Aut(\mathcal{P}|U_x)\longrightarrow Aut(\widetilde{\mathcal{P}}|\rho(U_x)$ are surjective. Finally, any integral manifold $\widetilde{\mathcal{V}}$ of $\widetilde{\mathcal{P}}$ (not necessarily of maximal dimension) gives rise, by a lift up, to the integral manifold $\rho^{-1}(\widetilde{\mathcal{P}})$ of $\mathcal{P}.$
\end{antena}

\vspace{2 mm}
\noindent
Conversely, the following statement also holds.

\vspace{2 mm}
\newtheorem{anterior}[CorollaryCounter]{Corollary}
\begin{anterior}
Let $\mathcal{P}$ be a Pfaffian system defined on the manifold M and let us assume that, for every $x\in M,$ there exists an open neighborhood $U_x$, a surmersion $\rho_x:U_x\longrightarrow V_x$ onto a manifold $V_x$ and a Pfaffian system $\widetilde{\mathcal{P}}_x,$ defined on $V_x,$ having null characteristics $(i.e.,~\mathcal{L}_c(\widetilde{\mathcal{P}}_x)=0)$ and such that $\mathcal{P}|U_x=\rho^*_x(\widetilde{\mathcal{P}}_x).$ The dimension of $\Delta_x$ is then constant $(i.e., \mathcal{P}~has~regular~characteristics),$ each $U_x$ is simple with respect to $\Delta$ and the family $(\rho_x)$ is a foliated atlas for the characteristic foliation of $\mathcal{P}.$
\end{anterior}

\vspace{2 mm}
\noindent
With the notations of the Corollary 2, we shall say that each quotient manifold $U/\Delta$ is a space of characteristic variables for $\mathcal{P}$ and the system $\widetilde{\mathcal{P}}$ will be called a local reduction of $\mathcal{P}$ to its characteristic variables. Furthermore, a complete set of local first integrals of $\Delta$ will be called a complete set, or system, of characteristic variables for $\mathcal{P}.$

\vspace{2 mm}
It is useful to examine the point-wise approach to the calculation of the characteristics, this being actually the standard approach. Based on the relation (7), we define, for every point $x\in M,$ the vector sub-space of $T_xM,$
\begin{equation*}
\tilde{\Delta}_x=\{v\in\Sigma_x~|~i(v)d\omega\in\mathcal{P}_x,~\forall~\omega\in\Gamma_{\ell}(\mathcal{P})\}~.
\end{equation*}
Since this space is defined by the linear equations $<v,i(w)d\omega>=0,$ with $w\in\Sigma_x$ and $\omega\in\Gamma_{\ell}(\mathcal{P}),$ its annihilator $\tilde{\Delta}_x^{\perp}$ is the subspace of $T_x^*M$ generated by
\begin{equation}
\{\omega_x,~i(w)d\omega_x~|~w\in\Sigma_x,~\omega\in\Gamma_{\ell}(\mathcal{P})\}~.
\end{equation}
The family $\tilde{\Delta}=(\tilde{\Delta}_x)$ is, in general, a singular vector sub-bundle of $TM$ (the fibres might not have, locally, constant dimensions) even if $\Delta$ is regular. The integer $s_x=codim~\tilde{\Delta}_x$ is called the \textit{class} of $\mathcal{P}$ or $\Sigma$ at the point $x\in M.$ We observe that $\Delta\subset\tilde{\Delta}$ and that $\mathcal{L}_c(\mathcal{P})$ is the pre-sheaf of all the local vector fields $\xi$ on \textit{M} taking values in $\tilde{\Delta}.$ Furthermore and on account of a semi-continuity argument, the equality $\Delta=\tilde{\Delta}$ holds locally if and only if either $\Delta$ or $\tilde{\Delta}$ have, locally, constant dimension. When this equality occurs, the integer $s=s_x,$ whose value is locally independent of the point \textit{x}, is called the \textit{class} of the Pfaffian  system $\mathcal{P}$ or of the distribution $\Sigma.$ Usually, the Darboux and the Cartan Theorems are proved under more restrictive assumptions on the class constancy though this is quite irrelevant. In any case, both distributions $\Delta$ and $\tilde{\Delta}$ have locally a constant dimension in open and everywhere dense subsets of \textit{M}, the subset corresponding to the later being contained in that corresponding to the former. One last remark is due namely, that the preceding corollaries admit the following additional comment. If, in the corollaries, one of the systems $\mathcal{P}$ or $\widetilde{\mathcal{P}}$ satisfies $\Delta=\tilde{\Delta}$ the same will also hold for the other.

\section{Classical contact structures.}
In this section we consider a distribution $\Sigma$ of co-dimension 1 namely, a field of hyperplanes or equivalently, a Pfaffian system of rank 1 defined on the manifold \textit{M}, and recall a few basic facts. 

\vspace{2 mm}
\newtheorem{anteriores}[LemmaCounter]{Lemma}
\begin{anteriores}
Let $\mathcal{P}$ be a Pfaffian system of rank 1. Then, for any point $x\in M,$ the co-dimension of $\tilde{\Delta}_x$ is odd and consequently the class of the system $\mathcal{P}$ is always odd.
\end{anteriores}

\vspace{2 mm}
\noindent
When $\tilde{\Delta}_x=0$ (resp. $dim~\tilde{\Delta}_x=1$), then forcibly the dimension of \textit{M} is odd (resp. even). The results that follow establish the relationship between the Darboux class of the local generators $\omega$ and the Cartan class of $\mathcal{P}.$

\vspace{2 mm}
\newtheorem{datena}[PropositionCounter]{Proposition}
\begin{datena}
Let $\mathcal{P}$ be a Pfaffian system of rank 1 defined on the manifold M. Then the following properties are equivalent:

\vspace{3 mm}
(a) The class of $\mathcal{P}$ is maximum (i.e., $\tilde{\Delta}_x=0$).

\vspace{2 mm}
(b) $dim~M=2n+1$ and any local generator $\omega\in\Gamma_{\ell}(\mathcal{P})-0$ satisfies the condition $\omega\wedge(d\omega)^n\neq 0$ everywhere.

\vspace{2 mm}
(c) $dim~M=2n+1$ and, for any point $x\in M,$ there exists an element $\omega\in\Gamma_{\ell}(\mathcal{P})-0,$ defined in a neighborhood of x, such that $\omega_x\wedge(d\omega_x)^n\neq 0.$
\end{datena}

\vspace{2 mm}
\noindent
Pfaffian systems of rank 1 and maximum class ($=dim~M$) are, in view of the property (b), the classical contact structures on odd dimensional manifolds defined by a covering $(U_{\alpha},\omega_{\alpha})$ where the $(\omega_{\alpha})$ are differential $1-$forms of maximum Darboux class, defined on the open sets $U_{\alpha},$ and satisfy the compatibility condition $\omega_{\beta}=f_{\beta\alpha}\omega_{\alpha}$ on the overlap $U_{\alpha}\cap U_{\beta}.$

\vspace{2 mm}
\newtheorem{datafolha}[PropositionCounter]{Proposition}
\begin{datafolha}
Let $\mathcal{P}$ be a Pfaffian system of rank $1$ on the manifold M. Then, the following properties are equivalent:

\vspace{3 mm}
(a) The class of $\mathcal{P}$ is equal to $dim~M-1$ or, equivalently, $dim~\tilde{\Delta}_x=1.$

\vspace{2 mm}
(b) $dim~M=2n+2$ and any local generator $\omega\in\Gamma_{\ell}(\mathcal{P}-0)$ satisfies everywhere the condition  $\omega\wedge(d\omega)^n\neq 0.$

\vspace{2 mm}
(c) $dim~M=2n+2$ and, for any $x\in M,$ there exists $\omega\in\Gamma_{\ell}(\mathcal{P}-0),$ defined in a neighborhood of x, such that $\omega_x\wedge(d\omega_x)^n\neq 0.$  

\vspace{2 mm}
(d) $dim~M=2n+2$ and, for any $x\in M,$ there exists $\omega\in\Gamma_{\ell}(\mathcal{P}-0),$ defined in a neighborhood of x, such that $(d\omega_x)^{n+1}\neq 0.$

\vspace{2 mm}
(e) $dim~M=2n+2$ and, for any $x\in M,$ there exists $\omega\in\Gamma_{\ell}(\mathcal{P}-0),$ defined in a neighborhood of x, such that $\omega_x\wedge(d\omega_x)^n\neq 0$ and $(d\omega)^{n+1}= 0.$
\end{datafolha}

\vspace{2 mm}
\noindent
Pfaffian systems of rank 1 and class equal to $=dim~M-1$ are, in view of the property (d), the classical contact structures on even dimensional manifolds defined by a covering $(U_{\alpha},\omega_{\alpha})$ where each $\omega_{\alpha})$ is a Pfaffian form of maximum Darboux class defined on the open sets $U_{\alpha},$ these forms satisfying the compatibility condition $\omega_{\beta}=f_{\beta\alpha}\omega_{\alpha}$ on the overlap. Any such structure is, locally, the inverse image of an odd classical contact structure.

\vspace{2 mm}
We next indicate a general result, consequence of the Cartan Theorem, concerning Pfaffian systems of rank 1 and constant class. Initially, we state the following \nocite{Cartan1914}

\vspace{2 mm}
\newtheorem{daffiti}[LemmaCounter]{Lemma}
\begin{daffiti}
Given a Pfaffian system $\mathcal{P}$ on the manifold M, the following properties are equivalent:

\vspace{3 mm}
(a) The rank of $\mathcal{P}$ is equal to 1 and the class is constant (i.e., $\Delta=\tilde{\Delta}$).

\vspace{2 mm}
(b) There exists, for any $x\in M,$ an open neighborhood U, a submersion $\rho:U\longrightarrow\mathcal{U}$ and an odd classical contact structure $\widetilde{\mathcal{P}},$ on $\mathcal{U},$ such that $\mathcal{P}|U=\rho^*\widetilde{\mathcal{P}}.$

\vspace{2 mm}
When these equivalent conditions are satisfied, the equality: $class~\mathcal{P}=class~\widetilde{\mathcal{P}}$ holds. In particular, the even classical contact structures are precisely those Pfaffian structures obtained, locally, as inverse images via submersions $\rho$ with $1-$dimensional fibres.
\end{daffiti}

\vspace{2 mm}
\newtheorem{ratafolha}[TheoremCounter]{Theorem}
\begin{ratafolha}
Let $\mathcal{P}$ be a Pfaffian system of rank $1$ on the manifold M. Then, the following properties are equivalent:

\vspace{3 mm}
(a) The class of $\mathcal{P}$ is constant and equal to $2p+1.$

\vspace{2 mm}
(b) Every $\omega\in\Gamma_{\ell}(\mathcal{P}-0)$ satisfies, everywhere, the conditions $\omega\wedge(d\omega)^p\neq 0$ and $\omega\wedge(d\omega)^{p+1}=0$  

\vspace{2 mm}
(c) For any $x\in M,$ there exists a form $\omega\in\Gamma_{\ell}(\mathcal{P}-0),$ defined in a neighborhood of x, such that $\omega_x\wedge(d\omega_x)^p\neq 0$ and $\omega_x\wedge(d\omega_x)^{p+1}=0.$

\vspace{2 mm}
(d) For any $x\in M,$ there exists a form $\omega\in\Gamma_{\ell}(\mathcal{P}-0),$ defined in a neighborhood of x, such that $\omega_x\wedge(d\omega_x)^p\neq 0$ and $(d\omega)^{p+1}=0.$

\vspace{2 mm}
(e) If $dim~M>2p+1,$ then for any $x\in M,$ there exists $\omega\in\Gamma_{\ell}(\mathcal{P}-0),$ defined in a neighborhood of x, such that $(d\omega_x)^{p+1}\neq 0$ and $\omega_x\wedge(d\omega_x)^{p+1}= 0.$
\end{ratafolha}

\vspace{2 mm}
\newtheorem{antenado}[CorollaryCounter]{Corollary}
\begin{antenado}
Two Pfaffian structures $(\mathcal{P},M)$ and $(\mathcal{P}',M'),$ of rank $1$ and constant class, are locally equivalent if and only if $dim~M=dim~M'$ and $class~\mathcal{P}=class~\mathcal{P}'.$ In particular, any such Pfaffian structure is transitive.
\end{antenado}

\vspace{2 mm}
We shall now illustrate the previous discussion by examining two among the most relevant classical contact structures namely, the Grassmann or canonical contact structure on the manifold of hyperplanes and the Liouville structure on the cotangent bundle.

\vspace{3 mm}
\noindent
\textbf{Example 1.} \textit{The Grassmannian bundle of hyperplanes.}

\vspace{2 mm}
Let \textit{M} be a manifold of dimension equal to $n+1$ and let us denote by $G(M)$ the Grassmannian bundle of all the hyperplanes contained in the various tangent spaces to the manifold \textit{M}. We indicate the elements of $G(M)$ by \textit{H} or by $H_x$ whenever it is desirable to specify the base point $x\in M$ ($H_x\subset T_xM$) and by $\pi:TM\longrightarrow M$ the standard projection. Local coordinates can be assigned to $G(M)$ as follows: Let $H_a$ be an element of $G(M)$ and $(U,x^i),~1\leq i\leq n+1,$ any coordinate system on \textit{M} with domain \textit{U} containing the point \textit{a}. There exists a sequence $1\leq i_1<~\cdots~<i_n\leq n+1$ for which the restricted linear forms $dx^{i_1}|H_a,~\cdots~,dx^{i_n}|H_a$ constitute a basis of $H^*_a.$ Let us assume, for convenience that $i_j=j,$ let us denote by $y$ the remaining coordinate $x^{n+1}$ and let us consider the subset $\mathcal{U}\subset\pi^{-1}U$ of all the hyperplanes $H$ for which the family $\{dx^i|H\},~1\leq i\leq n,$ is free. Coordinates $(x^i,y,p_j)$ can now be assigned to any element $H\in\mathcal{U}$ by simply taking the coordinate $(x^i,y)$ of $\pi(H)$ and writing
\begin{equation}
dy|H=\sum~p_idx^i~, \hspace{3 mm}1\leq i\leq n~.
\end{equation}
If $N\subset M$ is a sub-manifold of co-dimension 1 and if $T_aN\in\mathcal{U}$ then, by the inverse function theorem, \textit{N} is locally, in a neighborhood of \textit{a}, the graph of a function $y=f(x^1,~\cdots~,x^n)$ and $p_i(T_aN)=(\partial f/\partial x^i)_b,~b=(x^i(a)).$ We observe that $G(M)$ can be identified with the real projective bundle $RP(T^*M)$ since any hyperplane $H\subset T_xM$ is annihilated by a line in $T^*_xM.$ A simple computation shows that the coordinates $p_i$ are, up to a sign, the inhomogeneous projective coordinates in the fibres of $RP(T^*M).$  More precisely, $p_i=-\tilde{p}_i/\tilde{p}_{n+1},~1\leq i\leq n,$ where the $(\tilde{p}_1,~\cdots~,\tilde{p}_{n+1})$ are homogeneous coordinates on the fibres of $RP(T^*M).$

\vspace{2 mm}
We next consider the line bundle $\Phi$ on $G(M)$ whose fibre, at the point \textit{H}, is equal to $T_{\pi(H)}M/H.$ This bundle is simply the quotient $\pi^{-1}(TM)/\mathcal{E},$ where $\mathcal{E}$ is the canonical vector bundle on $G(M)$ with fibre $\mathcal{E}_H=H.$ The fundamental form $\Omega$ on $G(M)$ with values in $\Phi$ can now be defined as the quotient of $T\pi\hspace{2 mm}mod\hspace{1 mm}\mathcal{E}$ namely, $\Omega(v)=T\pi(v)\hspace{2 mm}mod\hspace{1 mm}H,$ for any $v\in T_HG(M).$ Taking the coordinates $(\mathcal{U};x^i,y,p_j)$ introduced earlier and observing that the vector field $\partial/\partial y,$ on $G(M),$ factors to an everywhere non-vanishing local section $[\partial/\partial y]$ of the bundle $\Phi,$ we infer from the relation (9) that
\begin{equation*}
\Omega~|~\mathcal{U}=\omega[\partial/\partial y]~,
\end{equation*}
where
\begin{equation*}
\omega=dy-\sum~p_idx^i ~.
\end{equation*}
Let $\Sigma=ker~\Omega\subset TG(M)$ and $\mathcal{P}=\Sigma^{\perp}.$ The Pfaffian system $\mathcal{P}$ is called the \textit{canonical contact structure} on $G(M)$ and, since it is generated on the open set $\mathcal{U}$ by the contact form $\omega,$ we infer from the Proposition 5 that $\mathcal{P}$ is an odd classical contact structure. A simple computation of the characteristic system will show that the class of $\mathcal{P}$ is equal to $2n+1.$ In fact, it suffices to observe that $\Sigma_H=T\pi^{-1}(H)$ and thereafter determine a set of free generators for $\tilde{\Delta}_H^{\perp}.$

\vspace{2 mm}
It is interesting to observe that $\Sigma_H$ is composed of those tangent vectors \textit{v} to the point \textit{H} corresponding to the infinitesimal motions tangent to \textit{H} \textit{i.e.}, the infinitesimal motions of the base point \textit{x} that are contained in \textit{H}.

\vspace{3 mm}
\noindent
\textbf{Example 2.} \textit{The Liouville structure on the cotangent bundle.}

\vspace{2 mm}
Let us now turn our attention to the Liouville structure on $T^*M.$ For this, we define the distribution $\overline{\Sigma}$ on the cotangent bundle $T^*M$ by assigning to each point $\mu$ the hyperplane $\overline{\Sigma}_{\mu}=(T_{\mu}q)^{-1}(ker~\mu)$ where $q:T^*M\longrightarrow M$ is the standard projection. Denoting by $\lambda$ the Liouville form on $T^*M$ defined by $\lambda_{\mu}=q^*\mu$ \textit{i.e.}, $\lambda_{\mu}(v)=<(T_{\mu}q)v,\mu>,$ it follows that $\overline{\Sigma}_{\mu}=ker~\lambda_{\mu}$ hence the Pfaffian system $\overline{\mathcal{P}}=\overline{\Sigma}^{\perp},$ defined on the manifold $T^*M,$ is generated by the Liouville form $\lambda.$ The distribution $\overline{\Sigma}$ as well as the Pfaffian system $\overline{\mathcal{P}}$ are regular on $T^*M-0,$ where $0$ stands for the null section, the latter being called the Liouville structure on the co-tangent bundle. Taking local coordinates $(x^i,p_ i),$ the Liouville form assumes the expression $\lambda=\sum~p_idx^i$ hence its Darboux class is equal to $2n+2$ ($dim~M=n+1$) and consequently, according to the Proposition 5, $\overline{\mathcal{P}}$ is an even classical contact structure on $T^*M.$ It is nevertheless interesting to have a closer look at the characteristic system. A straightforward calculation, using the above specified coordinates, will show that the characteristic subspace $\tilde{\Delta}_{\mu}$ is contained in the tangent space, at the point $\mu,$ to the fibre $T^*_xM,~x=q(\mu),$ and, since we can identify canonically $T_{\mu}T^*_xM$ with $T^*_xM,$ it also follows that $\tilde{\Delta}_{\mu}$ is the $1-$dimensional subspace of $T_{\mu}T^*_xM$ generated by $\mu.$ The tangent co-vector $\mu$ belongs therefore to $\tilde{\Delta}_{\mu}.$

\vspace{2 mm}
We shall try, next, to retrace the above considerations within a more geometrical context. Let \textit{L} be the Liouville vector field on $T^*M$ namely, the infinitesimal generator of the $1-$parameter group of homothetical transformations $(t,\mu)~\mapsto~h_t(u)=e^t\mu.$ A simple calculation shows, under the above-mentioned identification, that \textit{L} is the vertical vector field (tangent to the fibres) defined by $L_{\mu}=\mu$ and, therefore, that  $L=\sum~p_i~\partial/\partial p_i$ in local coordinates. The previous discussion also shows that \textit{L} generates $\widetilde{\Delta}=\Delta$ on $T^*M-0,$ a result that follows directly using the homothetical transformations. Let us also observe that $\theta(L)\lambda=\lambda$ and, since $<L,\lambda>=0,$ that the Liouville vector field is characteristic, hence the characteristic algebra of $\overline{\mathcal{P}}$ is precisely equal to the set of all the multiples $\{fL\}$ where \textit{f} is an arbitrary function.

\vspace{2 mm}
The maximal integral manifolds of the characteristic distribution $\Delta$ are the rays, emanating from the origin, in the various fibres $T^*M-0,$ hence the characteristic foliation admits a quotient manifold $(T^*M-0)/\Delta$ which is diffeomorphic to a sphere bundle over \textit{M}. The Pfaffian system $\overline{\mathcal{P}}$ factors then to a system $\widetilde{\mathcal{P}}$ on the quotient manifold, the class of $\widetilde{\mathcal{P}}$ being maximum and equal to $dim~(T^*M-0)/\Delta=2n+1.$ The system $\widetilde{\mathcal{P}}$ is therefore an odd classical contact structure as predicted locally by the Lemma 4. 

\vspace{2 mm}
We can however do slightly better. From the commutativity of the diagram

\vspace{2 mm}
\begin{equation*}
T^*M~\longrightarrow~RP(T^*M)
\end{equation*}
\vspace{1 mm}
\begin{equation*}
\downarrow\hspace{25 mm}\downarrow~
\end{equation*}
\vspace{1 mm}
\begin{equation*}
M~\leftarrow~Id~\rightarrow~M
\end{equation*}

\vspace{2 mm}
\noindent
and the sole definitions of $\Sigma$ and $\overline{\Sigma},$ it follows that $\overline{\Sigma}$ projects onto $\Sigma$ via the quotient map though the Liouville form $\lambda$ cannot project onto a generator of $\mathcal{P}$ (we should not forget that $RP(T^*M)$ is not orientable). However, restricting our attention to the open subsets where $p_{n+1}\neq 0,$ we can replace the generator $\lambda$ of $\overline{\mathcal{P}}$ by the local generator
\begin{equation*}
(1/p_{n+1})\lambda=dx^{n+1}+\sum~(p_i/p_{n+1})dx^i~,
\end{equation*}
the later form projecting onto $RP(T^*M)$ and giving rise to the standard generator $\omega=dy-\sum~p_idx^i$ of $\mathcal{P}$ since $p_i/p_{n+1}=-\tilde{p}_i/\tilde{p}_{n+1}$ and $x^{n+1}=y,$ the right hand side \textit{tildet}  coordinates being those considered in the first example.

\vspace{2 mm}
The system $\widetilde{\mathcal{P}}$ obtained by reducing $\overline{\mathcal{P}}$ to its characteristic variables sits on the twofold covering space $(T^*M-0)/\Delta$ of $G(M).$ It is rather interesting that $\widetilde{\mathcal{P}}$ can further be reduced to $\mathcal{P}.$ As a consequence of the last Corollary, we can also state the

\vspace{2 mm}
\newtheorem{anthras}[CorollaryCounter]{Corollary}
\begin{anthras}
Any odd classical contact structure is locally equivalent to the Grassmann structure on the manifold of hyperplanes and any even structure is locally equivalent to the Liouville structure on a co-tangent bundle.
\end{anthras}

\section{Cauchy characteristics.}
We consider again the canonical contact structure $\mathcal{P}$ defined, on the Grassmannian bundle $G(M),$ by the fundamental form $\Omega~,$ the manifold \textit{M} being of dimension $n+1~.$ We also denote by $\omega$ the scalar representative of $\Omega$ in a coordinate patch $(\mathcal{U};x^i,y,P_i).$

\vspace{2 mm}
Given an $n-$dimensional sub-manifold \textit{N} of \textit{M}, we denote by $N_1$ the set of all the tangent spaces $T_xN$ with $x\in N$ (it being desirable to distinguish $N_1$ from TN). The set $N_1$ is an $n-$dimensional sub-manifold of $G(M)$ diffeomorphic to \textit{N} via the projection $\pi$ and is transversal to the fibres of $\pi$ \textit{i.e.}, $T_HN_1\cap T_H\pi^{-1}(x)=0$ for any $H=T_xN.$ Manifolds such as $N_1$ will be called \textit{holonomic} sub-manifolds of $G(M)~.$ The following lemma is easily verified in local coordinates using the form $\omega~.$

\vspace{2 mm}
\newtheorem{posteriores}[LemmaCounter]{Lemma}
\begin{posteriores}
  (i) The subspace $\Sigma_H=ker_H\Omega$ is generated by all the contact elements $T_HN_1$ where $N_1$ is an arbitrary holonomic sub-manifold containing H. Equivalently, the canonical contact system $\mathcal{P}$ is generated by all the local $1-$forms that vanish simultaneously on every holonomic sub-manifold.
  
\vspace{2 mm}
~(ii) The sub-space $\Sigma_H$ contains $T_H\pi^{-1}(x)$ where $x=\pi(H)~.$
 
\vspace{2 mm}
(iii) An $n-$dimensional sub-manifold $\mathcal{N}$ of $G(M)$ transversal to the fibres of $\pi$ is locally holonomic if and only if it is an integral manifold of $\mathcal{P}$ i.e., $\iota^*\Omega=0$ (or, equivalently, $\iota^*\omega=0$) where $\iota:\mathcal{N}\longrightarrow G(M)$ is the inclusion map.
\end{posteriores}

\vspace{2 mm}
\noindent
The property \textit{(i)} shows that $\mathcal{P}$ admits, at any point, $n-$dimensional integral manifolds that are transverse to the fibres of $\pi$ and the property \textit{(ii)} shows that each fibre $\pi^{-1}(x)$ is also an $n-$dimensional integral sub-manifold. There are however $n-$dimensional integral manifolds other than those mentioned above, for example the sub-manifold defined by the equations:
\begin{equation*}
x^1=~\cdots~=x^{n-1}=0,~y=x^n,~p_n=1~.
\end{equation*}

\vspace{2 mm}
\noindent
More generally, let \textit{N} be a $q-$dimensional sub-manifold of \textit{M} , $0\leq q\leq n~,$ and denote by $N_1$ the set of all hyperplanes $H\in G(M)$ that are tangent to N namely, those verifying $T_xN\subset H$ where $x=\pi(H)\in N~.$ Clearly, $N_1$ is an $n-$dimensional sub-manifold of $G(M),$ the restricted projection $\pi:N_1\longrightarrow N$ is a submersion and, since $\iota^*\Omega=0~,~\iota:N_1\longrightarrow G(M)~,$ we infer that $N_1$ is also an integral sub-manifold of $\mathcal{P}.$ In particular, when $q=n~,$ $N_1$ becomes a holonomic sub-manifold and, when $q=0~,$ it is a discrete collection of fibres $\pi^{-1}(x),~x\in N~.$ We now show that such sub-manifolds are, generically, all the $n-$dimensional integral manifolds of $\mathcal{P}.$

\vspace{2 mm}
\newtheorem{costeriores}[LemmaCounter]{Lemma}
\begin{costeriores}
Let $\mathcal{N}$ be an $n-$dimensional integral sub-manifold of $\mathcal{P}$ and let us assume that $\pi|\mathcal{N}$ has constant rank. Then $\mathcal{N}$ is, locally, a sub-manifold of the form $N_1$ with $dim~N=rank~\pi|\mathcal{N}~.$
\end{costeriores}

\vspace{2 mm}
\noindent
\textbf{Proof.} We can assume, without loss of generality, that $\pi(\mathcal{N})=N$ is a connected $q-$dimensional sub-manifold of \textit{M} and that $\pi:\mathcal{N}\longrightarrow N$ is a submersion. Further, we can also assume that $0<q<n$ since the remaining two cases correspond precisely to the isolated fibres $\pi^{-1}(x)$ and to the holonomic sub-manifolds. Taking local coordinates $(\mathcal{U};x^i,y,p_i)$ on $G(M)~,$ the condition that $\omega=dy-p_idx^i$ vanishes on $\mathcal{N}$ implies that $\iota^*dy$ is a linear combination of the restricted differentials $\iota^*dx^i$ where $\iota:\mathcal{N}\longrightarrow G(M)$ is the inclusion map. Moreover, since the projection \textit{N} is $q-$dimensional, we can choose, locally, $n-q$ functions among the $x^i~,$ say for convenience $\{x^1,~\cdots~,x^{n-q}\}~,$ such that the differentials $\iota^*dx^i=d(x^i|\mathcal{N})~,$ $1\leq i\leq n-q~,$ are linearly independent, the remaining differentials $\iota^*dx^j,$ $j>n-q,$ and $\iota^*dy$ being linear combinations of these. It then follows, always locally, that, on $\mathcal{N}~,$ $x^j|\mathcal{N}=f^j(x^1|\mathcal{N}~,~\cdots~,~x^{n-q}|\mathcal{N})$ and $y|\mathcal{N}=g(x^1|\mathcal{N},~\cdots~,x^{n-q}|\mathcal{N})$ hence, on \textit{M} , the equations
\begin{equation}
x^j-f^j(x^i)=0\hspace{5 mm}and\hspace{5 mm}y-g(x^i)=0
\end{equation}
define a $q-$dimensional sub-manifold \textit{N} . Finally, the vanishing of $\omega$ on $\mathcal{N}$ coupled with the local expressions of the functions $x^j|\mathcal{N}$ and $y|\mathcal{N}$ in terms of the $x^i|\mathcal{N}~,$ provide further $n-q$ independent equations
\begin{equation}
\partial g/\partial x^i-(p_i+\sum~p_j~\partial f^j/\partial x^i)=0~,\hspace{5 mm}j>n-q~,
\end{equation}
that, together with the equations (10), define locally the sub-manifold $\mathcal{N}$ (for each $x\in N,$ the equations (11) define the fibre $(\pi|\mathcal{N})^{-1}(x)~.$ It is now straightforward to verify that the equations (11) define precisely, in coordinates, the set of all the hyperplanes contained in $\mathcal{U}$ that are tangent to \textit{N} . The proof is complete.

\vspace{2 mm}
We also emphasize that the integral manifolds of $\mathcal{P}$ have at most the dimension \textit{n} . In fact, we shall  prove in the sequel that $(n+1)-$dimensional integral contact elements do not exist, such elements being those linear sub-spaces \textit{W} for which $\omega|W=d\omega|W=0~.$

\vspace{2 mm}
(a) If \textit{W} is transverse to the fibres of $\pi,$ then the linear forms $\{dx^i|W,dy|W\},$ $1\leq i\leq n~,$ are linearly independent hence $\omega|W\neq 0~.$

\vspace{2 mm}
(b) If $W\supset T_H\pi^{-1}(x)~,$ where \textit{H} is the base point of \textit{W} , then $v_j=(\partial/\partial p_j)_H\in W$ and consequently $i(v^j)d\omega=-dx^j$ vanishes on \textit{W} . Since $\omega|W=0~,$ it also follows that $dy$ vanishes on \textit{W} hence this sub-space has at most dimension \textit{n} . More generally, let $p=dim~(T_H\pi^{-1}(x)\cap W)~,~p>0~,$ and let us denote by $(v^1,~\cdots~,v^p)$ a basis of the intersection. Setting $v^j=\sum~a^j_i\partial/\partial p_i,$ the matrix $a^j_i$ has rank \textit{p} and therefore the linear forms $\mu^j=i(v^j)d\omega=-\sum~a^j_idx^i,~1\leq j\leq p,$ are independent and vanish on \textit{W} . We next consider the sub-space $V=T\pi(W)\subset T_xM.$ Since $dim~V=n+1-p$ and since $dy|V=\sum~p_idx^i|V~,~p_i=p_i(H)~,$ it follows that a basis of $V^*$ can be chosen among the forms $\alpha^i=dx^i|V~,$ say $\{\alpha^1,~\cdots~,\alpha^{n+1-p}\}~.$ Taking vectors $w_1,~\cdots~,w_{n+1-p}\in W$ that satisfy the relations $<w_k,dx^i>=\delta^i_k~,~1\leq i,k\leq n+1-p~,$ we can write
\begin{equation*}
w_k=\partial/\partial x^k+*\partial/\partial y+\sum~*\partial/\partial p_i~,
\end{equation*}
and consequently,
\begin{equation*}
\eta_k=i(w_k)d\omega=i(w_k)\sum~dx^{\ell}\wedge dp_{\ell}\equiv dp_k\hspace{3 mm}mod~(dx^1,~\cdots~,dx^n)~.
\end{equation*}
We thus obtain a family of $n+1$ linearly independent $1-$forms $\{\mu^j,\eta_k\}$ vanishing on \textit{W}, where-after $dim~W\leq n~.$

\vspace{3 mm}
\noindent
\textbf{Remark.} The above indicated proof is, in fact, an argument involving only the $2n-$dimensional vector space $\Sigma_H$ and the exterior $2-$form $d\omega|\Sigma_H$ of maximum rank equal to $2n~.$ Abstracting the specific context involving the Grassmannian space $G(M)$ as well as the Pfaffian system $\mathcal{P}~,$ it actually proves a well known result concerning Lagrangean linear sub-spaces: Let \textit{V} be a $2n-$dimensional real or complex vector space and $\overline{\omega}$ an exterior $2-$form of rank $2n$ defined on \textit{V} . Then, the largest dimension of a linear sub-space $W\subset V$ on which $\overline{\omega}$ vanishes is equal to \textit{n} . This algebraic result, alone, limits the dimension of the integral sub-manifolds of $\mathcal{P}~.$ 

\vspace{2 mm}
We now turn our attention to the Cauchy problem for partial differential equations and, initially, restrict our discussion to first order systems since these already exhibit all the nuances of the theory and, to a great extent, simplify the language. By definition, a first order partial differential equations in one unknown function, on the manifold \textit{M} , is a sub-manifold $\mathcal{S}$ of co-dimension $1$ contained in $G(M).$ To simplify the terminology we shall replace the full expression simply by \textit{differential equation} and observe that, in the present context, it is far more befitting to consider contact element rather than first order jets of local sections as is usually done. A solution of $\mathcal{S}$ is an $n-$dimensional sub-manifold \textit{N} of \textit{M} such that $N_1\subset\mathcal{S}.$ Taking coordinates $(\mathcal{U};x^i,y,p_i)~,$ the sub-manifold $\mathcal{S}$ can be described, locally, by means of an equation $F(x^i,y,p_i)=0$ and, if \textit{N} is the graph of a certain function $y=f(x^1,~\cdots~,x^n)~,$ this sub-manifold will be a solution of $\mathcal{S}$ if and only if the function \textit{f} is a solution of the partial differential equation
\begin{equation}
F(x^i,y,\partial y/\partial x^i)=0~.
\end{equation}
 We indicate by $\Omega|\mathcal{S}$ the restriction of $\Omega$ to the sub-manifold $\mathcal{S}$ \textit{i.e.}, the form $\iota^*\Omega~,$ where $\iota:\mathcal{S}~\hookrightarrow~G(M)$ is the inclusion, and by $\mathcal{P}|\mathcal{S}$ the Pfaffian system on the manifold $\mathcal{S}$ generated by $\Omega|\mathcal{P}~.$ It is, according to the general theory of partial differential equations, the canonical contact Pfaffian system restricted to the equation. We finally denote by $\Sigma|\mathcal{S}$ the annihilator of the aforementioned Pfaffian system, this annihilator being equal to the kernel of $\Omega|\mathcal{S}$ or, equivalently, the intersection $\Sigma\cap T\mathcal{S}~.$ Inasmuch, if $\omega$ denotes a local generator of $\mathcal{P}~,$ then its restriction $\iota^*\omega$ to $\mathcal{S}$ becomes a local generator of $\mathcal{P}|\mathcal{S}~.$

\vspace{2 mm}
To integrate the differential equation $\mathcal{S}~,$ in the sense given to this problem by Sophus Lie, consists in determining the $n-$dimensional integral sub-manifolds of the Pfaffian system $\mathcal{P}|\mathcal{S}$ and, of course, the integral manifolds corresponding to the solutions of the differential equation $\mathcal{S}$ will be precisely those sub-manifolds transverse to the fibres of $\pi~.$ In the section 15, we shall examine the integration procedure for such Pfaffian systems, as conceived by Sophus Lie, via the geometry of contact transformations. At present, we simply compute the class of $\mathcal{P}|\mathcal{S}~,$ the method being essentially the same as that employed for $\mathcal{P}~.$ In the paragraphs that follow, the symbol $*$ indicates a coefficient that needs not be written explicitly.

\vspace{2 mm}
We observe initially that $\Omega|\mathcal{S}$ cannot vanish on any open subset $\mathcal{U}$ of $\mathcal{S}$ for then, $\mathcal{U}$ would be a $2n-$dimensional integral manifold of $\mathcal{P}$ which is never the case. However, since $\Omega|\mathcal{S}$ can eventually vanish on nowhere dense closed subsets, we shall assume that $\mathcal{P}|\mathcal{S}$ is a regular Pfaffian system on the whole manifold $\mathcal{S}$ in the sense that $(\Omega|\mathcal{S})_H\neq 0$ for all $H\in\mathcal{S}$ (\textit{i.e.}, the form $\Omega|\mathcal{S}$ never vanishes) or, equivalently, that $T_H\mathcal{S}+ker~\Omega_H=T_HG(M)~.$ We consider next the following two cases:   

\vspace{2 mm}
(a) $T_H\pi^{-1}(x)\subset T_H\mathcal{S}~.$ In this case, $v^j=(\partial/\partial p_j)_H\in(\Sigma|\mathcal{S})_H,~1\leq j\leq n,$ and $dim~T\pi(T_H\mathcal{S})=n$ hence the rank of the linear forms $i(v^j)d\omega|\mathcal{S}=-dx^j|\mathcal{S}$ and $\omega|\mathcal{S}\equiv dy|\mathcal{S}~mod~(dx^i|\mathcal{S})$ is precisely equal to \textit{n}. Further, since $T_H\pi^{-1}(x)\subset(\Sigma|\mathcal{S})_H~,$ $dim~T\pi[(\Sigma|\mathcal{S})_H]=n-1$ and we can choose vectors $w_k\in(\Sigma|\mathcal{S})_H,~1\leq k\leq n-1~,$ as well as a sequence of indices $1\leq i_1~<~\cdots~<~i_{n-1}\leq n$ in such a way that $<w_k,dx^{i_{\ell}}>=\delta^{\ell}_k~.$ It follows that
\begin{equation*}
i(w_k)d\omega|\mathcal{S}\equiv (dp_{i_k}+*dp_{i_n})|\mathcal{S}~mod~(dx^i|\mathcal{S})~,
\end{equation*}
where $i_n$ is the remaining index, hence the set of $2n-1$ linearly independent $1-$forms 
\begin{equation*}
\{i(v^j)d\omega|\mathcal{S}~,~i(w_k)d\omega|\mathcal{S}~,~\omega|\mathcal{S}\} 
\end{equation*}
generates the characteristic system $\tilde{\Delta}_H^{\perp}$ of $\mathcal{P}|\mathcal{S}~.$

\vspace{2 mm}
(b) $T_H\pi^{-1}(x)\not\subset T_H\mathcal{S}~.$ In this case, $dim~(T_H\pi^{-1}(x)\cap T_H\mathcal{S})=n-1~,$ $T\pi(T_H\mathcal{S})=T_xM$ and $dim~T\pi((\Sigma|\mathcal{S})_H)=n~.$ Choosing a basis $\{v^1,~\cdots~,v^{n-1}\}$ of the intersection $T_H\pi^{-1}(x)\cap(\Sigma|\mathcal{S})_H~,$ we obtain $n-1$ linearly independent $1-$forms $\mu^j=i(v^j)d\omega|\mathcal{S}$ since the $dx^i|\mathcal{S}$ are linearly independent on $\mathcal{S}~.$ Moreover, we can choose \textit{n} vectors $w_k\in(\Sigma|\mathcal{S})_H$ such that $<w_k,dx^i>=\delta^i_k$ and therefore obtain further \textit{n} linear forms $\eta=i(w_k)d\omega|\mathcal{S}~.$ The dimension of the intersection being equal to $n-1~,$ the rank of the forms $dP_k|\mathcal{S}$ becomes also equal to $n-1~,$ hence the rank of the forms
\begin{equation*}
\{\omega|\mathcal{S}\equiv dy|\mathcal{S}~mod~(dx^i|\mathcal{S})~,~\mu^j\equiv 0~mod~(dx^i|\mathcal{S})~,\hspace{15 mm}
\end{equation*}
\begin{equation*}
\hspace{60 mm}\eta_k\equiv dp_k|\mathcal{S}~mod~(dx^i|\mathcal{S})\}    
\end{equation*}
is at least equal to $2n-1~.$ Since, by the Lemma 3, the rank of the characteristic system is always odd and since $dim~\mathcal{S}=2n~,$ we infer that this rank is constant and equal to $2n-1~.$ One can prove, by a straightforward though rather long calculation on the manifold $\mathcal{S}~,$ that the rank of the forms $\{\mu^j,\eta_k\}$ is, modulo $\omega|\mathcal{S},$ precisely equal to $2(n-1)$ thus obtaining the desired rank for the characteristic system. However, if we follow Cartan and compute on the space of all the variables, namely on $G(M),$ then it is quite easy to arrive at the desired conclusion and obtain, along the way, a very convenient representation of the characteristic system. Taking local coordinates $(\mathcal{U};x^i,y,p_i)$ on $G(M)$ and representing locally  $\mathcal{S}$  by the equation $F=0,$ we observe that a vector $w\in T_HG(M),$ written as
\begin{equation*}
w=\sum~u^i\partial/\partial x^i+a\partial/\partial y+\sum~v_i\partial/\partial p_i~,
\end{equation*}
belongs to $(\Sigma|\mathcal{S})_H$ if and only if

\vspace{2 mm}
~(i) $<w,\omega>=0~,$ hence $a=\sum~p_iu^i$ and

\vspace{2 mm}
(ii) $<w,dF-\omega>=0~,$ hence
\begin{equation}
\sum(\frac{\partial F}{\partial x^i}+p_i\frac{\partial F}{\partial y})u^i+\sum~\frac{\partial F}{\partial p_i}v_i=0~.
\end{equation}
Since, by the regularity assumption on $\mathcal{P}|\mathcal{S},$ the forms $\omega$ and $dF$ are independent, the linear relation above is non-trivial and, furthermore, it is the only relation imposed on the coefficients $\{u^i,v_i\}.$ We infer that the space of linear forms
\begin{equation*}
\{i(w)d\omega=\sum~(u^idp_i-v_idx^i):w\in(\Sigma|\mathcal{S})_H\}
\end{equation*}
has rank $2n-1$ hence the linear system
\begin{equation*}
\{\omega,i(w)d\omega:w\in (\Sigma|\mathcal{S})_H\}    
\end{equation*}
has rank $2n~,$ when considered on the manifold $G(M),$ and consequently its kernel $V_H$ has dimension equal to 1. Moreover, since the linear sub-space $(\Sigma|\mathcal{S})_H$ is odd-dimensional, the restricted $2-$form $d\omega|(\Sigma|\mathcal{S})_H$ has a non-trivial annihilator that is clearly equal to $V_H$ hence the characteristic space $(\tilde{\Delta})_H$ of $\mathcal{P}|\mathcal{S}$ is also equal to $V_H$ and the restricted system
\begin{equation*}
\{\omega|\mathcal{S},i(w)d\omega|\mathcal{S}:w\in(\Sigma|\mathcal{S})_H\}~,
\end{equation*}
on the manifold $G(M),$ has rank $2n-1$ as desired. The characteristic distribution $\tilde{\Delta}=\Delta$ being of dimension 1 , it amounts to a system of ordinary differential equations. The relations
\begin{equation*}
\sum~u^idp_i-v_idx^i=0\hspace{5 mm}and\hspace{5 mm}dy-\sum~p_idx^i=0
\end{equation*}
together with (13) provide the system:
\begin{equation}
\frac{dx^1}{\partial F/\partial p_1}=\frac{dx^2}{\partial F/\partial p_2}=~\cdots~\frac{dx^n}{\partial F/\partial p_n}=
\end{equation}
\begin{equation*}
=\frac{-dp_1}{\partial F/\partial x^1+p_1\partial F/\partial y}=~\cdots~=\frac{-dp_n}{\partial F/\partial x^n+p_n\partial F/\partial y}=\frac{dy}{\Sigma~p_i\partial F/\partial p_i}~,
\end{equation*}

\vspace{2 mm}
\noindent
on the manifold $G(M),$ corresponding to the distribution $H~\mapsto~V_H~,$ hence its restriction to $\mathcal{S}$ provides the Cartan characteristics of the Pfaffian system $\mathcal{P}|\mathcal{S}~.$ These characteristics coincide, in view of the equation (14), with the well known characteristic stripes of the partial differential equation (12) \textit{i.e.}, with the Cauchy characteristics of this equation. The system (14) applies as well to the case (a) since neither assumptions (a) nor (b) are relevant in Cartan's argument. It should however be observed that the case (b) is, in fact, the natural context for the partial differential equation (12).

\vspace{2 mm}
The previous results extend easily to systems of partial differential equations. By definition, such a system of order 1 and on the manifold \textit{M} (in one unknown function) is a sub-manifold $\mathcal{S}$ of $G(M).$ When $codim~\mathcal{S}=q~,$ this sub-manifold is locally defined by a system of \textit{q} independent equations $\{F_{\alpha}=0\},$ these being called the local equations of $\mathcal{S}.$ A solution of $\mathcal{S}$ is an $n-$dimensional sub-manifold \textit{N} of \textit{M} such that $N_1\subset\mathcal{S}$ hence, if we are concerned with studying integrability problems, the systems of co-dimension $q>n+1$ become irrelevant. We shall therefore assume that $q\leq n+1~,$ the integer \textit{q} being called the rank of the system $\mathcal{S}~.$ Taking local coordinates $(\mathcal{U};x^i,y,p_i),$ the single relation (12) is now replaced by
\begin{equation*}
F_{\alpha}(x^i,y,\partial y/\partial x^i)=0~,\hspace{5  mm}1\leq\alpha\leq q~.
\end{equation*}
Similarly, we introduce the restricted form $\Omega|\mathcal{S}~,$ the restricted Pfaffian system $\mathcal{P}|\mathcal{S}~,$ the restricted distribution $\Sigma|\mathcal{S}$ and the restricted generator $\omega|\mathcal{S}~.$ We also restrict our attention to those equations $\mathcal{S}$ providing regular Pfaffian systems $\mathcal{P}|\mathcal{S}$ in the sense that they are locally trivial sub-bundles of $T^*\mathcal{S}$ and call such equations regular. When $q\leq n~,$ the restricted form $\Omega|\mathcal{S}$ cannot vanish on any open subset $\mathcal{U}$ of $\mathcal{S}$ otherwise this open set would be an integral manifold of $\mathcal{P}$ of dimension greater than $n~.$ We infer that the regularity of $\mathcal{P}|\mathcal{S}~,$ when $q\leq n~,$ simply means that $(\Omega|\mathcal{S})_H\neq 0$ for any $H\in\mathcal{S}~.$ However, when $q=n+1~,$ $\mathcal{P}|\mathcal{S}$ can be the null system as is necessarily the case when it admits solutions.   

\vspace{2 mm}
We say that $\mathcal{S}$ is integrable when every point $H\in\mathcal{S}$ belongs to an $n-$dimensional integral manifold of the Pfaffian system $\mathcal{P}|\mathcal{S}$ and let us now prove that such a system has regular characteristics. If $q=n+1~,$ integrability of $\mathcal{S}$ implies that $\mathcal{P}|\mathcal{S}=0~,$ the null system admitting regular characteristics namely, the manifold  $\mathcal{S}$ being a characteristic. Moreover, $\mathcal{S}$ is locally a holonomic sub-manifold of $G(M).$

\vspace{2 mm}
\noindent
\newtheorem{catafolha}[TheoremCounter]{Theorem}
\begin{catafolha}
Let $\mathcal{S}$ be a system of regular and integrable first order partial differential equations or rank $q < n+1$ defined on the manifold M. Then the associated Pfaffian system $\mathcal{P}|\mathcal{S},$ on the manifold $\mathcal{S},$ has regular characteristics and its class is equal to $2(n-q)+1.$
\end{catafolha}

\vspace{2 mm}
\noindent
\textbf{Proof.} We follow Cartan's argument and begin by taking a local coordinate system $(\mathcal{U};x^i,y,p_i)$ on $G(M)$ as well as local equations $F_{\alpha}=0$ for $\mathcal{S}~.$ Since the forms $\{\omega,dF_{\alpha}\}$ are linearly independent, the equation (13) is replaced by a system of independent equations
\begin{equation*}
\sum(\frac{\partial F_{\alpha}}{\partial x^i}+p_i\frac{\partial F_{\alpha}}{\partial y})u^i+\sum~\frac{\partial F_{\alpha}}{\partial p_i}v_i=0~.
\end{equation*}
hence the kernel $V_H$ of the system $\{\omega,i(w)d\omega:w\in(\Sigma|\mathcal{S})_H\}$ has dimension $q~.$ Further, the characteristic sub-space $\widetilde{\Delta}_H$ of $\mathcal{P}|\mathcal{S}$ being necessarily contained in $V_H~,$ it suffices to show, as a consequence of the integrability, that $dim~\widetilde{\Delta}_H\geq q~.$ In fact, $\mathcal{S}$ being integrable, the linear space $(\Sigma|\mathcal{S})_H$ contains an $n-$dimensional sub-space \textit{W} on which $d\omega$ vanishes. Let us write $\mu=d\omega|(\Sigma|\mathcal{S})_H~.$ Then, since $dim~(\Sigma|\mathcal{S})_H=2n-q~,$ we can choose a basis $\{v_1,~\cdots~,v_{2n-q}\}$ of $(\Sigma|\mathcal{S})_H$ in such a way that $x_j\in W$ for $j\leq n$ and therefore the characteristic sub-space $\tilde{\Delta}_H$ becomes the kernel of the linear forms $\{i(v_k)\mu:1\leq k\leq 2n-q\}.$ Observing now that $\mu|W=0~,$ we infer that $i(v_j)\mu|W=0$ for $j\leq n$ and that the forms $\{i(v_k)\mu|W:~k>n\}$ vanish on a sub-space of dimension at last equal to $n-(2n-q-n)=q$ hence $dim~\tilde{\Delta}_H\geq q$ as desired. The class of $\mathcal{P}$ is therefore equal to $(2n+1-q)-q=2(n-q)+1~.$

\vspace{2 mm}
We terminate this section by discussing the integration procedure of the system $\mathcal{P}|\mathcal{S}$ via its characteristics. We integrate, in the first place, the characteristic distribution so as to obtain, locally, the characteristic variables of $\mathcal{P}|\mathcal{S}~.$ In later sections, we shall discuss an algorithm especially devised for this purpose. Next and according to the Corollary 2, we reduce locally the system $\mathcal{P}|\mathcal{S}$ to the space of is characteristic variables and observe, in view of the Lemma 4, that the quotient system $\widetilde{\mathcal{P}|\mathcal{S}}$ is an odd classical contact structure defined on a space of dimension $2(n-q)+1~.$ It will then suffice to take an $(n-q)-$dimensional integral manifold $\widetilde{\mathcal{V}}$ of $\widetilde{\mathcal{P}|\mathcal{S}}~,$ this being achieved by reducing this system to its canonical form $\tilde{\mu}=0$ with $\tilde{\mu}=dZ-\sum~P_idX^i,~1\leq i\leq n-q~,$ (recall the previous discussion concerning the integral manifolds of the canonical contact structure $\mathcal{P}$ on $G(M)$) and lift it to an integral manifold $\mathcal{V}=\rho^{-1}(\widetilde{\mathcal{V}})$ of $\mathcal{P}|\mathcal{S}~,$ the dimension of $\mathcal{V}$ being equal to $(n-q)+q=n~.$ From a practical point of view, the above discussion boils down to finding a generator $\mu$ of class $2(n-q)+1$ for the system $\mathcal{P}|\mathcal{S}~,$ expressed in canonical form. 

\vspace{2 mm}
The above integration procedure also provides, as an extra dividend, a result that puts in evidence the full significance of the characteristics.

\vspace{2 mm}
\noindent
\newtheorem{catarolha}[CorollaryCounter]{Corollary}
\begin{catarolha}
Let $\mathcal{S}$ be a regular system of partial differential equations of rank $q\leq n+1$ on the manifold M and assume that the associated contact Pfaffian system $\mathcal{P}|\mathcal{S}$ has regular characteristics of dimension equal to $inf\{q,n\}.$ Then $\mathcal{S}$ is integrable.
\end{catarolha}

\vspace{2 mm}
\noindent
Regularity of the characteristics, with the appropriate dimension, is therefore a necessary and sufficient condition for integrability. As for the participation of Jordan and Hölder in the above discussion, it runs as follows. In the first place, we take the pseudo-algebra $\mathcal{CH}$ of all the local infinitesimal automorphisms (\textit{i.e.}, characteristic vector fields) of $\mathcal{P}|\mathcal{S}$ and factor locally the manifold \textit{M} to the space of the sole characteristic variables where we shall find, after factoring inasmuch the Pfaffian system $\mathcal{P}|\mathcal{S},$ the previously appointed odd canonical contact structure. According to the next section, all we have to do is to figure out the behaviour of an odd contact structure of class $2(n-q)-1$ in relation to a similar structure of class $2(n-q)+1$ and shall perceive that the ending result will be a continued Jordan-Hölder sequence with $1-$dimensional quotients and, once terminated the initial integration process, the remaining processes boil down to ordinary differential equations namely, the integration of the previously mentioned vector fields in the part (b) of the above discussion.  

\section{Contact transformations.}
In this section we study more carefully the geometry of the contact structure defined on the Grassmannian bundle $G(M)$ by the fundamental form $\Omega~.$ Several results are generic statements concerning odd dimensional classical contact structures and could as well be developped in the general context, others are more specific to the structure on $G(M).$

\vspace{2 mm}
Let $\mathcal{P}$ be the canonical contact structure on $G(M).$ The finite automorphisms of $\mathcal{P}$ are called the (first order) contact transformations and the infinitesimal automorphisms, the infinitesimal contact transformations or Lie vector fields. In local coordinates $(\mathcal{U};x^i,y,p_i),$ contact transformations $\varphi$ are defined by the equation $\varphi^*\omega\equiv\omega~mod~\omega$ and Lie vector fields $\xi$ by $\vartheta(\xi)\omega\equiv 0~mod~\omega~.$

\vspace{2 mm}
Among all the contact transformations, the simplest are those that arise from the transformations on the base space \textit{M}. Any local diffeomorphism $\varphi:U\longrightarrow U'$ of \textit{M} operates, via $T\varphi~,$ on the linear contact elements of \textit{M} defining thereafter a local diffeomorphism $\wp\varphi:\pi^{-1}(U)~\longrightarrow~\pi^{-1}(U')$ on $G(M)$ by setting $\wp\varphi(H_x)=T\varphi(H_x).$ The first statement in the Lemma 5 assures that $\wp\varphi$ is a contact transformation. The assignment $\varphi~\mapsto~\wp\varphi$ is obviously functorial \textit{i.e.}, $\wp Id=Id~,$  $\wp(\psi\circ\varphi)=\wp\psi\circ\wp\varphi~,$ $\wp(\varphi^{-1})=(\wp\varphi)^{-1}$ and local \textit{i.e.}, $\wp(\cup \varphi_{\alpha})=\cup\wp\varphi_{\alpha}$ where $\cup\varphi_{\alpha}$ denotes the transformation obtained by patching together the transformations $\varphi_{\alpha}$ that agree on the overlaps In particular, we infer that a differentiable $1-$parameter family $(\varphi)_t$ of local diffeomophisms of \textit{M} (resp. a local $1-$parameter group) gives rise to a differentiable $1-$parameter family $(\wp\varphi)_t$ on $G(M)$ (resp. a local $1-$ parameter group) hence any local vector field $\xi,$ defined on an open set \textit{U} of \textit{M} and generating the local $1-$parameter group $(\varphi)_t~,$ gives rise to the Lie vector field $\wp\xi=\frac{d}{dt}\wp\varphi_t|_{t=0}$ defined on the open set $\pi^{-1}(U).$ The assignment $\xi~\mapsto~\wp\xi$ is, of course, local with respect to the patching of vector fields. Moreover, since $\wp\xi$ can as well be defined by the above derivative where $(\varphi)_t$ is replaced by any differentiable $1-$parameter family satisfying $\varphi_0=Id$ and $\frac{d}{dt}\varphi_t|_{t=0}=\xi~,$ and since
\begin{equation*}
\xi+\eta=\frac{d}{dt}\varphi_t\circ\psi_t|_{t=0}~,
\end{equation*}
\begin{equation*}
\alpha\xi=\frac{d}{dt}\varphi_{\alpha t}|_{t=0},\hspace{5 mm}\alpha\in\textbf{R}~,
\end{equation*}
\begin{equation*}
[\xi,\eta]=\frac{d}{dt}\varphi_{-s}\psi_{-s}\varphi_s\psi_s|_{t=0},\hspace{5 mm}s=|t|^{1/2}~,
\end{equation*}
we infer that the above assignment is functorial with respect to the Lie algebra sheaf structure \textit{i.e.}~,
\begin{equation*}
\wp(\xi+\eta)=\wp\xi+\wp\eta,\hspace{3 mm}\wp(\alpha\xi)=\alpha\wp\xi,\hspace{3 mm}\wp[\xi,\eta]=[\wp\xi,\wp\eta]~.
\end{equation*}
It should be observed, however, that $\wp(f\xi)$ is not equal to $f\wp\xi$ unless \textit{f} is a constant function and the general expression is, of course, $\wp(f\xi)=f\wp\xi+\delta(f)\xi$ where $\delta$ is a certain Lie derivative. The transformations $\wp\varphi$ and $\wp\xi$ are called the prolongations of $\varphi$ and $\xi$ respectively, to the bundle $G(M)$ and project via $\pi$ onto the initial transformations ($\pi\circ\wp\varphi=\varphi\circ\pi$ and $\pi_*\wp\xi=\xi$).

\vspace{2 mm}
\noindent
\newtheorem{catarro}[LemmaCounter]{Lemma}
\begin{catarro}
A necessary and sufficient condition that a contact transformation $\varphi$ (resp. a Lie vector field $\xi$) be a prolongation is that it projects, via $\pi,$ onto a local diffeomorphism $\varphi_0$ (resp. onto a local vector field $\xi_0$) of M. This being the case, then $\varphi=\wp\varphi_0$ (resp. $\xi=\wp\xi_0$) on the common domain. 
\end{catarro}

\vspace{2 mm}
\noindent
The proof will be omitted since it is obvious.

\vspace{2 mm}
We consider again, as in the section 13, the line bundle $\Phi$ on the manifold $G(M)$ and define, for each Lie vector field $\xi~,$ the section $h_{\xi}=<\xi,\Omega>$ of $\Phi$ whose domain coincides with that of $\xi~.$ The section $h_{\xi}$ is called the \textit{contact hamiltonian} of $\xi~.$ According to the notations of the section 5, we denote by $\mathcal{L}(\mathcal{P})$ the set of all the Lie vector fields, $\mathcal{P}$ being the canonical contact structure on the manifold $G(M),$ and by $\Gamma_{\ell}(\Phi)$ the set of all the local sections of $\Phi~.$

\vspace{2 mm}
\noindent
\newtheorem{caçamba}[LemmaCounter]{Lemma}
\begin{caçamba}
The mapping $\mathcal{L}(\mathcal{P})\longrightarrow\Gamma_{\ell}(\Phi),~\xi~\mapsto~h_{\xi},$ is bijective. 
\end{caçamba}

\vspace{2 mm}
\noindent
We first argue locally. Taking a coordinate system $(\mathcal{U};x^i,y,p_i)$ on $G(M),$ we consider the local section $[\partial/\partial y]=\partial/\partial y~mod~\mathcal{E}~,$ of the bundle $\Phi~,$ induced by the vector field $\partial/\partial y$ on $G(M).$ This section is everywhere non-null and therefore trivializes the bundle $\Phi$ on the open set $\mathcal{U}~.$ Since $\Omega|\mathcal{U}=\omega[\partial/\partial y],$ the hamiltonian $h_{\xi}=<\xi,\omega>[\partial/\partial y]$ can be replaced by the scalar function $H_{\xi}=<\xi,\omega>$ called the hamiltonian of $\xi$ relative to $\omega~.$ Furthermore, since $\xi$ is an infinitesimal automorphism of $\mathcal{P}~,$ then $\theta(\xi)\omega=i(\xi)d\omega+di(\xi)\omega=\lambda\omega$ hence $(i(\xi)d\omega+dH_{\xi})\wedge\omega=0~.$ A straightforward computation will also show that, conversely, given any differentiable function \textit{f} , the equations
\begin{equation}
i(\xi)\omega=f\hspace{5 mm}and\hspace{5 mm}(df+i(\xi)d\omega)\wedge\omega=0
\end{equation}
have a unique solution $\xi~,$ this solution being, of course, a Lie vector field since the second equation is simply $[\theta(\xi)\omega]\wedge\omega=0~.$ More precisely, the coordinate expression of $\xi$ is given by
\begin{equation}
\xi=-\sum~\frac{\partial f}{\partial p_i}\frac{\partial}{\partial x^i}+(f-\sum~\frac{\partial f}{\partial p_i}p_i)\frac{\partial}{\partial y}+\sum~(\frac{\partial f}{\partial x^i}+p_i\frac{\partial f}{\partial y})\frac{\partial}{\partial p_i}~. 
\end{equation}
The injectivity of the mapping in the Lemma 8 is now obvious since, assuming that that $h_{\xi}=0~,$ the scalar hamiltonian $H_{\xi}~,$ efined on any coordinate patch, also must vanish and therefore $\xi=0~.$ We observe that the injectivity follows inasmuch from the fact that the system $\mathcal{P}$ has null characteristics. If $f=<\xi,\omega>=0~,$ then $\xi$ is a characteristic vector field of $\mathcal{P}$ and consequently $\xi=0~.$ In order to prove the surjectivity, we take a local section $\sigma$ of $\Phi$ defined on an open set \textit{U}, a coordinate system $(\mathcal{U};x^i,y,p_i),$ trivialize the bundle $\Phi|\mathcal{U}$ by means of the section $[\partial/\partial y]$ and consider the unique solution $\xi$ of the equations (15), where \textit{f} and $\omega$ are determined respectively by $\sigma|\mathcal{U}=f[\partial/\partial y]$ and $\Omega|\mathcal{U}=\omega[\partial/\partial y]$ restricted to the open set $\mathcal{U}\cap U.$ We claim that two such solutions $\xi$ and $\overline{\xi}~,$ defined via two coordinate systems $(\mathcal{U};x,y,p)$ and $(\overline{\mathcal{U}};\overline{x},\overline{y},\overline{p}),$ do agree on the overlap defining, thereafter, a unique Lie vector field, on the open set \textit{U}, whose hamiltonian is the given section $\sigma~.$ We leave the details of the proof to the reader since it only involves straightforward calculations.

\vspace{2 mm}
The Lie bracket defined in $\mathcal{L}(\mathcal{P})$ can now be transported, via the isomorphism of the Lemma 8, to a Lie algebra bracket on the space $\Gamma_{\ell}(\Phi)$ namely, the Lagrange bracket for contact hamiltonians. Given a coordinate system $(\mathcal{U};x,y,p),$ hamiltonians can be replaced by differentiable functions and the Lagrange bracket becomes a Lie algebra bracket on functions. If $\xi$ is a Lie vector field and $f=<\xi,\omega>$ the corresponding hamiltonian, then it is easy to verify that
\begin{equation*}
[f,g]=\xi g-g\partial f/\partial y~.
\end{equation*}
We next observe that
\begin{equation*}
<\xi,df>=<\xi,df+i(\xi)d\omega>=<\xi,\lambda\omega>=\lambda f~,
\end{equation*}
hence $<\xi,df>_H=0$ whenever $f(H)=0$ and of course, under these conditions, $<\xi,\omega>_H=0~.$ Expressed in other terms, the Lie vector field $\xi$ with hamiltonian $h\in\Gamma_{\ell}(\Phi)$ is tangent to the variety $h^{-1}(0)$ and, at these points, takes values in the contact structure $\Sigma~.$ We observe that $h^{-1}(0)$ needs not be a regular sub-manifold of $G(M)$ since the function \textit{f} can admit singularities. Nevertheless, for any $H\in h^{-1}(0),$ the curve $exp~t\xi(H)$ is entirely contained in $h^{-1}(0)$ since the relation $<\xi,df>=\lambda f$ is equivalent to $\xi f=\lambda f$ and therefore equivalent, when $H\in G(M)$ is fixed, to the linear homogeneous equation
\begin{equation*}
\frac{d}{dt}f(exp~t\xi (H))=\lambda (exp~t\xi (H))f(exp~t\xi (H))~,
\end{equation*}
in the unknown function $f(exp~t\xi (H)).$ If $H\in h^{-1}(0),$ then $f(H)=0=f(exp~0\xi (H))$ hence $f(exp~t\xi (H))\equiv 0~.$ We infer that $exp~t\xi (H)$ is an integral curve of the system $\mathcal{P}$ contained in the level variety $h^{-1}(0).$

\vspace{2 mm}
\noindent
\newtheorem{cacatua}[LemmaCounter]{Lemma}
\begin{cacatua}
Let $\xi$ be a Lie vector field with hamiltonian \textit{h} and $\mathcal{N}$ a regularly embedded sub-manifold of $G(M)$ contained in $h^{-1}(0).$ Let us further assume that $\mathcal{N}$ is an integral manifold of $\mathcal{P}$ (i.e., $\Omega|\mathcal{N}=0$) and that $\xi_H\not\in T_H\mathcal{N}$ for any $H\in\mathcal{N}$ (i.e., the vector field $\xi$ is transversal to $\mathcal{N}$). Under these assumptions, the union of all the flow lines of $\xi$ containing points of $\mathcal{N}$ is, for small values of \textit{t}, a regularly embedded integral sub-manifold of $\mathcal{P}~.$ More precisely, we can choose, for any $H\in\mathcal{N}~,$ an open neighborhood $\mathcal{V}_H\subset\mathcal{N}$ and a constant $\epsilon_H$ such that $\cup_H\cup_{t,K}~exp~t\xi (K),$ $H\in\mathcal{N}~,$ $K\in\mathcal{V}_H~,$ $|t|<\epsilon_H~,$ is an integral manifold $\widetilde{\mathcal{N}}$ of $\mathcal{P}$ whose dimension is equal to $dim~\mathcal{N}+1~.$
\end{cacatua}

\vspace{2 mm}
\noindent
\textbf{Proof.} It will suffice to argue locally and glue together the local integral manifolds thus obtained. Let us take a point $H\in\mathcal{N}$ and choose appropriately a neighborhood $\mathcal{V}$ of \textit{H} in $\mathcal{N}$ as well as a constant $\epsilon >0$ such that $exp~t\xi(K),~K\in\mathcal{V}~,~|t|<\epsilon~,$ is defined. The transversality condition on $\xi$ implies that the mapping $]\epsilon,\epsilon[\times\mathcal{V}\longrightarrow G(M),$ $(t,K)~\mapsto~exp~t\xi(K),$ is an immersion for small \textit{t} whereupon its image is a regularly embedded sub-manifold for, eventually, a smaller constant $\epsilon$ and a smaller neighborhood $\mathcal{V}~.$ We can therefore assume that the image is in fact a sub-manifold $\widetilde{\mathcal{N}}$ that is furthermore contained in $h^{-1}(0)$ since \textit{h} is the hamiltonian of $\xi$ and $\mathcal{V}\subset h^{-1}(0).$ To prove that $\widetilde{\mathcal{N}}$ is an integral sub-manifold of $\mathcal{P}~,$ we first observe that the curves $\exp~t\xi(K)$ are integral curves of $\mathcal{P}$ and, since $\xi$ is an infinitesimal automorphism of this Pfaffian system, we infer that $\mathcal{V}_t=exp~t\xi(\mathcal{V})$ is, for any fixed value of \textit{t} , an integral manifold of $\mathcal{P}~.$ Next, we observe that the tangent space to $\widetilde{\mathcal{N}}~,$ at any point $L=\exp~t\xi(K),$ is generated by the tangent space $T_L\mathcal{V}_t$ together with the vector $\xi_L$ tangent to the curve $\exp~t\xi(K).$ It then follows readily that the form $\Omega$ vanishes on $T_L\widetilde{\mathcal{N}}$ and therefore that $\widetilde{\mathcal{N}}$ is an integral manifold of $\mathcal{P}~.$

\vspace{2 mm}
We now examine more closely the flow of the Lie vector field $\xi$ with hamiltonian \textit{h} (or \textit{f}). From the local expression (16), we infer that $\xi_H=0$ if and only if
\begin{equation}
f(H)=\frac{\partial f}{\partial p_i}(H)=(\frac{\partial f}{\partial x^i}+p_i\frac{\partial f}{\partial y})(H)=0~,
\end{equation}
a result that can be obtained as well, from the equation (15), as follows: We observe that these equations provide point-wise conditions on the vector $\xi_H$ namly, given any function \textit{f} and any point $H\in G(M),$ there exists a unique vector $v\in T_HG(M)$ for which
\begin{equation*}
i(v)\omega_H=f(H)\hspace{5 mm}and\hspace{5 mm}(df_H+i(v)d\omega_H=0)~.
\end{equation*}
In fact, let $v'$ be a second vector satisfying the above condition. Then $i(v-v')\omega_H=0$ \textit{i.e.}, $(v-v')\in\Sigma_H$ and $(i(v-v')d\omega_H)\wedge=0$ hence $(v-v')=0$ since $d\omega_H|\Sigma_H$ has maximal rank. We infer that $\xi_H=0$ if and only if
\begin{equation*}
f(H)=0\hspace{5 mm}and\hspace{5 mm}df_H\wedge\omega_H=0\hspace{2 mm}(or,~(i(\xi_H)d\omega_H)\wedge\omega_H=0)~,
\end{equation*}
the later condition asserting that $ker~df_H\supset ker~\omega_H=\Sigma_H.$ Therefore, if $df_H\wedge\omega_H\neq 0$ everywhere, the flow lines of $\xi$ are $1-$dimensional sub-manifolds of $G(M).$ Replacing now the function \textit{f} by $gf$, the equations $i(g~\!\xi)\omega=(gf)(H)$ and $(d(gf)+i(g~\!\xi)d\omega)\wedge\omega=0$ are satisfied at any point $H\in f^{-1}(0)$ hence the Lie vector field $\eta$ with hamiltonian $gf$ coincides with $g~\!\xi$ along the variety $h^{-1}(0).$ When \textit{g} is everywhere non-zero, the flow lines of $\xi$ and $g~\!\xi$ provide the same $1-$dimensional sub-manifolds, the flows differing only by their parametrizations.

\vspace{2 mm}
Coming back to partial differential equations of order 1, we defined them in the section 14 as being sub-manifolds $\mathcal{S}$ of co-dimension 1 in $G(M)$ and assumed that $\mathcal{P}|\mathcal{S}$ is a regular Pfaffian system on this same sub-manifold. Equivalently, this amounts to say that the restricted fundamental form $\Omega|\mathcal{S}$ never vanishes. Such a sub-manifold is defined locally by a non-singular equation $f=0$ (\textit{i.e.}, $df_H\neq 0$ for any $H\in\mathcal{S}$), any other regular equation, in the neighborhood of the same point $H_0$ being of the form $gf=0$ with $g(H)\neq 0$ for any $H\in\mathcal{S}.$ Such a local equation $f=0$ gives rise to a Lie vector field $\xi$ on $G(M)$ with hamiltonian equal to \textit{f}, this vector field being tangent to the sub-manifold $\mathcal{S}$ and therefore providing an infinitesimal automorphism of the Pfaffian system $\mathcal{P}|\mathcal{S}~.$ Further, since $<\xi,\omega>=f$ vanishes on $\mathcal{S}$ (locally equal to $f^{-1}(0)$), we infer that $\xi|\mathcal{S}$ is a characteristic vector field of $\mathcal{P}|\mathcal{S}.$ Finally, since the regularity condition on $\mathcal{P}|\mathcal{S}$ simply means that $(df\wedge\omega)_H~,~H\in\mathcal{S}$ is everywhere non-zero and since the Cauchy characteristics of $\mathcal{P}|\mathcal{S}$ are $1-$dimensional sub-manifolds, we infer that the vector field $\xi|\mathcal{S}$ is also everywhere non-zero and that its flow lines are open subsets of the characteristics, any other local regular equation $\tilde{f}=gf=0$ providing the same flow lines. In conclusion, the Cauchy characteristics of $\mathcal{P}|\mathcal{S}$ are unions of flow lines of Lie vector fields whose hamiltonians arise from local regular equations for the sub-manifold $\mathcal{S}.$ 

\vspace{2 mm}
We can now restate, in the context of contact transformations, the results of section 14 concerning the Cauchy characteristics of a single equation.

\vspace{2 mm}
\noindent
\newtheorem{catavento}[TheoremCounter]{Theorem}
\begin{catavento}
Let $\mathcal{S}$ be a first order partial differential equation on the manifold M, denote by $\xi$ the Lie vector field with hamiltonian f and assume that $f=0$ is a regular local equation of $\mathcal{S}.$ Then $\xi$ is tangent to $\mathcal{S}$ and the flow lines of the restricted vector field $\xi|\mathcal{S}$ are open subsets of the ($1-$dimensional) Cauchy characteristics of $\mathcal{P}|\mathcal{S},$ independently of the particular choice of the hamiltonian f. Any maximal (=n) dimensional integral manifold of the Pfaffian system above is the union of such flows and given any $(n-1)-$dimensional regularly embedded integral sub-manifold $\mathcal{N}$ of $\mathcal{P}|\mathcal{S}$ that is transversal to the Cauchy characteristics, it can be extended to an $n-$dimensional integral manifold by taking the union of all the flow lines of $\xi$ beginning at the points of $\mathcal{N}.$
\end{catavento}

\vspace{2 mm}
It is often useful to handle a version of the previous theorem when the function \textit{f} is replaced by the hamiltonian \textit{h} . With this in view, we claim that it suffices to determine conditions under which the equation $h=0$ is equivalent to a regular local equation $f=0~,$ where \textit{h} is a local section of the line bundle $\Phi,$ $\xi$ the Lie vector field with hamiltonian $h=<\xi,\Omega>~,$ $f=<\xi,\omega>$ and $h=f~\![\partial/\partial y].$

\vspace{2 mm}
We recall (\cite{Kumpera1972}) that, given any vector bundle \textit{E} with base manifold \textit{M} and projection $\pi~,$ the following short sequence is exact:
\begin{equation*}
0\longrightarrow T^*M\otimes E~\xrightarrow{\iota}~J_1E~\xrightarrow{\beta}~E\longrightarrow 0~,
\end{equation*}
where $J_1E$ is the first order jet bundle of the local sections of the fibration  $\pi:E~\longrightarrow~M~.$ Each $1-$jet identifies with an $n-$dimensional linear contact element of E, transverse to the fibres of $\pi$ ($n=dim~M$). In particular, the line bundle $\Phi$ with base space $G(M)$ provides the following short exact sequence:
\begin{equation*}
0\longrightarrow T^*G(M)\otimes\Phi~\xrightarrow{\iota}~J_1\Phi~\xrightarrow{\beta}~\Phi\longrightarrow 0~.
\end{equation*}

\vspace{2 mm}
\noindent
Given a hamiltonian $h\in\Gamma_{\ell}(\Phi)~,$ its image $\mathcal{H}$ is a sub-manifold of $\Phi$ and therefore $\mathcal{H}_1=\cup~T_y\mathcal{H}~,~y\in\mathcal{H}$ is also a sub-manifold of $J_1\Phi~.$ The set of points $H\in G(M)$ where \textit{h} vanishes is equal to the projection, in $G(M),$ of $\mathcal{H}\cap O_{\Phi}~,$ the second term being the null section, hence, at these points, $\beta(T_{h(H)}\mathcal{H})=0$ and therefore $T_{h(H)}\mathcal{H}=\iota(u_H)$ with $u_H\in T^*G(M)\otimes\Phi~.$ It is easy to check that $h=0$ is a regular equation if and only if $u_H\neq 0$ at any point $H\in G(M)$ where $h(H)=0~.$

\vspace{2 mm}
We now extend the previous results to systems of partial differential equations and state, firstly, some auxiliary lemmas, establishing thereafter the desired results. The remark following the the relation (10) and the fact that a linear combination, with constant coefficients, of Lie vector fields is again a Lie vector field implies the following

\vspace{2 mm}
\noindent
\newtheorem{piranha}[LemmaCounter]{Lemma}
\begin{piranha}
Let $\{\xi_{\alpha}\}$ be a family of Lie vector fields with associated hamiltonians $\{f_{\alpha}\},$ consider a fixed element $H\in G(M)$ and assume that the family $\{(df_{\alpha}\wedge\omega)_H\}$ is linearly independent. Then the family $\{\xi_{\alpha}(H)\}$ is also linearly independent and, if moreover $f_{\alpha}(H)=0$ for all $\alpha,$ the above sufficient condition becomes as well necessary. 
\end{piranha}

\vspace{2 mm}
\noindent
Using the local expression of the Lagrange bracket and an argument involving linear homogeneous ordinary differential equations similar to that preceding the Lemma 7, we also derive the following result.

\vspace{2 mm}
\noindent
\newtheorem{piracema}[LemmaCounter]{Lemma}
\begin{piracema}
Let $\{\xi_{\alpha}\}$ be a family of Lie vector fields with associated hamiltonians $\{f_{\alpha}\},$ and let $\{g_{\beta}\}$ be an arbitrary family of local functions on $G(M).$ In order that the family of vector fields $\{\xi_{\alpha}\}$ be tangent to the variety $\mathcal{W}$ defined by the equations $g_{\beta}=0$ (i.e., $<\xi_{\alpha},dg_{\beta}>=0$) it is necessary and sufficient that $[f_{\alpha},g_{\beta}]|\mathcal{W}=0$ for all the indices $\alpha$ and $\beta.$ Under these conditions and restricting our attention to finite subsets of indices $\alpha,$ the image of the exponential map $(t^{\alpha})~\mapsto~exp(\sum~t^{\alpha}\xi_{\alpha})(H)$ initiating at points $H\in\mathcal{W}$ is entirely contained in $\mathcal{W}.$
\end{piracema}

\vspace{2 mm}
\noindent
Again, the remark following the relation (17) and the local expression of the Lagrange bracket yield the following two lemmas.

\vspace{2 mm}
\noindent
\newtheorem{pirapora}[LemmaCounter]{Lemma}
\begin{pirapora}
Let $\{f_{\alpha}\}$ be a finite family of independent functions defined locally on $G(M),$ $\{\xi_{\alpha}\}$ the corresponding family of Lie vector fields, $\mathcal{W}$ the regularly embedded sub-manifold defined by the equations  $\{f_{\alpha}=0\},$ $\{g_{\beta}\}$ an arbitrary family of functions vanishing on $\mathcal{W}$ and $\{\eta_{\beta}\}$ the corresponding family of Lie vector fields. Under these conditions:

\vspace{2 mm}
(a) Each $\eta_{\beta}|\mathcal{W}$ is a linear combination of the restricted fields $\xi_{\alpha}|\mathcal{W}$ and, setting $g_{\beta}=\lambda^{\alpha}_{\beta}f_{\alpha},$ $\eta_{\beta}|\mathcal{W}=\lambda^{\alpha}_{\beta}\xi_{\alpha}|\mathcal{W}~.$ 

\vspace{2 mm}
(b) The condition $[f_{\alpha},f_{\alpha'}]|\mathcal{W}=0$ carries over to $[g_{\alpha},g_{\alpha'}]|\mathcal{W}=0~.$
\end{pirapora}

\vspace{3 mm}
\noindent
We observe that, unless all the Lagrange brackets vanish on $\mathcal{W}~,$ the vector fields $\xi_{\alpha}$ and $\eta_{\beta}$ are not necessarily tangent to $\mathcal{W}~,$ the notations $\xi_{\alpha}|\mathcal{W}$ and $\eta_{\beta}|\mathcal{W}$ simply indicating vector fields \textit{along} the sub-manifold $\mathcal{W}~.$

\vspace{2 mm}
\noindent
\newtheorem{piraquê}[LemmaCounter]{Lemma}
\begin{piraquê}
Let $\xi$ and $\eta$ be two Lie vector fields with hamiltonians f and g respectively and let us assume that $f(H)=g(H)=[f,g]_H=0,$ where $H\in G(M)$ is a given point. Then, $[\xi,\eta]_H$ is a linear combination of $\{\xi_H,\eta_H\}$ and, in coordinates, $[\xi,\eta]_H=(\partial g/\partial y)_H\xi_H-(\partial f/\partial y)_H\eta_H.$
\end{piraquê}

\vspace{3 mm}
\noindent
We next prove a result that extends the Lemma 9.

\vspace{2 mm}
\noindent
\newtheorem{pirajuba}[LemmaCounter]{Lemma}
\begin{pirajuba}
Let $\xi_{\alpha},~1\leq\alpha\leq\ell,$ be a family of linearly independent Lie vector fields, $f_{\alpha}$ the corresponding family of hamiltonians and $\mathcal{N}$ a regularly embedded sub-manifold of $G(M)$ contained in the variety $\mathcal{W}$ defined by the equations $\{f_{\alpha}=0\}.$ We assume, furthermore, that $\mathcal{N}$ is an integral manifold of $\mathcal{P},$ that $[f_{\alpha},f_{\beta}]|\mathcal{W}=0$ and that the linear subspaces $\Xi_H\subset T_HG(M)$ generated by the vectors $\{\xi_{\alpha}(H)\}$ are transversal to $\mathcal{N}$ at every point $H\in\mathcal{N}$ (i.e., $\Xi_H\cap T_H\mathcal{N}=0).$ Then, the union of all the $\ell-$dimensional flows $exp\sum t^{\alpha}~\!\xi_{\alpha}(H)$ initiating at points $H\in\mathcal{N}$ is, for small values of the parameters $(t^{\alpha})$ (in a sense analogous to that of the Lemma 9), a regularly embedded integral sub-manifold $\widetilde{N}$ of $\mathcal{P}$ contained in $\mathcal{W},$ its dimension being equal to $dim~\mathcal{N}+\ell.$
\end{pirajuba}

\vspace{2 mm}
\noindent
\textbf{Proof.} The condition $[f_{\alpha},f_{\beta}]|\mathcal{W}=0$ assures that the vector fields $\xi_{\alpha}$ are tangent to $\mathcal{W}$ and that each flow $exp~t^{\alpha}\xi_{\alpha}(H),$ initiating at a point $H\in\mathcal{W},$ is contained in $\mathcal{W}$ (the Lemma 11 in the special case when $f_{\alpha}=g_{\alpha}$). The transversality condition then implies that the union of all these flows initiating at points of $\mathcal{N}$ is, for small values of the parameters $t^{\alpha},$ a regularly embedded sub-manifold $\widetilde{\mathcal{N}}$ of dimension equal to $dim~\mathcal{N}+\ell.$ Further, the Lemma 10 asserts that the differential forms $\{df_{\alpha}\wedge\omega\}$ sre linearly inndependent at every point of $\mathcal{W}$ hence, \textit{a fortiori}, the differentials $\{df_{\alpha}\}$ are also linearly independent along $\mathcal{W}.$ It follows that $\mathcal{W}$ is a regularly embedded sub-manifold of $G(M).$ By the Lemma 13, the distribution $\Xi,$ defined on $\mathcal{W}$ by the vector fields $\xi_{\alpha}|\mathcal{W},$ is integrable and, of course, the flow $exp\sum t^{\alpha}~\!\xi_{\alpha}(H)$ is, in a neighborhood of \textit{H}, an open subset of the integral leaf of $\Xi$ containing the point \textit{H}. We next observe that each Lie vector field $\xi_{\alpha}$ is an infinitesimal automorphism of the system $\mathcal{P}$ and that its restriction $\xi_{\alpha}|\mathcal{W}$ is a characteristic vector field of $\Xi.$ The local $1-$parameter group generated by any linear combination, with constant coefficients, $\sum a^{\alpha}~\!\xi_{\alpha}(H)$ transforms integral manifolds of $\mathcal{P}$ into integral manifolds and preserves the leaves of $\Xi.$ Finally, let $K\in\widetilde{\mathcal{N}}$ be an element of the form $K=exp\sum t^{\alpha}~\!\xi_{\alpha}(H),$ with $H\in\mathcal{N},$ and denote by $\phi_t$ the local $1-$parameter group generated by the vector field $\sum a^{\alpha}~\!\xi_{\alpha}(H).$ Then, since the  transversality is preserved by diffeomorphisms, $T_K\widetilde{\mathcal{N}}=T\phi_1(T_H\mathcal{N})\oplus T\phi_1(\Xi_H)$ and on account of the above remarks, it follows that $\Omega$ vanishes on $T\phi_1(T_H\mathcal{N})=T_K\phi_1(\mathcal{N})$ as well as on $T\phi_1(\Xi_H)=\Xi_K~,$ hence it also vanishes on $T_K\widetilde{\mathcal{N}},$ the sub-manifold $\widetilde{\mathcal{N}}$ being therefore an integral manifold of $\mathcal{P}$ contained in $\mathcal{W}.$

\vspace{2 mm}
\noindent
\textbf{Remark.} The existence of the integral manifold $\widetilde{\mathcal{N}}$ being only assured locally, we can relax the Lagrange bracket condition by requiring that $[f_{\alpha},f_{\beta}]|\mathcal{W}$ vanishes only in a neighborhood of $\mathcal{N},$ and we can as well only require that the vector fields $\{\xi_{\alpha}\}$ be linearly independent along the sub-manifold $\mathcal{N}~.$ We could also compute directly the tangent space, at the point \textit{K} , to the flow $exp\sum t^{\alpha}~\!\xi_{\alpha}(H)$ by using the formula of the Proposition 4.1 in \cite{Kumpera1982}. However, the Lemma 13 seems to be a more convenient technical device.

\vspace{2 mm}
We finally need the following linear algebraic result:

\vspace{2 mm}
\noindent
\newtheorem{piracanjuba}[LemmaCounter]{Lemma}
\begin{piracanjuba}
Let $\mu$ be an exterior $2-$form of rank $2n$ defined on a real or complex vector space V, H a hyperplane of V, $\mu_1$ the induced form on H and $A=\{v\in V|i(v)\mu=0\}$ the annihilator of $\mu.$

\vspace{2 mm}
(a) If $dim~V=2n$ then $rank~\mu_1=2n-2~.$

\vspace{2 mm}
(b) If $dim~V=2n+1$ and if $A\cap H=0,$ then $rank~\mu_1=2n~.$

\vspace{2 mm}
(c) If $dim~V=2n+1$ and if $A\cap H\neq 0,$ then $rank~\mu_1=2n-2~.$
\end{piracanjuba}

\vspace{2 mm}
\noindent
\textbf{Proof.} (a). Since $dim~H=2n-1,$ the form $\mu_1$ has a non-trivial annihilator $A_1.$ However, $v\in V~\mapsto~i(v)\mu\in V^*$ being an isomorphism, the kernel of the restricted map $v\in V~\mapsto~(i(v)\mu)|H\in H^*$ is $1-$dimensional and therefore $dim~A_1=1~.$

\vspace{2 mm}
(b) and (c). Under the present hypotheses, $dim~A=1~.$ If $A\cap H=0~,$ then the annihilator $A+1$ is necessarily trivial otherwise any non-null vector $v\in A_1$ would also belong to \textit{A}. If $A\cap H\neq 0$ or, equivalently, if $A\subset H~,$ then $A\subset A_1$ and, factoring all the data to the quotient space $V/A~,$ we recover the data of the first case. 

\vspace{2 mm}
\noindent
\newtheorem{picanha}[TheoremCounter]{Theorem}
\begin{picanha}
Let $\mathcal{S}$ be a system of regular first order partial differential equations or rank $q<n+1,$ on the manifold M, and $\{f_{\alpha}=0\}$ a set of local equations of $\mathcal{S}.$ Let us also denote by $\xi_{\alpha}$ the Lie vector field with hamiltonian $f_{\alpha}$ and assume that $[f_{\alpha},f_{\beta}]|\mathcal{S}=0$ for all the indices $\alpha$ and $\beta.$ Then,

\vspace{2 mm}
(a) The vector fields $\{\xi_{\alpha}\}$ are tangent to $\mathcal{S}$ and the restrictions $\{\xi_{\alpha}|\mathcal{S}\}$ are everywhere linearly independent. 

\vspace{2 mm}
(b) Each $\xi_{\alpha}|\mathcal{S}$ is a characteristic vector field for the system $\mathcal{P}|\mathcal{S}$ and the family $\{\xi_{\alpha}|\mathcal{S}\}$ generates the characteristic distribution of $\mathcal{P}|\mathcal{S}.$

\vspace{2 mm}
(c) Each $q-$dimensional flow $exp\sum t^{\alpha}~\!\xi_{\alpha}(H),~H\in\mathcal{N},$ is locally (i.e., for small values of the parameters $t^{\alpha}$) independent of the particular choice for the local equations of $\mathcal{S},$ being therefore intrinsically associated to the system $\mathcal{S}.$ Moreover, it is locally an open subset of the Cauchy characteristic containing the point $H.$

\vspace{2 mm}
(d) Any maximal ($=n$) dimensional integral manifold $\mathcal{N}$ of the Pfaffian system $\mathcal{P}|\mathcal{S}$ is the union of flows $exp\sum t^{\alpha}~\!\xi_{\alpha}(H)$ with $H\in\mathcal{N}.$

\vspace{2 mm}
(e) Given an $(n-q)-$dimensional regularly embedded integral manifold $\mathcal{N}$ of $\mathcal{P}|\mathcal{S},$ transversal to the Cauchy characteristics, it extends to an $n-$dimensional integral manifold $\widetilde{\mathcal{N}}$ by taking the union of all the flows $exp\sum t^{\alpha}~\!\xi_{\alpha}(H)$ initiating at points $H\in\mathcal{N}.$
\end{picanha}

\vspace{2 mm}
\noindent
\textbf{Proof.} The regularity condition $\Omega|\mathcal{S}_H\neq 0~,~H\in\mathcal{S},$ is obviously equivalent to the linear independence of the forms $\{\omega_H,(df_{\alpha})_H\}$ and a simple calculation shows that the later property is equivalent, inasmuch, to the linear independence of the forms $\{(df_{\alpha}\wedge\omega)_H\}$ hence, since $\mathcal{S}$ is regular, we infer from the Lemma 10 that the vectors $\{(\xi_{\alpha})_H\}$ are linearly independent along $\mathcal{S}.$ The tangency, to the manifold $\mathcal{S},$ of the vector fields $\{\xi_{\alpha}\}$ follows, by the Lemma 11, from the bracket condition $[f_{\alpha},f_{\beta}]|\mathcal{S}=0~.$ Observing that each $f_{\alpha}$ vanishes on $\mathcal{S},$ we also infer that $(\xi_{\alpha})_H\in(\Sigma|\mathcal{S})_H$ for any $H\in\mathcal{S},$ hence the restrictions $\xi_{\alpha}|\mathcal{S}$ are characteristic vector fields of the system $\mathcal{P}|\mathcal{S}.$ The remaining statements are now direct consequences of the Theorem 6 for $dim~\mathcal{S}=2n+1-q$ and therefore the dimension of the characteristic distribution of $\mathcal{P}|\mathcal{S}$ is equal to $(2n+1-q)-(2n-2q+1)=q~.$

\vspace{2 mm}
It is often useful and always instructive to carry out a complete proof independently of the Theorem 6 and the argument runs as follows. Denote by $\Xi$ the distribution, on the manifold $\mathcal{S},$ generated by the vector fields $\xi_{\alpha}|\mathcal{S}.$ Then $\Xi$ is intrinsically associated to the system $\mathcal{S}$ for, by the Lemma 12, any other set of local equations for $\mathcal{S}$ provides the same distribution. We also infer from this Lemma that the bracket condition is as well intrinsical. Next, we know from the proof of the Lemma 14 that each flow $exp\sum t^{\alpha}~\!\xi_{\alpha}(H),~H\in\mathcal{N}$ is, locally, an open subset of the integral leaf of $\Xi$ containing \textit{H} hence two such flows, associated to different sets of local equations for $\mathcal{S},$ will agree in a neighborhood of \textit{H}, this settling part of the item (c). Clearly, the item (e) is a restatement of the Lemma 14. Since $\mathcal{P}$ and therefore $\mathcal{P}|\mathcal{S}$ admit integral manifolds of dimensions at most equal to \textit{n} and since a flow is the union of all the integral curves, initiating at \textit{H} , of the characteristic vector fields $\sum a^{\alpha}~\!\xi_{\alpha}$ (within the adequate bounds for the coefficients $a^{\alpha}$), the statement (d) is a direct consequence of the Lemma 14, in the special case where $\ell=1~,$ for these characteristic fields must then be tangent to $\mathcal{N}.$ Finally, to prove that the flows $exp\sum t^{\alpha}~\!\xi_{\alpha}(H),~H\in\mathcal{S},$ are locally open subsets of the Cauchy characteristics, it is enough to prove that the vector fields $\{\xi_{\alpha}|\mathcal{S}\}$ generate the characteristic distribution of $\mathcal{P}|\mathcal{S}.$ For this, we consider the systems $\mathcal{R}_1=\{f_1=0\},$ $\mathcal{R}_2=\{f_1=f_2=0\},$ $\cdots$, $\mathcal{R}_q=\mathcal{R}=\{f_{\alpha}=0\},$ and argue as follows: The restriction $d\omega|\Sigma_H$ is a $2-$form of rank $2n$ defined on the vector space $\Sigma_H$ of dimension $2n$ and $(\Sigma|R_1)_H$ is a hyperplane of $\Sigma_H$ ($\Sigma|R_1$ means the characteristic system of $R_1$). Therefore, according to the Lemma 15, the rank of $(d\omega|\Sigma|R_1)_H$ is equal to $2n-2~.$ Similarly, $(\Sigma|R_2)_H$ is a hyperplane in the $(2n-1)-$dimensional vector space $(\Sigma|R_1)_H$ and therefore, the rank of $d\omega|(\Sigma|R_2)_H$ is either $2n-2$ or $2n-4~.$ Proceeding inductively, we infer that the rank of $d\omega|(\Sigma|R_p)_H$ is at least equal to $2n-2p$ hence, at the final stage, the rank of $d\omega|(\Sigma|R)_H$ is at least equal to $2n-2q~.$ It follows that the dimension of the characteristic sub-space $\widetilde{\Delta}_H$ of $\mathcal{P}|\mathcal{S}$ is, at most, equal to $dim~(\Sigma|R)_H-(2n-2q)=(2n-q)-(2n-2q)=q~.$ Since $\Xi_H\subset\tilde{\Delta}_H$ and since $dim~\Xi_H=q~,$ we conclude that $\Xi=\widetilde{\Delta}.$ The reader will have noticed that the present argument is nothing but a fancy way to deal with the Cartan argument introduced in the proof of the Theorem 6.

\vspace{2 mm}
\noindent
\textbf{Remark.} We could, of course, enlarge integral manifolds of $\mathcal{P}|\mathcal{S}$ by taking, step by step, either individual Lie vector fields $\xi_{\alpha}$ or sub-families $\{\xi_{\alpha_i}\}$ since such sub-families generate, in view of the Lemma 13, integrable sub-distributions of $\Xi~.$

\vspace{2 mm}
\noindent
When $q=1~,$ the above theorem reduces to the Theorem 7, the bracket condition then becoming  trivial. When $q=n~,$ the $n-$dimensional distribution $\Sigma|\mathcal{S}$ is, under the hypotheses of the Theorem, generated by its characteristic vector fields and therefore is completely integrable.

\vspace{2 mm}
We now examine the case $q=n+1$ excluded in the statement of the Theorem. Since the manifold $\mathcal{S}$ has, in this case, the dimension \textit{n} , the system admits $n-$dimensional solutions, containing any of its points, if and only if $\Omega|\mathcal{S}=0~,$ the regularity hypothesis $\Omega|\mathcal{S}\neq 0$ being void of sense.

\vspace{2 mm}
\noindent
\newtheorem{piracicaba}[LemmaCounter]{Lemma}
\begin{piracicaba}
Let $\mathcal{N}$ be an $n-$dimensional integral manifold of the system $\mathcal{P}$ and $\{f,g\}$ two functions vanishing on $\mathcal{N}.$ Then $[f,g]$ also vanishes on $\mathcal{N}.$
\end{piracicaba}

\vspace{2 mm}
\noindent
In fact, the corresponding Lie vector fields $\xi$ and $\eta$ must be tangent to $\mathcal{N}$ otherwise we could enlarge this manifold to an $(n+1)-$dimensional integral manifold of $\mathcal{P}.$ If, for example, $\xi_H$ were transverse to $\mathcal{N}$ then $(df\wedge\omega)_H\neq 0$ hence $df_H\neq 0$ and the level variety $f^{-1}(0)$ would become a manifold, in a neighborhood of \textit{H} , containing $\mathcal{N}.$ We could then, according to the Lemma 9, enlarge the integral manifold $\mathcal{N}$ this. however, being excluded. The tangency of $\xi$ and $\eta$ implies that of $[\xi,\eta]$ and therefore the hamiltonian $[f,g]=<[\xi,\eta],\omega>$ vanishes on $\mathcal{N}$ since $\omega|\mathcal{N}=0~.$

\vspace{2 mm}
\noindent
\newtheorem{piracica}[LemmaCounter]{Lemma}
\begin{piracica}
Let $\mathcal{N}$ be an $n-$dimensional sub-manifold of $G(M)$ defined locally by the independent equations $\{f_{\alpha}=0\}$ and let us assume that $[f_{\alpha},f_{\beta}]|\mathcal{N}=0~.$ Then $\Omega|\mathcal{N}=0$ and, furthermore, all the Lie vector fields $\xi_{\alpha}$ with hamiltonians $\{f_{\alpha}\}$ are tangent to $\mathcal{N},$ their restrictions generating $T\mathcal{N}.$
\end{piracica}

\vspace{2 mm}
\noindent
The proof will be omitted for being almost a repetition. The above Lemma shows that a system $\mathcal{N}$ of rank $n+1$ and satisfying the above bracket condition is in fact an integral manifold of $\mathcal{P}.$ The associated Pfaffian system $\mathcal{P}|\mathcal{N}$ is therefore null, its characteristics are regular and the characteristic distribution, equal to $T\mathcal{N},$ is generated by the vector fields $\xi_{\alpha}|\mathcal{N}.$

\vspace{2 mm}
\noindent
\newtheorem{pamonha}[CorollaryCounter]{Corollary}
\begin{pamonha}
In order that an $n-$dimensional sub-manifold $\mathcal{N}$ of $G(M)$ be an integral manifold of $\mathcal{P},$ it is necessary and sufficient that $[f_{\alpha},f_{\beta}]|\mathcal{N}$ $=0$ whatever the choice of the local equations $\{f_{\alpha}=0\}$ for $\mathcal{N}.$
\end{pamonha}

\vspace{3 mm}
\noindent
We can now settle the integrability problem for the system $\mathcal{S}$ namely, the existence of $n-$dimensional integral manifolds for the associated Pfaffian system $\mathcal{P}|\mathcal{S}.$

\vspace{2 mm}
\noindent
\newtheorem{alcatra}[TheoremCounter]{Theorem}
\begin{alcatra}
A regular system $\mathcal{S}$ is integrable if and only if $[f_{\alpha},f_{\beta}]|\mathcal{S}=0$ whatever the choice of the local equations $\{f_{\alpha}=0\}$ for $\mathcal{S}.$ 
\end{alcatra}

\vspace{2 mm}
\noindent
When $rank~\mathcal{S}=q=n+1~,$ the result is a restatement of the last Corollary. Assuming that $q\leq n~,$ the necessity of the condition follows from the Lemma 16 and the sufficiency from the Theorem 8 for, under the stated condition of regularity, the Pfaffian system $\mathcal{P}|\mathcal{S}$ has regular characteristics with the appropriate dimension (\textit{cf.}, Corollary 6).

\vspace{2 mm}
\noindent
\textbf{Remark.} We infer, from the Lemma 12, that the integrability condition, expressed in terms of the Lagrange bracket, is intrinsically associated to the equation $\mathcal{S}.$ It could as well be restated as follows. The bracket of any two functions, vanishing on $\mathcal{S},$ also vanishes on $\mathcal{S}.$ Geometrically speaking, the condition states that, in a neighborhood of any point $H\in\mathcal{S},$ there exist \textit{q} Lie vector fields $\xi_{\alpha}$ with independent hamiltonians vanishing on $\mathcal{S}$ and such that $[\xi_{\alpha},\xi_{\beta}]_y\in\Sigma_y$ whenever $y\in\mathcal{S}.$ Since the Lagrangean bracket $[f,g]$ reduces, on the variety defined by the equations $f=g=0~,$ to the Jacobi bracket (\textit{cf.}, section 16), it follows that the above bracket condition boils down, when $q=n~,$ to the classical Frobenius integrability criterion.

\vspace{2 mm}
We end up this section by discussing the infinitesimal contact automorphisms of the partial differential equations under consideration. The results will become useful in later sections.

\vspace{2 mm}
\noindent
Clearly, a Lie vector field $\xi$ is, in restriction to $\mathcal{S},$ an infinitesimal automorphidm of the Pfaffian system $\mathcal{P}|\mathcal{S}$ if and only if it is tangent to $\mathcal{S}.$ Such a vector field will be called an infinitesimal contact automorphism of the equation $\mathcal{S}.$ The Lemma 11 then specializes to the following result.

\vspace{2 mm}
\noindent
\newtheorem{peruibe}[LemmaCounter]{Lemma}
\begin{peruibe}
A Lie vector field $\eta$ with hamiltonian g is an infinitesimal contact automorphism of the equation $\mathcal{S}$ if and only if $[g,f_{\alpha}]|\mathcal{S}=0$ for any choice of local equations $\{f_{\alpha}=0\}$ for $\mathcal{S}.$ 
\end{peruibe}

\vspace{2 mm}
\noindent
Given an element $H\in\mathcal{S},$ we denote by $\overline{g}_H$ the Lie algebra of all the germs of infinitesimal contact automorphisms of $\mathcal{S},$ at the point \textit{H}, and set
\begin{equation*}
g_H=\overline{g}_H|\mathcal{S}~,\hspace{10 mm}h_H=\{g:[g,\overline{f}_{\alpha}]|\mathcal{S}=0\}~,
\end{equation*}
where \textit{g} is a germ of a function at the point \textit{H} and $\{\overline{f}_{\alpha}=0\}$ indicates the germs, at \textit{H}, of a set of local equations for $\mathcal{S}$ at \textit{H}.

\vspace{2 mm}
\noindent
\newtheorem{alcatrada}[TheoremCounter]{Theorem}
\begin{alcatrada}
The Lie algebra $g_H$ is isomorphic to $h_H/(f_{\alpha})^2,$ where $(f_{\alpha})$ denotes the ideal (in the sense of associative algebras) generated, in the algebra of all the germs, by the functions $f_{\alpha}.$
\end{alcatrada}

\vspace{2 mm}
\noindent
\textbf{Proof.} Let $\eta$ be an infinitesimal contact automorphism of $\mathcal{P}|\mathcal{S}$ with hamiltonian \textit{g} and let us assume that $\eta|\mathcal{S}=0~.$ Then $g|\mathcal{S}=0$ and, writing locally $g=\sum~\lambda^{\alpha}f_{\alpha}~,$ where $\{f_{\alpha}\}$ is a set of local equations for $\mathcal{S},$ we infer that $0=(dg\wedge\omega)_H=\sum~\lambda^{\alpha}(H)(df_{\alpha}\wedge\omega)_H$ whenever $H\in\mathcal{S},$ hence $\lambda^{\alpha}|\mathcal{S}=0~.$ Conversely, if $g\in(f_{\alpha})^2,$ then $\eta|\mathcal{S}=0~.$ It should be noted that the infinitesimal contact automorphisms of $\mathcal{P}|\mathcal{S}$ coincide with those of its characteristic system.

\vspace{2 mm}
In all the preceding discussion, the scalar hamiltonians \textit{f} can always be replaced by $\Phi-$valued hamiltonians \textit{h} except eventually when expressions involving $dh$ and $h^2$ should occur. The differential $dh$ is meaningful and takes the role of $df$ only at those points where \textit{h} vanishes. However, the expression $h^2$ is meaningless and the ideal $(f_{\alpha})^2$ has to be replaced by the subset of all the germs of $\Phi-$valued sections \textit{h} vanishing to order 1 on $\mathcal{S}.$

\vspace{2 mm}
Among the many Lie vector fields that we had the occasion to deal with in this section, let us recall those introduced in the very beginning and obtained, by prolongation, from vector fields coming from the base space \textit{M}. In view of the Lemma  and the local expression (16), the hamiltonians \textit{f} of such prolonged vector fields are characterized, locally, by the condition $\partial^2f/\partial p_i\partial p_j\equiv 0$ hence are semi-linear functions with respect to the variables $p_i~.$

\section{Brackets.}
Given an odd classical contact structure $\mathcal{P}$ on a manifold \textit{M} of dimension $2n+1~,$ we know that the defining Pfaffian system is locally generated by a $1-$form $\omega$ with Darboux class equal to $2n+1~.$ The generator $\omega$ can be chosen globally (\textit{i.e.}, with domain equal to \textit{M}) if and only if \textit{M} is orientable. The pair $(M,\omega),$ where $\omega$ is global and of constant Darboux class equal to $dim~M,$ is known as a \textit{Pfaffian structure} and it induces, of course, a contact structure $\mathcal{P}.$ We derive, in this section, some results concerning Pfaffian structures and introduce, in particular, the associated Jacobi and Poisson brackets. When global generators are not available, as is the case for the canonical contact structure on $G(M),$ these results will be applied on the domains of the local generators $\omega.$

\vspace{2 mm}
Let $(M,\omega)$ be a Pfaffian structure and $\mathcal{P}$ the induced contact structure. Since the dimension of the characteristic system $\Delta$ of $d\omega$ is equal to 1 (\textit{i.e.}, $dim~\Delta_x=1$), there exists a canonically defined vector field $\xi_0$ satisfying the equations
\begin{equation*}
i(\xi_0)d\omega=0\hspace{5 mm}and\hspace{5 mm}i(\xi_0)\omega=1~.
\end{equation*}
This vector field is an infinitesimal automorphism of $\omega$ \textit{i.e.}, $\theta(\xi_0)\omega=0$ hence, \textit{a fortiori}, an automorphism of $\mathcal{P}$ namely, the Lie vector field with hamiltonian equal to the constant function 1~. The distribution $\Sigma=ker~\omega$ and the characteristic distribution $\Delta$ of $d\omega$ define a direct sum decomposition $TM=\Delta\oplus\Sigma~.$ Denoting by $\pi_{\Sigma}$ the projection upon the factor $\Sigma~,$ any vector field $\xi$ on \textit{M} can be decomposed according to the formula
\begin{equation}
\xi=f\xi_0+\pi_{\Sigma}\xi~,  
\end{equation}
where $f=i(\xi)\omega~.$ In particular, $\xi$ is a Lie vector field if and only if $(df+i(\pi_{\Sigma}\xi)d\omega)\wedge\omega=0$ and, this being the case, $df+i(\pi_{\Sigma}\xi)d\omega=(i(\xi_0)df)\omega~.$ Furthermore, the vector field $\xi$ is an infinitesimal automorphism of $\omega$ if and only if $i(\xi_0)df=0$ hence if and only if the function \textit{f} is a first integral of the characteristic distribution $\Delta.$ Denote by $\mu$ the linear map that assigns to each vector field $\eta$ tangent to $\Sigma$ (\textit{i.e.}, $\eta_x\in\Sigma_x$) the Pfaffian form $i(\eta)d\omega~.$ Since $i(\xi_0)[i(\eta)d\omega]=0~,$ we infer that $\mu$ establishes an isomorphism between the module of vector fields tangent to $\Sigma$ and the module of $1-$forms vanishing on the characteristic distribution $\Delta.$ A simple calculation provides the formula
\begin{equation*}
\xi=f\xi_0+\mu^{-1}(<\xi_0,df>\omega-df)    
\end{equation*}
for any Lie vector field $\xi~,$ where \textit{f} is the hamiltonian of $\xi~.$ In this particular setting, the Lagrange bracket is given by
\begin{equation}
[f,g]=i(\xi)dg-gi(\xi_0)df~,
\end{equation}
where $\xi$ is the Lie vector field with hamiltonian \textit{f} . Introducing the Jacobi bracket
\begin{equation}
\{f,g\}=i(\xi)dg-fi(\xi_0)dg~,
\end{equation}
the formula above reduces to the more symmetrical expression
\begin{equation}
[f,g]=\{f,g\}+f~\!i(\xi_0)dg-g~\!i(\xi_0)df~,
\end{equation}
and it comes out that
\begin{equation*}
\{f,g\}=d\omega(\xi,\eta)~,
\end{equation*}
where $\eta$ is the Lie vector field with hamiltonian \textit{g} or, equivalently,
\begin{equation*}
\{f,g\}~\!\omega\wedge(d\omega)^n=n~\!df\wedge dg\wedge\omega\wedge(d\omega)^{n-1}~.
\end{equation*}
In adapted local coordinates $(\mathcal{U};x^i,y,p_i)$ \textit{i.e.}, coordinates reducing $\omega$ to its canonical form $\omega=dy-p_idx^i~,$ we obtain the expressions:
\begin{equation}
\{f,g\}=\sum~\frac{\partial f}{\partial p_i}(\frac{\partial g}{\partial x^i}+p_i\frac{\partial g}{\partial y})-\sum~\frac{\partial g}{\partial p_i}(\frac{\partial f}{\partial x^i}+p_i\frac{\partial f}{\partial y})~,
\end{equation}
\begin{equation*}
\xi_0=\partial/\partial y\hspace{5 mm}and\hspace{5 mm}i(\xi_0)dg=\partial g/\partial y~.
\end{equation*}
It is useful to remark, as a consequence of the expression (21), that the Lagrange bracket $[f,g]$ reduces to the Jacobi bracket $\{f,g\}$ on the sub-variety defined by the equations $f=g=0$ or whenever \textit{f} and \textit{g} are first integrals of the characteristic system $\Delta.$ In the later case, $\frac{\partial f}{\partial y}=\frac{\partial g}{\partial y}=0$ and the local expression (22) shows that the Jacobi bracket is simply the Poisson bracket with respect to the $2-$form $d\omega.$ Geometrically speaking, this can be described as follows: The characteristic system $\Delta$ of $d\omega$ being integrable, we can consider, locally, the quotient manifold $\overline{M}$ of \textit{M} modulo the characteristics, factor $d\omega$ to a closed $2-$form $\overline{dw}$ of maximum rank on  $\overline{M}$ and a first integral \textit{f} to a function $\overline{f}.$ Next, we recall that a Lie vector field $\xi$ with hamiltonian \textit{f} is an infinitesimal automorphism of $\omega$ if and only if the function \textit{f} is a first integral of $\Delta.$ Such a vector field $\xi$ preserving $\omega$ will also preserve $d\omega,$ thereafter preserving as well the characteristic system of $d\omega~.$ We infer that $\xi$ factors to a vector field $\overline{\xi}$ defined on $\overline{M}$ and the relation $df+i(\pi_{\Sigma}\xi)d\omega=0$ implies that $d\overline{f}=-i(\overline{\xi})\overline{d\omega}$ hence $\overline{\xi}$ is a hamiltonian vector field with symplectic hamiltonian $\overline{f}.$ The Poisson bracket $(f,g),$ relative to $d\omega,$ factors to the usual Poisson bracket $(\overline{f},\overline{g})$ with respect to the symplectic form $\overline{d\omega}.$

\vspace{2 mm}
Let us finally consider a Lie vector field $\xi$ with (contact) hamiltonian \textit{f} and let \textit{g} be a generic function on the manifold \textit{M}. According to the formula (20),
\begin{equation*}
\{f,g\}=(\xi-f~\!\xi_0)g~,
\end{equation*}
hence the \textit{adjoint action} $g~\mapsto~\{f,g\},$ relative to a fixed function \textit{f} , is equal to the action of the vector field $\xi_f=\xi-f~\!\xi_0~.$

\vspace{2 mm}
\noindent
\newtheorem{peruibinho}[LemmaCounter]{Lemma}
\begin{peruibinho}
Let $\xi$ and $\eta$ be two Lie vector fields, defined on the manifold M, with hamiltonians f and g respectively and let us assume that $\{f,g\}_x=0$ at some point $x\in M.$ Then $[\xi,\eta]_x$ is a linear combination of the vectors $\xi_x$ and $\eta_x$ that, translated in terms of coordinates, yields the relation
\begin{equation*}
[\xi,\eta]_x=(\partial g/\partial y)_x\xi_x-(\partial f/\partial y)_x\eta_x~.
\end{equation*}
\end{peruibinho}

\vspace{2 mm}
\noindent
We can now derive a classical result, due to Sophus Lie (see the references), concerning the integrability of \textit{complete} systems.

\vspace{2 mm}
\noindent
\newtheorem{almofada}[TheoremCounter]{Theorem (Lie)}
\begin{almofada}
Let $(f_i),~1\leq i\leq q,$ be a family of functions defined on the manifold M and satisfying the Jacobi bracket condition $\{f_i,f_j\}=0$ for all the indices i and j. Then the system of linear partial differential equations
\begin{equation*}
\{f_1,g\}=0,~\cdots~,\{f_q,g\}=0~,
\end{equation*}
in the unknown function g is complete.
\end{almofada}

\vspace{2 mm}
\noindent
To say that the system is complete means that the (eventually singular) distribution generated by the vector fields $\overline{\xi}_i=\xi_i-f_i\xi_0$ is integrable, $\xi_i$ being the Lie vector field with hamiltonian $f_i$ and $\xi_0$ the vector field naturally associated to the canonical contact structure, as defined previously. When the Lie vector fields $\xi_i$ generate, at each point $x\in M,$ a linear sub-space transverse to $\xi_0~,$ then the linear independence of the vector fields $\overline{\xi}_i$ is equivalent to that of the $\xi_i~.$ We can therefore apply the criterion provided by the Lemma 10. Assuming now that the vector fields $\overline{\xi}_i$ are independent, the completeness, in the sense given to it by Sophus Lie, means the existence of $dim~\!M-q$ independent solutions $g_j$ of the system $\{f_i,g\}=0~,~1\leq i\leq q~.$

\section{Higher order differential equations.}
Let \textit{M} be a manifold of dimension \textit{n} , $x_0\in M$ a given though arbitrary point and $(N_1,N_2)$ a couple of sub-manifolds of \textit{M} both containing the point $x_0$ and both having their dimensions equal to $q~.$ We say that the two sub-manifolds have a \textit{contact} of order $k,~k\geq 0~,$ at the point $x_0$ when there exists a local diffeomorphism 
\begin{equation*}
\varphi:U_1\subset N_1~\longrightarrow~U_2\subset N_2~,
\end{equation*}
defined on an open neighborhood $U_1$ of the point $x_0$ and having as target an open neighborhood $U_2$ of the same point, such that the two mappings $\iota_1$ and $\iota_2\circ\varphi$ have a contact of order \textit{k} at the point $x_0,$ where $\iota_j$ denotes the inclusion map. Stated more explicitly in terms of $k-$jets, the above contact requirement simply means that $j_k\iota_1(x_0)=j_k(\iota_2\circ\varphi)(x_0)~,$ the notation $j_kf(x)$ indicating the $k-$jet of the mapping \textit{f} taken at the point \textit{x} . The above defined relation is, of course, an equivalence relation, the equivalence classes being called the $k-$th order tangent contact elements to the manifold \textit{M}. Further details and results on this subject can be found in \cite{Kumpera1987}. The set of all the $k-$th order contact elements with a given dimension \textit{q} is endowed with a differentiable manifold structure, denoted by $G_k^q(M),$ and the "base point" map $\pi:G_k^q(M)~\longrightarrow~M$ is a locally trivial fibre bundle. First order contact elements are simply $q-$dimensional linear sub-spaces contained in the tangent spaces to the manifold \textit{M} hence, in particular, $G(M)=G_1^{n-1}(M),$ where $n=dim~M~.$ We observe that $G_k^q(M),$ also known as the $k-$th order Grassmann bundle, admits an everywhere dense open subset namely, the bundle of all the  $k-$jets of local sections of a fibration (submersive surjection) $M~\longrightarrow~M_0~.$ For this, it suffices to map the $k-$th order jet $j_k\sigma(y)$ of a local section $\sigma~,$ at the point $y\in M_0~,$ into the $k-$th order equivalence class of the sub-manifold $im~\!\sigma$ at the point $x=\sigma(y).$ We list below some useful properties verified by these higher order contact elements.

\vspace{2 mm}
(a) For any $r<q,$ every $r-$dimensional contact element of order \textit{k} extends to a not necessarily unique $q-$dimensional contact element of the same order \textit{k} .

\vspace{2 mm}
(b) For any $h<k,$ every $h-$th order contact element of dimension \textit{q} extends to a not necessarily unique $k-$th order contact element of the same dimension \textit{q} .

\vspace{2 mm}
(c) We can also combine and iterate the above two extension procedures and the concept of \textit{inclusion} as well as that  of \textit{iteration}, for contact elements, is not only well defined but has moreover a wide scope of possibilities.

\vspace{2 mm}
(d) Any local diffeomorphism $\varphi:U~\longrightarrow~U'$ of \textit{M} extends (prolongs) canonically to a local diffeomorphism $\wp\varphi:\mathcal{U}~\longrightarrow~\mathcal{U'}$ of the Grassmann bundle $G_k^q(M)$ and, actually, this prolongation procedure can also be carried out for any local diffeomorphism $\varphi$ initiating on a manifold \textit{M} and terminating on a different manifold $M'$.

\vspace{3 mm}
Let us next define the $k-$th order \textit{canonical contact structure} namely, a Pfaffian system on the manifold $G_k^q(M)$ that will play a fundamental role in the study of the integrability problem for $k-$th order partial differential  equations. We call \textit{holonomic} the contact elements as defined earlier and introduce now a notion of \textit{semi-holonomy}. Given a $q-$dimensional sub-manifold $N\subset M$ and fixing a point $x_0\in M,$ let us denote by $H_{k,x_0}\in (G^q_k)_{x_0}$ the $k-$th order contact element determined by the above data. Considering the sub-manifold $N_{k-1}\subset G^q_{k-1}$ and, subsequently, taking the first order contact element determined by this sub-manifold at the corresponding point $H_{k-1,x_0}~,$ let us indicate by $(H^q_1,H^q_{k-1,x_0})$ the corresponding first order contact element thus obtained and call it \textit{semi-holonomic}. This element is clearly well determined since the initially given contact element provides all the necessary derivatives up to the order \textit{k} in order to be so. We thus obtain a map 
\begin{equation*}
sh_{k-1,k}:G^q_k(M)~\longrightarrow~G^q_1(G^q_{k-1}(M))
\end{equation*}
that we call the \textit{semi-holonomic inclusion} and a similar procedure also provides the more general semi-holonomic inclusion
\begin{equation*}
sh_{k-h,k}:G^q_k(M)~\longrightarrow~G^q_h(G^q_{k-h}(M)).
\end{equation*}

\vspace{2 mm}
We can now introduce the \textit{canonical contact structure}, of order 1 , on each manifold $G^q_k$ by proceeding recurrently on the order \textit{k} . For $k=1~,$ define $\mathcal{P}=\mathcal{P}_1$ by setting $\mathcal{P}_{H_x}=(T\pi_{H_x})^*(H_x)^{\perp}$ and in general, for $\mathcal{P}=\mathcal{P}_k~,$
\begin{equation*}
\mathcal{P}_{H_{k,x}}=(T\rho_{k,k-1})^*(sh_{k-1,k}(H_{k,x}))^{\perp}+(T\rho_{k-1,k})^*\mathcal{P}_{k-1,H_{k-1,x}}~,
\end{equation*}
where $\rho_{h,k}$ denotes the natural projection of the $k-$th order contact elements onto the $h-$th order contacts elements and the tangent maps are taken at the appropriate points. The main properties of these Pfaffian systems can be stated in the following theorem, where $N_k$ denotes the sub-manifold of $G^q_k(M)$ obtained by taking all the $k-$th order contact elements defined by the sub-manifold $N\subset M.$

\vspace{2 mm}
\newtheorem{almofadinha}[TheoremCounter]{Theorem (Lie)}
\begin{almofadinha}
The following two properties describe the basic attributes and importance of the canonical contact structures at any order:

\vspace{2 mm}
(a) A sub-manifold $\overline{N}\subset G^q_k(M)$ is of the form $N_k$ if and only if it is an integral manifold of $\mathcal{P}_k$ transverse to the projection $\pi$ onto M.

\vspace{2 mm}
(b) A local diffeomorphism $\overline{\varphi}$ of $G^q_k(M)$ is of the prolongued form $\wp\varphi$ if and only if it is projectable onto a local diffeomorphism $\varphi$ on M and if, moreover, it is an automorphism of the contact structure $\mathcal{P}_k.$  
\end{almofadinha}

\vspace{2 mm}
\noindent
Similar properties are true inasmuch for the infinitesimal automorphisms.

\vspace{2 mm}
We denote, as previously, by $\mathcal{P}|\mathcal{S}$ or $\mathcal{P}_k|\mathcal{S}$ the restriction of the contact system to a differential equation $\mathcal{S}$ of order \textit{k} and adopt as well the notation $\mathcal{N}_k$ to indicate the manifold composed by all the $k-$th order contact elements that contain the contact element determined by the sub-manifold $\mathcal{N}$ of \textit{M}. This being so, we can now reformulate all the results of the previous sections in terms of the differential equations of order \textit{k}. However, we shall not enter here into the details since these would be essentially a copy of what is already written, the index 1 being replaced by \textit{k} .

\section{The local equivalence problem.}
We discuss here the equivalence problem in the realm of Ehresmann's prolongation spaces of which, Jet spaces and Grassmannian contact bundles are special cases (\cite{Ehresmann1958},\cite{Ehresmann1967},\cite{Ehresmann1966}). All the emphasis is given to the first order canonical contact structure associated to a partial differential equations since the local or, inasmuch, the global equivalence of these Pfaffian systems will entail the corresponding equivalence for the equations. A fundamental ingredient in our approach is the consideration of \textit{merihedric} prolongation spaces (\textit{prolongements mériédriques}) as defined and used by Élie Cartan in many of his writings. Though rather absent in the recent literature, the reader will find interesting examples in \cite{Kumpera1999}, \cite{Kumpera2014} and \cite{Kumpera2015}. As for the calculations, we try to reduce them to the strict minimum, the main concern being the  determination of a \textit{fundamental set of invariants} associated to a differential system that will characterise as well as describe the local equivalences.

\vspace{2 mm}
Let us recall, briefly, some very standard facts so as to fix the notations and subsequently discuss, in more detail, the aims and techniques of the next two sections.

\vspace{2 mm}
\noindent
As is usual, a \textit{fibration} is a surjective map $P~\longrightarrow~M$ of maximum rank and we define a \textit{multi-fibration} as being a sequence, finite or infinite in length,

\begin{equation*}
\cdots~\longrightarrow~P_\mu~\longrightarrow~P_{\mu-1}~\longrightarrow~\cdots~\longrightarrow~P_1\longrightarrow P_0 ~,
\end{equation*}

\vspace{3 mm}
\noindent
where all the projections $\rho_{\beta,\alpha}:P_\alpha\longrightarrow P_\beta$ , $\alpha\geq\beta$ , are fibrations, the map $\rho_{\beta,\alpha}$ being the composite of the successive projections. Following Cartan, we say that a multi-fibration is a \textit{merihedric} prolongation space with respect to the basic manifold $P_0$ when there is defined a prolongation algorithm $\wp$ such that any local transformation $\varphi:U~\longrightarrow~V~,$ on the manifold $P_0~,$ can be prolonged (extended) to a local transformation $\varphi_{\mu}=\wp_{\mu}(\varphi)$ operating on the manifold $P_{\mu}~,$ the usual commutativity relations with respect to the projections, the composition of transformations and the passage to the inverses being preserved. However, we shall only require, with the obvious notations, that the open sets $U_{\mu}$ and  $V_{\mu}$ project surjectively onto $U_{\nu}$ and $V_{\nu}$ whenever $\nu\leq\mu$ without requiring that the former be the inverse images of the latter. Next, we require that to every local sub-manifold $\mathcal{N}$ of $P_0$ corresponds, at each level $\mu~,$ a \textit{prolonged} sub-manifold $\wp\mathcal{N}$ and that these sub-manifolds project one upon the other. Finally we require that, at each level $P_\mu~,~\mu\geq 1~,$ a Pfaffian system $\mathcal{C}_{\mu}$ be defined and that it satisfies the following properties:

\vspace{3 mm}
(a) $\rho^*_{\mu-1,\mu}\mathcal{C}_{\mu-1}\subset\mathcal{C}_{\mu}~,$

\vspace{2 mm}
(b) for any local sub-manifold $\mathcal{N}$ of $P_0~,$ the space $\wp_{\mu}\mathcal{N}$ is an integral sub-manifold of $\mathcal{C}_{\mu}$ and, conversely, any local sub-manifold $\overline{\mathcal{N}}$ of $P_{\mu}~,$ transverse to the fibres over $P_0$ and integral for the previous system is, locally, the prolongation of a sub-manifold in the base space,    

\vspace{2 mm}
(c) the $\mu-$th prolongation of any local transformation of $P_0$ is an automorphism of $\mathcal{C}_{\mu}~,$ the converse being verified just locally.

\vspace{3 mm}
\noindent
A similar statement also holds for vector fields (\textit{infinitesimal transformations}) and their prolongations, where the operations are now the linear operations together with the Lie derivative and the bracket of vector fields. We also require that this \textit{infinitesimal}  prolongation procedure be compatible with the \textit{finite} prolongation of the elements of 
a local $1-$parameter group. 

\vspace{2 mm}
\noindent
Besides the Grassmannians, among the most required prolongation spaces are, of course, the sequences of Ehresmann's Jet spaces of arbitrary orders associated to a given fibration $\pi:P~\longrightarrow~M~,$ where we consider $k-$jets of local sections. No details are provided on Jet spaces since these are assumed to be a well known subject.

\vspace{2 mm}
Let us finally indicate a prolongation procedure for differential equations \textit{i.e.}, sub-manifolds $\mathcal{S}$ contained in merihedric prolongation spaces by adopting exactly the same definitions as those given in the case of Jet spaces. In fact, let us take a sub-manifold $N\subset P_0$ with its dimension equal to that prescribed for the solutions of a given equation $\mathcal{S}\subset P_{\mu}$ and say that it is $\ell-$admissible, at the point $x\in N,$ when it satisfies the following two properties:
 
\vspace{3 mm}
(a) $\wp_{\mu}(N)_x\in\mathcal{S}~,$ where the index \textit{x} indicates the point of $\wp_{\mu}(N)$ that corresponds to \textit{x} and, at that point,
 
\vspace{2 mm}
(b) the sub-manifold $\wp_{\mu}(N)$ is tangent, to order $\ell~,$ to $\mathcal{S}.$
 
\vspace{3 mm}
\noindent
We define, thereafter, the $\ell-$th prolongation $\wp_{\ell}(\mathcal{S})$ of $\mathcal{S}$ as being the totality of the $(\mu+\ell)-$th order contact elements obtained by means of the $\ell-$admissible sub-manifolds and calculated at the admissibility points. All the properties of the standard prolongation algorithm for differential equations, defined in Jet spaces, transcribe in this more general setting and, in particular, the $\ell-$th prolongation of $\mathcal{S}$ can be obtained by iterated first order prolongations. 
 
\vspace{2 mm}
In studying the equivalence, local or global, of two differential equations, we shall always consider both equations defined over the same basic manifold for, if the two equations were situated above distinct basic manifolds, we could always shift one onto the other by taking a local or global diffeomorphism whose prolongations will also transport all the other data. Concerning the basic properties and operations relative to differential equations and their prolongations, we refer the reader to \cite{Kumpera1972}. 

\vspace{2 mm}
Let us now start examining the notion of \textit{equivalence}. This notion must respect the solutions of the equations \textit{i.e.}, a solution of one of the equations must be transformed into a solution of the other hence must respect the sophisticated structure of the Grassmannian spaces.

\vspace{2 mm}
\newtheorem{lagrangean}[DefinitionCounter]{Definition}
\begin{lagrangean}
Two $k-$th order differential equations $\mathcal{S}$ and $\mathcal{S}'$ are said to be locally absolutely equivalent in the neighborhoods of $X\in\mathcal{S}$ and $X'\in\mathcal{S}'$ when there exists a local transformation $\varphi:U~\longrightarrow~U'~,$ defined on the base manifold $P_0~,$ such that its prolongation $\wp_k\varphi$ transforms X into X' and, further, becomes a diffeomorphism when restricted to the appropriate open subsets of $\mathcal{S}$ and $\mathcal{S}'.$ The equivalence is global when $\varphi$ is global.
\end{lagrangean}

\vspace{2 mm}
\noindent
Obviously, $\wp_k\varphi^{-1}$ is also a local equivalence and any one of them transforms solutions into solutions. If, for instance, $\mathcal{N}$ is a solution of $\mathcal{S},$ then $\varphi(\mathcal{N})$ is a solution of $\mathcal{S}'.$  

\vspace{2 mm}
\noindent
The notion of absolute equivalence seems too restrictive when confronted with applications and we give below a broader and more suitable definition. Needless to say that both definitions are stated in Cartan's \textit{Mémoire} \cite{Cartan1914}.

\vspace{2 mm}
\newtheorem{laguère}[DefinitionCounter]{Definition}
\begin{laguère}
Two differential equations are said to be locally equivalent when they admit prolongations of certain orders that are locally absolutely equivalent. Global equivalence has a similar definition.
\end{laguère}

\vspace{2 mm}
\noindent
We observe that this definition of equivalence enables us to confront equations with different orders in which case the two prolongations are also of distinct orders.

\vspace{2 mm}
The first example occurs in the case of \textit{ordinary differential equations}. In this case, when the image spaces have the same dimension (the same number of dependent variables), any two equations are always locally equivalent at non-singular points. When singularities are present, the matter becomes considerably more involved. As a second example, we mention that all the integrable Pfaffian systems with the same rank and co-rank are locally equivalent. As for the global equivalence, in the present case as well as in the case of ordinary differential equations, this is a beautiful and still unresolved problem in topology and a very interesting situation can be found in \cite{Almeida1985}. Let us now have a glance at non-integrable Pfaffian systems, a rather nice example being provided by the \textit{Flag Systems}. These are systems of co-rank equal to 2 (when the characteristics vanish) hence the integral manifolds are just $1-$dimensional curves, \textit{l'intégrale générale ne dépendant que de fonctions arbitraires d'une seule variable}, as would claim Cartan. However, the reader might delight himself in examining the many local models exhibited in \cite{Kumpera2014}.

\vspace{2 mm}
These Flag Systems provide a very nice illustration of the merihedric context and, further, are prolongation spaces of finite length. Before doing so, we recall (\cite{Kumpera1999}) that, in Jet spaces, the \textit{total derivative} of functions can be trivially extended to a \textit{total Lie derivative} of differential forms and, with this in mind, we shall give attention to Pfaffian systems instead of partial differential equations, which was by far Cartan's preference. In this respect, we also recall that the total Lie derivative of the $k-$th order canonical contact structure, defined on a Jet space $J_k\pi$ or on a Grassmannian bundle $G^q_k(M)~,$ is precisely the corresponding contact structure on the $(k+1)-$st level (where we shall have to add, for obviuos reasons, the pullback of the contact structure defined on the first level, $k=1$). Let us then return to the sequence
\begin{equation*}
(P_\ell,S)~\longrightarrow~(P_{\ell-1},\overline{S}_1)~\longrightarrow~\cdots~\longrightarrow~(P_0,\overline{S}_{\ell}) ~,
\end{equation*}
exhibited in \cite{Kumpera2014}, p.8, $l.-3,$ \textit{S} being a Flag System of length $\ell$ defined on the space $P=P_{\ell},$ having null characteristics and where $\overline{S}_{\ell-\nu}$ denotes a Pfaffian systems canonically isomorphic to the successive $\nu-$th derived systems of \textit{S} though defined on a more appropriate spaces so as to display, as well, null characteristics. The above finite sequence is a prolongation space together with Pfaffian systems, at each level, naturally associated to the initially given system \textit{S}. The terminating system $\overline{S}_\ell$ vanishes and $\overline{S}_{\ell-1}$ is a  Darboux system in $3-$space. As shown in \cite{Kumpera2014}, every local or infinitesimal automorphism of this Darboux system extends (prolongs) canonically to an automorphism of $(P_{\ell-\nu},\overline{S}_{\nu}),$ this correspondence becoming, moreover, an isomorphism of \textit{automorphism} pseudo-groups as well as of \textit{automorphism} pseudo-algebras. We are not to be concerned with local sections since $P_1$ can be considered as the total (base) space, whereas the base space $P_0$ collapses to an open set, eventually to a point. We claim that each Pfaffian system $\overline{S}_{\nu}$ is the prolongation of $\overline{S}_{\nu+1}.$ As is, not only it seems that we are walking backwards but worse, the prolongation algorithm as described previously is somehow missing. Nevertheless, in the passage from $\overline{S}_{\nu+1}$ to $\overline{S}_{\nu}~,$ three options are possible (\cite{Kumperaa1982}) and, upon choosing the appropriate one, we can thereafter apply the prolongation method as described earlier. Furthermore, if we choose the appropriate coordinates so as that the form $\omega=dy-zdx$ generates the Darboux system $\overline{S}_{\ell-1}~,$ then total derivatives reappear again and Cartan, as always, is perfectly right in his claims. 

\section{The local equivalence of differential systems}
We shall now examine a few aspects of this local equivalence problem where most enlightening examples are provided by the study of geometrical structures since, in most cases, the local equivalence is examined on the level of the defining equations for these structures (\cite{Kumpera2017}).

\vspace{3 mm}
\noindent
$(a)~Ordinary~Differential~Equations.$

\vspace{2 mm}
An ordinary differential equation of order 1 in \textit{n} unknown functions is simply a vector field on a manifold \textit{P}, where we take the fibratiom $\pi:P~\longrightarrow~\mathbf{R}$ with $dim~P=n+1~.$ A $k-$th order equation is a contact vector field (infinitesimal automorphism of the corresponding contact structure) defined on $J_k\pi~$ or, inasmuch, on $G^1_k(P).$ As mentioned earlier, any two ordinary differential equations are always locally equivavalent at two non-singular points. As for the global equivalence, this is a topological problem that involves the nature of the global solutions and, of course, compactness is the main ingredient. We also observe that two ordinary differential equations of distinct orders can be locally, \textit{viz.} globally, equivalent in view of Cartan's definitions. As for the local equivalence in neighborhoods of singular points, there is an immense literature on this subject hence we shall not enter here into the details relative to this matter.

\vspace{2 mm}
\noindent
$(b)~Linear~Partial~Differential~Equations.$

\vspace{2 mm}
Such equations are defined, most conveniently, as being vector sub-bundles of the $k-$th order Jet spaces of local sections of vector bundles \textit{i.e.}, the \textit{total space} of the initial fibration is the total space of a vector bundle and, consequently, the corresponding Jet space is as well a vector bundle. Nevertheless, $k-$th order contact elements defined on the total spaces of vector bundles can do as well since they also have a vector bundle structure. In this case, it is most appropriate to restrict the local equivalences to linear morphisms of vector bundles though two linear equations can actually be equivalent via a non-linear local transformation. This is not a flaw but an immense virtue since it enables us to replace a non-linear equation by a linear one. A first systematic investigation on the local equivalence of linear differential systems was initiated by Drach, Picard and Vessiot, whereupon resulted the well known \textit{Théorie de Picard-Vessiot} (\textit{cf.} \cite{Kumpera2016}).

\vspace{3 mm}
\noindent
$(c)~General~Differential~Equations.$

\vspace{2 mm}
Who undoubtedly most contributed not only with results but mainly with ideas and methods, in this general equivalence problem, was Élie Cartan and the reader is invited to examine the table of contents of the \textit{O\!euvres Complètes, Partie II}. Very significant contributions were also brought by Bernard Malgrange, in particular for the case of elliptic equations, as well as by Masatake Kuranishi for the case of the defining equations of a Lie pseudo-group namely, those equations that have the additional algebraic structure of a groupoid. Here again, we shall not enter into further details since the subject is by far too extensive.

\vspace{3 mm}
\noindent
$(d)~Pfaffian~Systems.$

\vspace{2 mm}
Initial contributions were, of course, brought by Pfaff himself and later expanded by contemporary Sophus Lie (\textit{Gasammelte Abhandlungen, Vol.6}) and Mark von Weber. But again, the most striking contributions were placed forward by no other than Élie Cartan who gave us much insight into the local equivalence of non-integrable systems. As mentioned earlier, any two integrable Pfaffian systems are \textit{everywhere} locally equivalent if and only if they have the same rank and co-rank (the underlying manifolds ought to have the same dimension). It is also a straightforward consequence of the linearity, in the tangent and co-tangent bundles, that any two Pfaffian systems, integrable or not, are always equivalent, \textit{to first order}, if and only if their ranks and co-ranks are equal.   

\vspace{3 mm}
\noindent
$(e)~Exterior~Differential~Systems.$

\vspace{2 mm}
\noindent
All the merits belong here to Élie Cartan who, in fact, introduced such systems in the study of various geometrical problems (\cite{Cartan1945}). Nevertheless, important contributions are also due to Masatake Kuranishi.

\section{The local equivalence of Lie groupoids}
Since we are interested in the local equivalence of partial differential equations, the sole Lie groupoids to be considered in the sequel are those whose total spaces are sub-manifolds in some jet space and we begin here with some general considerations.

\vspace{2 mm}
\noindent
Given the fibration $\pi=p_1:P=M\times M~\longrightarrow~M$ and the corresponding $k-$th order Jet bundle $J_k\pi~,$ we denote by $\Pi_kM$ the groupoid of all the $k-$jets of invertible sections of $\pi$ (\cite{Kumpera2015}). In other terms, $\Pi_kM$ is the set of all the invertible $k-$jets with source and target in \textit{M} and is an open dense subset of $J_k\pi~.$ The Lie groupoids we shall be interested in are the sub-groupoids of $\Pi_kM$ that are inasmuch sub-manifolds (not necessarily regularly embedded) of the ambient groupoid and, together with the differentiable sub-manifold sructure, become Lie groupoids in the sense of Ehresmann (\textit{cf.}, the above citation). Such $k-$th order Lie groupoids will, most often, be denoted by $\Gamma_k~.$ We next observe that the standard prolongation of order \textit{h} of a Lie groupoid, in the sense of partial differential equations, is again a Lie groupoid at order $k+h$ and all the usual properties relative to differential equations transcribe for Lie groupoids. Moreover, it can be shown that the above-mentioned properties still hold in the more general context of merihedric prolongations. In considering the local equivalences, we proceed in analogy with what is usually done in Lie group theory and restrict our attention to open neighborhoods of the unit elements of the groupoids. More precisely, we look for conditions under which there exists a local contact transformation that maps isomorphically an open neighborhood of the unit elements in one of the groupoids onto a similar neighborhood of the other. The term "isomorphically" means that it preserves the algebraic as well as the differentiable properties and "an open neighborhood of the unit elements" refers to a neighborhood of an open set in the units sub-manifold. Since a contact transformation, as above, operating on the sub-manifolds is by necessity induced by a local contact transformation operating in the ambient Jet space, we ultimately consider, in our setting, a local contact transformation
\begin{equation*}
\varphi:\mathcal{U}~\longrightarrow~\mathcal{U}'
\end{equation*}
that transforms accordingly the respective intersections with the grou-poids and where $\mathcal{U}$ and $\mathcal{U}'$ are open sets in $\Pi_kM.$

\vspace{2 mm}
Let us next consider the equation $\mathcal{R}$ composed by all the elements $X\in\Pi_kM$ whose sources belong to the units of the first groupoid and the targets to the units of the second. Then, of course, any local equivalence will be a local solution of this equation the converse being also true and where a local transformation of the base space operates on $\Pi_kM$ via the conjugation by its $k-$jets. Let us next observe that the first groupoid operates, to the right, on the space $\mathcal{R}$ and that the second groupoid operates to the left. So as to render the data more homogeneous and regular, we shall assume that both these actions are transitive, such an assumptions being much in accordance with the problem considered. In fact, most relevant examples do conform to this requirement and, more important, most contexts where equivalence is sought for also do comply with it. It should be noted that whenever the local equivalence is verified and a local contact transformation is put forward, then one transitivity assumptions gives rise to the other. It now remains to investigate the equation $\mathcal{R}$ as well as its prolongations and determine the desired invariants that will characterize the local equivalence. A first obvious remark is that the dimensions of both groupoids inasmuch as the dimensions of their units sub-spaces must be equal and this also entails the equality, in dimensions, of both the $\alpha$ and the $\beta$ fibres. The standard as well as the merihedric prolongation algorithms do preserve some data and operations but not all required for our purpose. In fact, the prolongation of a groupoid is still, algebraically, a groupoid but its differentiable structure is not guaranteed. Furthermore, the action of the prolonged groupoid on the prolonged equation is not forcibly transitive. On the other hand and concerning the equation $\mathcal{R}~,$ similar statement can also be affirmed. Nevertheless, it should be emphasized that the above stated transitivity assumption implies that this equation inherits, from either groupoid, a manifold structure rendering it a sub-manifold of $\Pi_kM~.$ In what follows, we list several necessary conditions for our problem to have a solution.

\vspace{3 mm}
\noindent
$\mathbf{1.~Differentiability.}$ The iterated prolongations of both groupoids are also Lie groupoids and sub-manifolds of the respective jet spaces.

\vspace{3 mm}
\noindent
$\mathbf{2.~Dimension.}$ Both prolonged groupoids have, at each order, the same dimension. We do not need to bother about the dimensions of the unit spaces since these spaces are equal (diffeomorphic) to those of the initially given groupoids,

\vspace{3 mm}
\noindent
$\mathbf{3.~Transitivity.}$ Both prolonged groupoids operate, at each order, transitively on the prolonged equation. Consequently, each prolonged equation is a sub-manifold of the corresponding jet space. 

\vspace{3 mm}
\noindent
$\mathbf{4.~Symbols.}$ The symbols of both prolonged groupoids have, at each order, the same dimension. 

\vspace{3 mm}
\noindent
$\mathbf{5.~\delta-cohomology.}$ Both iterated prolongations of the groupoids have equidimensional $\delta-$cohomologies and thereafter become simultaneously $2-$acyclic (\textit{cf.}\cite{Kumpera1972}).

\vspace{3 mm}
\noindent
We now assume that the above five properties do hold for the given two Lie groupoids of jets. On account of the transitivity of the groupoid actions, it follows that the iterated prolongations of the equation $\mathcal{R}$ become also $2-$acyclic at orders not greater than those for the groupoids. Finally, if we further assume that all the data is real analytic, then the Cartan Theorem\footnote{most improperly called the Cartan-Kähler Theorem} guarantees the desired local equivalence. We also observe that the knowledge of Jordan-Hölder resolutions for the two groupoids, \textit{a priori} assumed to be equivalent, gives rise to the possibility of constructing composition series, for the equation $\mathcal{R}~,$ that eventually become Jordan-Hölder. Unfortunately, we cannot make the same statement in the $C^{\infty}$ realm since, under these weaker assumptions, the whole Cartan Theory breaks down. There are examples of non-analytic though involutive (2-acyclic) equations that do not possess any  solution. In this last realm, much further work is awaiting to be done.

\vspace{2 mm}
The local equivalence of analytic Lie groupoids entails the local equivalence of Lie pseudo-groups since these are defined with the help of such groupoids, their defining equations.

\section{Integral contact elements}
Let $\mathcal{P}$ be a not necessarily integrable Pfaffian system of constant rank \textit{r} defined on the manifold \textit{M} of dimension \textit{n} , denote by
\begin{equation*}
\{\omega^1,\omega^2,~\cdots~,\omega^r\}
\end{equation*}
a free system of local generators, defined in a neighborhood of a 
generic point $x_0~,$ and by $\Sigma$ the annihilator of $\mathcal{P}$ in $TM~,$ thus obtaining a field of $(n-r)-$dimensional first order contact elements annihilated by $\mathcal{P}.$

\vspace{2 mm}
Let us take an arbitrary point $x_0\in M$ and try to determine, \textit{de proche en proche} as would say Cartan, an ascending chain of integral contact elements at the point $x_0~.$ We firstly take a one-dimensional subspace $l\subset \Sigma_{x_0}~,$ a line in the terminology of Cartan, and examine when such a linear subspace can be included in a higher dimensional \textit{integral} contact element. The sole condition $l\subset\Sigma_{x_0}$ is not sufficient and the entire subsequent discussion reposes on the following remark:

\vspace{3 mm} 
\textit{Given an integral sub-manifold $\mathcal{I}$ of the Pfaffian system $\mathcal{P}$ and an arbitrary local section $\omega$ of this system, the restrictions $\iota^*\omega$ and $\iota^*d\omega$ vanish identically on $\mathcal{I}~,$ where $\iota:\mathcal{I}~\hookrightarrow~M$ is the inclusion.}

\vspace{2 mm}
\noindent
It suffices, of course, to check the above condition for the sole generators $\omega^i$ and furthermore, this condition transcribes by the following two requirements: $i(v)\omega^i=0~,~v\in l~,~v\neq 0~,$ and $i(v)d\omega^i\equiv 0\hspace{3 mm}mod~\{\omega^i\},$ the latter condition also rewriting by $(d\omega^i)_{x_0}\equiv 0\hspace{3 mm}mod~\{\omega^i\}.$

\vspace{2 mm}
We now fix a given $1-dimensional$ linear integral contact element $E_1$ issued at the point $x_0$ (\textit{i.e.}, contained in $\Sigma_{x_0}$) and look for the conditions under which a plane contact element $E_2$ ($dim=2$) containing the given linear element $E_1$ is also integral. In order to find such an element, it suffices to look for a vector $w\in T_{x_0}M$ linearly  independent from \textit{v} and verifying the two conditions $i(w)\omega^i=0$ and $i(w)i(v)d\omega^i=d\omega(v,w)=0~.$ We shall say that \textit{w} is \textit{in involution} with \textit{v} (in Cartan's terminology, the line $\lambda$ generated by \textit{w} is in involution with the line \textit{l}). Linearity then implies that all the vectors \textit{w} that are in involution with \textit{v} (adding, of course, the null vector) form a linear subspace $\widetilde{E}_1$ containing all the $2-$dimensional integral subspaces considered above. 

\vspace{2 mm}
Let us next fix a $2-dimensional$ integral contact element $E_2$ (a linear $2-$dimensional subspace of $\widetilde{E}_1$) and consider the vectors
\textit{z} that are in involution with all the vectors belonging to $E_2.$ It suffices, of course, to consider those vectors that are in involution, simultaneously, with those belonging to any given basis of $E_2$ and we shall denote by $\widetilde{E}_2$ the linear subspace consisting of all these vectors that obviously is contained in $\widetilde{E}_1.$

\vspace{2 mm}
Fixing a $3-dimensional$ contact element $E_3$ that is contained in  $\widetilde{E}_2$ and contains $E_2~,$ we can pursue this nice \textit{Origami} game and construct an ascending chain $E_0\subset E_1\subset E_2\subset E_3~ \cdots~,$ each element satisfying the involutiveness condition with respect to the preceding one. We include the null space $E_0$ just to please Cartan where it figures, in his writings, as the one point set $E_0=\{x_0\}~.$ Dimensions being finite, there will be an integer $\rho_k$ for which the above process breaks down namely, there will not exist any other non-null vector, in involution with $E_{\rho_k}~,$ other than those already obtained previously. 

\vspace{2 mm}
\newtheorem{uvw}[DefinitionCounter]{Definition}
\begin{uvw}
 According to Cartan, we shall say that the integer $s_1=n-r-\rho_k$ is the \textit{character} (caractère) of $\mathcal{P}$ at the point $x_0$ and relative to the above ascending chain.
\end{uvw}

\vspace{2 mm}
\noindent
Described more geometrically, the character is simply the difference between the values of $dim~\Sigma_{x_0}$ and $dim~E_{\rho_k}~.$ We next describe Cartan's most ingenious method for determining these maximal integral elements. Let us take a tangent contact element \textit{E} contained in some $\Sigma_{x_0}$ as well as a second sub-space $E'\subset E~.$ The sub-space $E'$ is said to be \textit{characteristic} with respect to $E$ whenever any element of $E'$ is in involution with every element of $E.$ The sum of all the characteristic sub-spaces is then a maximal characteristic sub-space called the \textit{characteristic element} of $E$ with respect to $\mathcal{P}.$ Though $\mathcal{P}$ is assumed to have null characteristics, these characteristic elements need not be null. Next, given any two sub-spaces \textit{E} and \textit{F} contained in $\Sigma_{x_0}~,$ we shall say that they are \textit{conjugate} when their characteristic elements are in involution, this meaning that any vector belonging to one of them is in involution with all the vectors belonging to the other. In order to obtain an integral element of dimension $\mu\leq\rho_k~,$ it suffices to take $\mu$ independent vectors $\{v_1,~\cdots~,v_{\mu}\}$ for which the corresponding sub-spaces $\widetilde{E}_{1,i}~,~1\leq i\leq\mu~,$ are pairwise conjugate and then simply take their intersection.

\vspace{2 mm}
We next define $s_j=dim~\widetilde{E}_j~,$ for  $1\leq j\leq\rho_k$ and at the given point $x_0~$. Obviously, $s_1-s_2\geq s_2-s_3\geq~\cdots~\geq s_{\rho_k-1}-s_{\rho_k}$ since, at each step, we have a lesser choice of vectors to look for and, furthermore, $\rho_k=n-r$ if and only if $\mathcal{P}$ is integrable in a neighborhood of $x_0.$ We should also observe that this point-wise technique is most appropriate for the study of singular Pfaffian systems in view of obtaining the maximal integral sub-manifolds whose dimensions can vary not only from point to point but, in fact, also at the same point since these maximal integral elements depend upon the choice of the chain of intermediate integral elements. Moreover, two distinct maximal integral elements can have different dimensions as well as a common intersection whose dimension is greater than $1~.$ This is the case, for example, occuring with the systatical systems (systèmes systatiques, \cite{Cartan1901}) where two distinct maximal integral manifolds can have, at a point $x_0~,$ a first order tangency of dimension $\geq 2$ or even a common sub-manifold of such a dimension. The integers $s_j$ are called the successive \textit{enlarged characters} of $\mathcal{P}$ at the point $x_0~.$ We finally call the attention of the reader to the way Cartan embraced his thoughts. Invariably and persistently, he would always pursue the \textit{co-variant} trail and, in his \textit{Mémoire} \cite{Cartan1901} as well as previously in \cite{Cartan1900}, this shows up as soon as he describes the above-mentioned contact elements. In fact, these are depicted with the help of convenient linear combinations of the differentials $dx^i$ relative to a local coordinate system $(x^1,~\cdots~,x^n)$ defined in a neighborhood of the point $x_0~.$ In other terms, Cartan employs linear bases for the sub-spaces, in $T^*M~,$ that annihilate the \textit{contra-variant} subspaces $E_j$ and $\widetilde{E}_j$ considered previously ([\textit{nous}] \textit{regardons $dx_1,dx_2,~\cdots~,dx_r$ comme les paramètres directeurs d'une droite issue de ce point}, pg.4, \textit{l.}8). In so doing, the calculations become much simpler and completely straightforward though, unfortunately, our \textit{contra-variant} intuition of the geometric world around us suffers with it.

\section{Integral sub-manifolds}
We first observe that any Pfaffian system whith rank not exceeding $n-1$ always admits integral $1-$dimensional sub-manifolds, at a given point $x_0~,$ since it suffices to integrate any local vector field $\xi$ that is annihilated by the system and that satisfies $\xi_{x_0}\neq 0~.$ The above assertion can fail for singular systems for, in this case, it might happen that $\Sigma_{x_0}=0$ at a given point $x_0~.$ Hence, there is in principle no integrability condition pending for the dimension 1. Secondly, we observe that the dimensions of all the previously considered contact elements are \textit{lower semi-continuous} (the dimensions tend to increase) and, further, that the \textit{pertinence relation} $\in$ is a \textit{closed} relation. Consequently, these properties enable us to choose locally, in a neighborhood of $x_0~,$ a finite set of linearly independent vector fields $\xi_{\mu}~,~1\leq\mu\leq\rho_k~,$ such that

\vspace{3 mm}
\textbf{1.} The vectors $\{\xi_{\mu}(x_0)\}$ generate $E_{\rho_k}$ and

\vspace{2 mm}
\textbf{2.} The vector fields $\xi_{\mu}$ generate, in that neighborhood, an integrable distribution (field of contact elements).

\vspace{3 mm}
\noindent
We can then integrate this local distribution and obtain an integral sub-manifold $\mathcal{I}$ of the initially given Pfaffian system $\mathcal{P}~,$ this integral manifold being of maximal dimension relative to the point $x_0$ and to the ascending chain. More than that, we can actually obtain a \textit{sub-foliation} where the leaf containing $x_0$ is a maximal integral sub-manifold, at that point.\footnote{These sub-foliations having been obtained by means of an entirely naïve argument, it is hardly comprehensible why the author of \cite{Stormark2000} goes into so much trouble in defining \textit{maximal involutions} and thereafter \textit{complete subsystems} (pg.2 of the introduction and a subsequent later section). Furthermore, the presence of the Cauchy-Kovalevskaya theorem, valid only in the real analytic realm, is still more incomprehensible.} We should however observe that it is not unique relative to the point $x_0$ unless, of course, its dimension be equal to the co-rank of $\mathcal{P}.$ When this is not the case, there are many possible choices for the chains. Since the element $E_{\rho_k}$ is maximal, relative to the sub-chain terminating in $E_{\rho_k-1}~,$ this element is unique, always relative to the sub-chain, though the above indicated integral sub-manifold needs not be locally unique, relative to the complete chain. Furthermore, choosing the local vector fields $\xi_{\mu}$ in such a way that the vectors $\{\xi_{\mu}(x_0),~\mu\leq j\}$ generate $E_j~,$ we can obtain, by successive integrations, an ascending chain of integral sub-manifolds respectively tangent to the elements of the linear chain, each one contained in the next manifold. We finally observe that the singularities as well as the singular integral sub-manifolds can be located inasmuch as evaluated by computing the integer $\rho_k$ at the various points of the manifold \textit{M} and corresponding to the various ascending chains of integral elements.

\section{The derived systems}
Given a Pfaffian system $\mathcal{P}~,$ we recall that the \textit{characteristic system} associated to $\mathcal{P}$ is the Pfaffian system generated by all the differential 1-forms $\omega$ and $i(\xi)d\omega$ where $\xi$ is an arbitrary local vector field annihilated by $\mathcal{P}$ and $\omega$ an arbitrary local section of the given system. The contra-variant counterpart namely, the annihilator of the characteristic system is defined as follows: At each point $x\in M~,$ we consider the sub-space, of the tangent space, whose elements are the vectors \textit{v} such that $i(v)\omega=0$ and $i(w)i(v)d\omega=0$ for any vector \textit{w} annihilated by $\mathcal{P}$ and any local $1-$form $\omega$ belonging to the system. In other terms, the above \textit{characteristic} sub-space is defined as the set of all the vectors \textit{v} annihilated by $\mathcal{P}$ that further verify $d\omega(v,w)=0$ for any vector \textit{w} also annihilated by $\mathcal{P}$ and any local section $\omega$ of $\mathcal{P}.$ We denote by $\mathcal{CH}$ this characteristic system and observe that its purpose resides in (loosely speaking) finding the minimum number of independent function by means of which the initially given system can be re-written or, equivalently, in finding the local fibration with the largest dimensional fibres with respect to which the given system factors to the quotient \textit{i.e.,} to the base space of this fibration onto an equivalent \textit{quotient} system (\textit{cf.} \cite{Kumpera1982}). Inasmuch, the \textit{derived system} of $\mathcal{P},$ denoted by $\mathcal{P}_1~,$ is the Pfaffian system generated by the forms $\omega~,$ local sections of $\mathcal{P},$ for which $d\omega\equiv 0\hspace{2 mm}mod~\mathcal{P}$ \textit{i.e.}, generated by those local forms that satisfy the integrability condition with respect to $\mathcal{P}.$ The derived system is, of course, a Pfaffian sub-system of the initially given system. Iterating this procedure, we obtain a descending chain of Pfaffian systems terminating by an integrable system (its derived system coincides with the system) that, eventually, reduces to the null system. It is obvious that we can state \textit{contra-variant} definitions for both the characteristic as well as the derived systems. However, such definitions are rather clumsy, lacking all the subtle elegance of Cartan's definitions. 

\vspace{2 mm}
We can now apply the above inquiry, concerning integrable elements, to both characteristic and derived systems though, presently, we shall only highlight the most relevant facts and, among them, those facts that were placed in evidence by Cartan (pg.9, n.6 and the following pages). We first observe, by simply inspecting the contra-variant definition, that the characteristic system is equal, at each point, to the union of those  $1-$integral elements (or tangent vectors) of $\mathcal{P}$ that are in involution with all the other $1-$integral elements of $\mathcal{P}.$ The characteristic system is always integrable (involutive) independently of the nature of $\mathcal{P}$ and, if moreover it is regular, we can factor locally its integral foliation. Furthermore, we can also factor all the differential forms of $\mathcal{P}$ \textit{modulo} the leaves and thereafter, the initially given system factors as well to an equivalent quotient system  $\widetilde{\mathcal{P}}.$ We infer that the integral manifolds of the initial system are the pullbacks (inverse images) of the integral manifolds in the quotient and observe that the quotient system has null characteristics.

\vspace{2 mm}
We next give a glance at the derived system. From a \textit{contra-variant} point of view, the contact element determined by $\mathcal{P}_1$ at the point $x_0$ is the linear subspace generated by all the vectors annihilated by $\mathcal{P}$ together with all the brackets $[\xi,\eta]_{x_0}$ where $\xi$ and $\eta$ are two local vector fields annihilated by $\mathcal{P}.$ In the \textit{contra-variant} mode, the chain of successive derived spaces or systems is increasing whereas, \textit{co-variantly}, it decreases.

\vspace{2 mm}
\noindent
Obviously,
\begin{equation*}
\mathcal{P}_1\subset\mathcal{P}\subset\mathcal{CH}
\end{equation*}
hence, in what concerns the corresponding pseudo-groups of local automorphisms and taking into account the functoriality of the constructions, we obtain the following two inclusions: 
\begin{equation}
\Gamma_{\mathcal{P}}\subset\Gamma_{\mathcal{CH}}\hspace{10 mm}and\hspace{10 mm}\Gamma_{\mathcal{P}}\subset\Gamma_{\mathcal{P}_1}~.
\end{equation}
We can then proceed to extend the above second inclusion so as to obtain a composition series and, eventually, a \textit{Jordan-Hölder resolution}. Making use of the previous results as well as those in \cite{Kumpera2016}, we might then simplify the integration process of $\mathcal{P}.$

\section{The pseudo-groups of local automorphisms}
The hypotheses and the notations being as above, we shall now have a closer look into the integration process of the Pfaffian system $\mathcal{P}.$ According to Lie and Cartan, we first have to choose convenient intermediate pseudo-groups, in the present case these will be associated to structures, in such a way that the second inclusion in (23) extends hopefully to a Jordan-Hölder resolution. As a first step, we examine the local automorphisms of the successive derived Pfaffian systems and compare these pseudo-groups with each other as well as with the automorphisms of the initial system $\mathcal{P}.$

\vspace{2 mm}
Let $\mathcal{P}$ be a Pfaffian system defined on the $n-$dimensional manifold \textit{M} and let us assume, at least for the time being, that it is regular \textit{i.e.}, that $rank~\mathcal{P}_x=rank~\mathcal{P}_y$ for any two arbitrary points $x,y\in M$ where, in the \textit{co-variant} jargon, the \textit{rank} is equal to the \textit{dimension}. We then consider a local diffeomorphism (referred to as a local transformation)
\begin{equation*}
\varphi:\mathcal{U}~\longrightarrow~\mathcal{V}~,
\end{equation*}
where $\mathcal{U}$ and $\mathcal{V}$ are open subsets in \textit{M} and the transformation $\varphi$ is assumed to be differentiable. To say that $\varphi$ preserves or leaves invariant the Pfaffian system $\mathcal{P}$ means, of course, that $\varphi^*(\mathcal{P}_y)=\mathcal{P}_x$ as soon as $\varphi(x)=y~.$ We shall also refer to $\varphi$ as being an \textit{automorphism} of $\mathcal{P}.$ Observing that the linear automorphism
\begin{equation*}
\varphi^*_x:\mathcal{P}_y~\longrightarrow~\mathcal{P}_x
\end{equation*}
is entirely determined by its $1-$jet $j_1\varphi(x)$ (in fact, it is essentially the same thing), we infer that any local transformation $\varphi$ for which $j_1\varphi(x),$ $x\in\mathcal{U}~,$ transposes to an automorphism pertaining (restricting) to $L(\mathcal{P}_y,\mathcal{P}_x),$ belongs to $\Gamma_{\mathcal{P}}$ and consequently, this pseudo-group is a Lie pseudo-group of order 1.\footnote{The groupoid of first order jets is a differentiable manifold, the source and target maps are submersions (more precisely, surmersions) and the composition operation as well as the passage to inverses are both differentiable mappings, the latter being a diffeomorphism.}

\vspace{2 mm}
Our next task is to show that the Lie groupoid composed by all the $1-$jets of the elements belonging to $\Gamma_{\mathcal{P}}$ is an involutive first order partial differential equation in the sense of Élie Cartan (\cite{Kumpera1972}). Let us denote by $\textbf{G}^1_{\mathcal{P}}$ this first order groupoid associated to $\Gamma_{\mathcal{P}}$ and let us recall that the \textit{standard prolongation} of the differential equation $\textbf{G}^1_{\mathcal{P}}$ is the set of all second order jets $j_2\psi(x)~,$ of local maps $\psi$ of \textit{M}, such that the \textit{Ehresmann flow}
\begin{equation*}
y~\longmapsto~j_1\psi(y)
\end{equation*}
is tangent, at first order and at the point \textit{x}, to the equation $\textbf{G}^1_{\mathcal{P}}.$ This means, in denoting by $\widetilde{\psi}$ the above flow, that the image sub-manifold of $\widetilde{\psi}$ is tangent to $\textbf{G}^1_{\mathcal{P}}$ at the point $\widetilde{\psi}(x)$ ($\widetilde{\psi}_*(T_xM)\subset T_{\widetilde{\psi}(x)}\textbf{G}^1_{\mathcal{P}}$). We can define inasmuch the standard prolongations of the higher order jet groupoids associated to the pseudo-group $\Gamma_{\mathcal{P}}$ and a straightforward argument shows the following lemma:

\vspace{2 mm}
\newtheorem{jhil}[LemmaCounter]{Lemma}
\begin{jhil}
All the $k-$th order groupoids, $k>1~,$ associated to the pseudo-group $\Gamma_{\mathcal{P}}$ are the successive standard prolongations of the first order groupoid $\textbf{G}^1_{\mathcal{P}}.$
\end{jhil}

\vspace{2 mm}
\noindent
More important, we want to show that the first order groupoid as well as all the higher order groupoids are \textit{involutive equations} in the sense of Cartan. The above Lemma being a first step, we now have to show that all the linear \textit{symbols} of these equations are involutive \textit{i.e.,} that they are $2-acyclic.$ A rather long juggling with the techniques found in \cite{Kumpera1972}  including the application, to the previous equations, of the \textit{Spencer linear complex} will eventually prove the following result.\footnote{We do not provide the details of the proof since it is rather long and the result, in itself, is just auxiliary, not of first order importance.}

\vspace{2 mm}
\newtheorem{klm}[LemmaCounter]{Lemma}
\begin{klm}
All the groupoids associated to the Lie pseudo-group $\Gamma_{\mathcal{P}}$ are involutive.
\end{klm}

\vspace{2 mm}
\noindent
Let us next consider the derived sequence associated to $\mathcal{P}.$ It is a descending chain
\begin{equation*}
\mathcal{P}\supset\mathcal{P}_1\supset~\cdots~\supset\mathcal{P}_{\mu}=\mathcal{P}_{\mu+1}~,
\end{equation*}
the last term being integrable or eventually null and, by the naturality of the constructions, we conclude that
\begin{equation*}
\Gamma(\mathcal{P}_{\mu})\supset~\cdots~\supset\Gamma(\mathcal{P}_1)\supset\Gamma(\mathcal{P})
\end{equation*}
hence
\begin{equation}
\textbf{G}^k_{\mathcal{P}_{\mu}}\supset~\cdots~\supset\textbf{G}^k_{\mathcal{P}_1}\supset\textbf{G}^k_{\mathcal{P}}~,
\end{equation}
for all integers $k\geq 1~.$ In general, there is no reason for the previous sequence to be a composition series (each term being normal in the preceding one) though we shall exhibit a rather outstanding situation where not only this is true but, furthermore, the series becomes \textit{Jordan-Hölder}. 

\vspace{2 mm}
Let us now assume that the Pfaffian system $\mathcal{P}$ is a \textit{flag system} \textit{i.e.,} $rank~\mathcal{P}_i=rank~\mathcal{P}_{i+1}+1$ for all \textit{i}. Then, with the notations of the reference \cite{Kumpera2014}, the finite and infinitesimal automorphisms of any of the factored derived systems  $\overline{\mathcal{P}}_{\mu}$ (\textit{modulo} their characteristic variables) is canonically equivalent (isomorphic), via a natural \textit{merihedric} prolongation algorithm, to those of $\mathcal{P}$ and consequently all the corresponding higher order groupoids are canonically isomorphic. However, when these derived systems are considered on their \textit{ambient} space \textit{M}, new automorphisms do appear namely, those of the associated characteristic systems that are also automorphisms of the derived systems in evidence. However, a characteristic system being always integrable, its automorphisms are, locally, those that leave invariant (\textit{i.e.,} permute) the fibres of the fibration obtained by integrating the characteristics and, still more locally, those that leave invariant (permute) the parallel spaces to one of the components in a product. In terms of local coordinates, we consider all the local transformations in an $n-$space, a \textit{simple} pseudo-group according to Cartan, and just retain those transformations that maintain a certain number of coordinates depending only upon these restricted  coordinates, the resulting pseudo-group being \textit{simple} as well. The quotients of these pseudo-groups and, inasmuch, the quotients of the groupoids are therefore simple and the sequence
\begin{equation*}
\textbf{G}^k_{\mathcal{P}_{\mu}}\supset~\cdots~\supset\textbf{G}^k_{\mathcal{P}_1}\supset\textbf{G}^k_{\mathcal{P}}
\end{equation*}
is Jordan-Hölder. The same argumentation holds on the infinitesimal level and the integration  of the system $\mathcal{P}$ will run as follows: Each quotient manifold is $1-$dimensional and the quotient Pfaffian system is invariant by the entire first order groupoid of jets hence is equal to the total co-tangent bundle (with $1-$dimensional fibre) whereafter the connected integral manifolds are the points. Taking successively the resolvent systems, we end up concluding, as already expected, that the maximal integral manifolds of the Flag System  $\mathcal{P}$ are just its integral curves.

\vspace{2 mm}
\noindent
Let $E_{\rho_k}$ be a maximal integral contact element of $\mathcal{P}$ at the point $x_0$ and $\mathcal{I}_{\rho_k}$ a local integral sub-manifold of maximal dimension, tangent to $E_{\rho_k}~.$ Then, of course, $\mathcal{I}_{\rho_k}$ is also an integral sub-manifold of the derived system, though not forcibly maximal, and will be contained in maximal integral manifolds (not all) of this derived system.

\section{Pfaffian systems whose characters are equal to one}
We now assume that $\rho_k=n-r-1~,$ where we recall that $n=dim~M$ and $r=rank~\mathcal{P}.$ In other terms, the dimension of $E_{\rho_k}$ and consequently also that of $\mathcal{I}_{\rho_k}$ are one unit less than the dimension of $\Sigma_{x_0},$ the contact element at the point $x_0$ annihilated by $\mathcal{P}_{x_0}.$ Let us show that $rank~\mathcal{P}_1=r-1$ and provide a criterion for the integrability of this derived system. Cartan and, earlier, von Weber already gave such a criterion that Cartan claimed to be \textit{un fait très remarquable}.  

\vspace{2 mm}
Firstly, since $\mathcal{P}$ is not integrable (involutive) on account of the maximality of $E_{\rho_k},$ it follows that $rank~\mathcal{P}_1\leq r-1~.$ On the other hand, since the dimension of $\mathcal{I}_{\rho_k}$ is equal to $n-r-1$ and all the forms $\omega^i$ vanish on this integral manifold, we can assume relabeling if necessary by $\{\varphi^i\}$ the local basis of $\mathcal{P},$ that $\mathcal{P}_1$ is generated by $\{\varphi^2,~\cdots~,\varphi^r\}~.$ In fact, choosing a local section $\varphi_1$ of $\mathcal{P}$ that verifies $d\varphi_1|\Sigma_{x_0}\neq 0~,$ we can argue, together with Cartan, without performing the calculations.

\vspace{2 mm}
\noindent
Secondly and returning to the consideration of the \textit{integral} contact elements, we consider at the point $x_0$ a maximal integral contact element $E_{\rho_k}.$ As evidenced previously, this element is equal to the tangent space, at $x_0~,$ of an integral sub-manifold \textit{W} of the system $\mathcal{P}.$ Locally, this sub-manifold is integral with respect to the forms $\omega^i$ of a local basis of $\mathcal{P}$ and some additional independent forms $\varphi^s$ that restrict to a single form $\varphi$ in the case when the \textit{character} of the system is equal to one. But then, for the Pfaffian system $\mathcal{Q}$ generated by these forms, we can write
\begin{equation*}
d\omega^i\equiv c^i\mu^i\wedge\varphi\hspace{5 mm}mod~\mathcal{P},
\end{equation*}
with certain functions $c^i.$ Consequently, the forms belonging to the sub-module of $\Gamma(\mathcal{P})$ defined by the equation $~\sum~c^i\omega^i=0~$ are precisely those that belong to $\mathcal{P}_1.$ This system is under certain conditions integrable namely, when its \textit{gender} is larger than 1. This invariant, already introduced by von Weber, has to do with the "number" of differential $1-$forms independent from those forming $\mathcal{P}$ and that enter in the expressions of $d\omega\hspace{2 mm}mod~\mathcal{P},~\omega\in\Gamma{\mathcal{P}}$ (local sections). We shall consider this invariant only further but define it right below. Taking a maximal integral manifold of $\mathcal{P}_1,$ when this system is integrable, that contains $x_0~,$ we can restrict all the data to this manifold, the investigation of the properties of $\mathcal{P}$ becoming considerably facilitated. We shall also examine, later, those Pfaffian systems with characters equal to two or more but presently let us exhibit Cartan's \textit{co-variant} version of the former discussion.

\vspace{2 mm}
\newtheorem{abc}[DefinitionCounter]{Definition}
\begin{abc}
The \textbf{gender} (genre) of a local exterior differential form $\Omega~,$ with respect to a Pfaffian system $\mathcal{P},$ is the smallest integer \textbf{h} such that $\Omega^{h+1}\equiv 0\hspace{2 mm}mod~\mathcal{P}$ (the exponent refers to wedge products). The \textbf{gender} of a family of local exterior differential forms is the maximum value of the gender of its elements. 
\end{abc}

\vspace{2 mm}
When $\mathcal{P}$ is locally generated by $\{\omega^1,~\cdots~,\omega^r\}~,$ then the above condition can be restated by the equality $\omega^1\wedge~\cdots~\wedge\omega^r\wedge\Omega^{h+1}=0$ for the exterior forms defined on open sets contained in the domain of the generators $\omega^i.$ The \textbf{gender} of a Pfaffian system $\mathcal{P}$ is equal to the gender of the family $d\omega$ where $\omega$ is an arbitrary local section of $\mathcal{P}$ and, of course, it is given by the gender of any system of generators. We also observe that the system $\mathcal{P}$ is integrable if and only if it is of gender zero (\cite{Cartan1901}).

\section{The Cartan co-variant approach}
We start by considering the linear $1-$dimensional elements belonging to $\mathcal{P}$ namely, all the $1-$dimensional sub-spaces of $\mathcal{P}$ at the point $x_0$ and inquire which should be considered as 1-dimensional \textit{integral elements}. Let $\omega$ be a generator of such a linear element and let us imagine that its annihilator $ker~\omega_{x_0}$ is an integral element for the Pfaffian system generated by $\omega~.$ Then the condition reads $(d\omega)_{x_0}\equiv 0\hspace{2 mm}mod~\omega_{x_0}.$ It seems however more adequate, for reasons that shall be clarified hereafter, 
to assume the stricter requirement for the integral elements namely,
\begin{equation*}
(d\omega)_{x_0}\equiv 0\hspace{2 mm}mod~\mathcal{P}.
\end{equation*}
In terms of local generators, the above condition reads
\begin{equation*}
(d\omega)_{x_0}\equiv 0\hspace{2 mm}mod~\{\omega^1_{x_0},~\cdots~,\omega^r_{x_0}\}.
\end{equation*}

\noindent
Next, we inquire which should be considered as being the $2-$dimensional integral elements containing a given (fixed) $1-$dimensional element. Though entering in conflict with the notations adopted previously in the \textit{contra-variant} discussion, we shall keep these notations in order to conform to Cartan's writing. We choose a given (fixed) though arbitrary linear integral element $E_1$ and consider the sub-space of all those linear integral elements that are \textit{in involution} with the given element. Inasmuch as in the contra-variant setting, not all the linear integral elements are suitable and the condition for this to be so is, of course, the integrability condition involving the two elements. More precisely, the co-vector $\omega_{x_0}~,$ or the sub-space generated by it, is in involution with $E_1$ when, by definition,
\begin{equation*}
\omega_{x_0}\equiv 0\hspace{3 mm}mod~\mathcal{P},\hspace{3 mm}(d\omega)_{x_0}\equiv 0\hspace{3 mm}mod~\mathcal{P}\hspace{3 mm}and 
\end{equation*}
\begin{equation*}
d\omega_0\equiv 0\hspace{3 mm}mod~\mathcal{P}~, 
\end{equation*}
where $\omega_0$ is chosen so as to induce a non vanishing co-vector belonging to $E_1~.$ In other words, this simply means that $\omega\in\mathcal{P}$ and that
\begin{equation*}
d\omega\wedge\varphi^1\wedge~\cdots~\wedge\varphi^r=0~,
\end{equation*}
the same holding for $\omega_0~.$ We can now choose $E_2$ containing $E_1~,$ only composed of co-vectors that are in involution with $E_1$ and consider the space of all those co-vectors, belonging to $\mathcal{P},$ that are simultaneously in involution with all the elements of $E_2~.$ Continuing this process, we choose an element $E_3~,$ containing $E_2~,$ only composed by co-vectors that are in involution with all the elements of $E_2$ and so forth. The $character~s_1$ defined by Cartan as being the integer satisfying $n-r-s_1=dim~\widetilde{E}_1,$ where the latter is the space of all the integral elements in involution with $E_1,$ is then equal to $\rho_k$ in the contra-variant setting and both situations, the co-variant and the contra-variant, are mirror images one of the other with respect to $\mathcal{P}.$ Furthermore, the result of the previous section concerning the derived system has, in this context, a very simple proof. In the next sections, we discuss Pfaffian systems with characters larger than one that provide a wide spectrum of different situations. We shall nevertheless keep within the Cartan co-variant setting since, as already mentioned earlier, the discussion as well as the calculations become straightforward and simple though we stop understanding anything since we are rowing against the tide.

\vspace{2 mm}
\noindent
From what was shown in the previous section, we can state the following lemma. 

\vspace{2 mm}
\newtheorem{kiss}[LemmaCounter]{Lemma}
\begin{kiss}
A necessary and sufficient condition that the character of a Pfaffian system $\mathcal{P}$ be equal to one is that the rank of its derived system $\mathcal{P}_1$ be one unit less than its own rank.
\end{kiss}

\section{Integration of Pfaffian systems with character equal to one}
We start by taking a local basis $\{\omega^i\}$ for the Pfaffian system $\mathcal{P}$ such that $\{\omega^2,~\cdots~,\omega^r\}$ generates its associated derived system $\mathcal{P}_1,$ complete this basis to a local co-frame by adding the 1-forms $\{\overline{\omega}^j\}$ and continue to assume that the characteristics of $\mathcal{P}$ are null.

\vspace{2 mm}
Let us first prove the result that Cartan considered so much surprising, recalling that \textbf{h} denotes the \textit{gender} of a Pfaffian system.

\vspace{2 mm}
\newtheorem{blss}[LemmaCounter]{Lemma}
\begin{blss}
A necessary and sufficient condition that the first derived system $\mathcal{P}_1$ of a Pfaffian system $\mathcal{P},$ of character $1~,$ be integrable is that $\textbf{h}\geq 2~.$
\end{blss}

\vspace{2 mm}
\noindent
Adapting the generators of $\mathcal{P}$ in such a way that the forms $\{\omega^i\},~i\geq 2,$ generate $\mathcal{P}_1~,$ we have by definition
\begin{equation}
d\omega^i\equiv 0\hspace{3 mm}mod~\{\omega^1,~\cdots~,\omega^r\},\hspace{3 mm}i\geq 2~,
\end{equation}
and write, to begin,
\begin{equation}
d\omega^2\equiv\chi\wedge\omega^1\hspace{3 mm}mod~\{\omega^2,~\cdots~,\omega^r\}~,
\end{equation}
where $\chi$ can be taken to be a linear combination of the forms $\overline{\omega}^j~.$ A second differentiation then yields
\begin{equation*}
d\chi\wedge\omega^1+\chi\wedge d\omega^1\equiv 0\hspace{3 mm}mod~\{\omega^2,~\cdots~,\omega^r\}~,    
\end{equation*}
hence
\begin{equation*}
\chi\wedge d\omega^1\equiv 0\hspace{3 mm}mod~\{\omega^1,~\cdots~,\omega^r\}~,    
\end{equation*}
on account of (25). We next assume that $\chi$ and $\omega^1$ are independent. Since $\chi$ only contain terms in $\overline{\omega}^j,$ we can write $\chi\wedge d\omega^1=\mu+\eta$ where each term of $\mu$ contains some element $\omega^i$ but $\eta$ only contains terms that are expressed by means of the $\overline{\omega}^j.$ Moreover, the term $\eta$ is equal to the product $\chi\wedge\overline{\eta}~,$ where $\overline{\eta}$ is the sum of all the terms in $d\omega^1$ that are written only with the help of the forms $\overline{\omega}^j.$ It then follows that the exterior product $d\omega^1\wedge d\omega^1$ presents, in its unique term not containing any form $\omega^i,$ a double product $\chi\wedge\chi$ hence this term vanishes and the double product above becomes equal to:
\begin{equation*}
d\omega^1\wedge d\omega^1\equiv 0\hspace{3 mm}mod~\{\omega^1,~\cdots~,\omega^r\}~.   
\end{equation*}
Consequently the gender of $\mathcal{P}$ is at most equal to 1 contradicting the initial hypothesis. We infer that $\chi$ and $\omega^1$ are dependent and that the relation (26) reduces consequently to
\begin{equation*}
d\omega^2\equiv 0\hspace{3 mm}mod~\{\omega^2,~\cdots~,\omega^r\}~.
\end{equation*}
The same argument being valid for all the other forms generating $\mathcal{P}_1,$ the integrability of this derived system then follows. 

\vspace{2 mm}
We can now proceed with the integration of the system $\mathcal{P}$ namely, construct integral manifolds of dimension $n-r-1$ that are tangent to given maximal integral contact elements, also of the same dimension, since the character of $\mathcal{P}$ is assumed to be equal to one. Fixing the point $x_0~,$ we denote by $\mathcal{L}$ the integral leaf of $\mathcal{P}_1$ that contains the given point. Then $dim~\mathcal{L}=n-r+1$ and, if we assume as previously that the forms  $\omega^i,~i\geq 2~,$ generate $\mathcal{P}_1,$ the restricted system $\mathcal{P}|\mathcal{L}$ is generated by the single 1-form $\omega^1|\mathcal{L}~.$ Obviously, any maximal integral element issued at the point $x_0$ is contained in $T_{x_0}\mathcal{L}$ and any maximal integral manifold of $\mathcal{P},$ containing $x_0~,$ is contained in $\mathcal{L}~,$ being as well an integral manifold of the restricted system. Moreover, the maximality properties are preserved under restriction since the restricted system is not integrable, otherwise there would be an $(n-r)-$dimensional integral manifold of $\mathcal{P}$ passing through $x_0~.$ Since the \textit{gender} of a Pfaffian system is a point-wise notion involving only the exterior algebra of the co-tangent spaces at given points, we infer that the gender of the restricted system is the same as that of the initial system. On the other hand, the restricted system $\mathcal{P}|\mathcal{L}$ is not integrable and has rank equal to one hence is a Darboux system (\cite{Kumpera2014}). The definition \textit{per se} of the gender of this restricted system will then assert that $\mathcal{P}|\mathcal{L}$ is a Darboux system with Cartan or, inasmuch, Darboux class equal to  $2\textbf{h}+1~,$ hence generated locally by a form
\begin{equation*}
dz^{h+1}+p_1dz^1+p_2dz^2+~\cdots~+p_hdz^h~,    
\end{equation*}
where $\it{h}=\textbf{h}$ is the gender at a point\footnote{Cartan preferred minus signs.}. Putting all the stacks together, we infer that $\mathcal{P}$ admits local models of the form:
\begin{equation*}
dy^1\hspace{64 mm}
\end{equation*}
\begin{equation*}
dy^2\hspace{64 mm}
\end{equation*}
\begin{equation*}
..........\hspace{63 mm}
\end{equation*}
\begin{equation*}
dy^{r-1}\hspace{59 mm}
\end{equation*}
\begin{equation*}
dz^{h+1}+p_1dz^1+p_2dz^2+~\cdots~+p_hdz^h
\end{equation*}
with respect to suitable local coordinates and where the functions $y^1,y^2,~\cdots~,y^{r-1}$ are chosen to be independent first integrals of the (integrable) derived system (\textit{cf.} \cite{Cartan1901}, pg.28). The solutions of the above system will provide, locally, the maximal integral manifolds of $\mathcal{P}.$ If by any chance the given Pfaffian system admits non-vanishing characteristics\footnote{Cartan never imposed any restriction on the characteristics since, in his \textit{modus operandi}, nothing much would change.}, then it suffices to chose, for $y,~z$ and $p,$ first integrals of the characteristic system.

\vspace{2 mm}
Let us finally examine what happens when $\textbf{h}=1~.$ We can in this case integrate, starting from a given point $x_0$ on $\mathcal{L},$ all the vector fields belonging to the annihilator of $\mathcal{P}|\mathcal{L}~.$ Denoting by 
\begin{equation*}
t~\longmapsto~(z^1(t),~\cdots~,z^{\mu}(t))~,\hspace{5 mm}\mu=n-r+1~,
\end{equation*}
a generic integral curve commencing at $x_0~,$ we plot all the points obtained by considering the images of the above integral curves conditioned to a given specific, though arbitrary, non-trivial relation among the targets $z^i.$ A detailed proof of this statement is rather long and shall be omitted. Fortunately, the author is backed up by Cartan's fantastic intuition (\cite{Cartan1901}, pg.28-29, where he refers to von Weber \cite{Weber1898}).

\section{Pfaffian systems whose characters are larger than one}
We begin by examining systems with characters equal to 2 and, to simplify, we assume right away that their characteristics are null. Much will then depend upon the \textit{gender} of the system and offers several options. We exhibit firstly conditions under which the character being equal to 2 implies that the rank of the first derived system is equal to $r-2~,$ the general rule not being as strict in the present case as it was previously. It seems furthermore worthwhile, at this point, to make a few comments since our statement is, apparently, in disagreement with Cartan's claims. Assume, for a moment, that the Pfaffian system $\mathcal{P}$ has non-trivial characteristics and denote by $\mu$ the difference in the dimensions of the initially given manifold \textit{M} and the (local) quotient manifold $\overline{M}~,$ modulo the characteristics of the system $\mathcal{P}.$ Then $\mathcal{P}$ factors to an equivalent system $\overline{\mathcal{P}},$ both systems have the same rank and, consequently, the annihilator distributions $\Sigma$ and $\overline{\Sigma}$ differ, in their point-wise dimensions, by $\mu,$ the dimension of $\overline{\Sigma}$ being lesser since the characteristic distribution is contained in $\Sigma$. Similar statements hold for the derived system and its quotient in $\overline{M}$ though the latter can still have non-trivial characteristics thus explaining the discrepancy concerning the rank of $\overline{\mathcal{P}_1}$ and the dimensions of the maximal integral elements of $\mathcal{P}$ and $\overline{\mathcal{P}}.$ We also observe that the passage to the derived systems is functorial and therefore compatible with quotients \textit{i.e.}, $\overline{\mathcal{P}_1}=(\overline{\mathcal{P}})_1~.$

\vspace{3 mm}
In the sequel we shall only give a brief account on Cartan's results so as to promptly continue with our main discussion. In what concerns the Pfaffian systems of character 2, the main result is the following:

\vspace{2 mm}
\newtheorem{buss}[TheoremCounter]{Theorem (Cartan)}
\begin{buss}
The rank of the derived system of a Pfaffian system $\mathcal{P}$ of character two is always two units less than the rank of the system except when $\mathcal{P}$ admits characteristic elements with dimensions not less than  $dim~\Sigma_{x_0}-3$ in which case the rank lowers by at least three units.
\end{buss}

\vspace{2 mm}
\noindent
The proof is similar to that of the previous lemma as soon as we choose two independent local sections of the system whose differentials do not vanish on the kernel and an entirely analogous result also holds in general. We shall therefore only give some attention to the second part of the statement that is far more delicate. Let us then assume, the Pfaffian system being of character 2, that the rank of the derived system $\mathcal{P}_1$ is equal to $r-\mu$ where $\mu\geq 3~.$ This being so, we can choose a local basis $\{\omega^1,~\cdots~,\omega^{\mu},~\cdots~,\omega^r\}$ for $\mathcal{P}$ and in a neighborhood of the point $x_0$ in such a way that the forms $\{\omega^1,~\cdots~,\omega^{\mu}\}$ do not belong to $\mathcal{P}_1~,$ hence inasmuch any non-trivial linear combination of these forms cannot belong to the derived system, the remaining forms generating $\mathcal{P}_1~.$ This being so and since $d\omega^j\not\equiv 0\hspace{3 mm}mod~\mathcal{P},~j\leq\mu,$ we can select $\mu$ independent vectors $v_j$ in $\Sigma_{x_0}$ not belonging, of course, to the characteristic sub-space of $\mathcal{P}$ at the point $x_0$ in such a way that the family
\begin{equation*}
\{\omega^1,~\cdots~,\omega^r,i(v_1)d\omega^1,~\cdots~,i(v_{\mu})d\omega^{\mu}\}
\end{equation*}
is a local basis for the characteristic system of $\mathcal{P}.$ The rank $r+\mu$ of this characteristic system must be equal to \textit{n}, hence $\mu =n-r~,$ otherwise the system $\mathcal{P}$ would have non-trivial  characteristics contradicting our assumption. But then, the system $\mathcal{P}$ will be of rank equal to the dimension of \textit{M}, its annihilator is the null distribution and the integral leaves are the points of \textit{M}, a rather uninteresting foliation. Within the scope of our assumptions, we can therefore restate the above theorem as follows:

\vspace{2 mm}
\newtheorem{boss}[TheoremCounter]{Theorem}
\begin{boss}
Let $\mathcal{P}$ be a Pfaffian system with null characteristics and character equal to 2. Then the rank of its  derived system is two units less than the rank of the system.
\end{boss}

\vspace{2 mm}
\noindent
Let us now say a word on the number 3. Apart from being an extraordinary prime number that brings much luck, let us argue as follows:

\vspace{2 mm}
\noindent
We select among all the possible non-trivial linear combinations of the forms $\{\omega^1,~\cdots~,\omega^{\mu}\}$ one such that has a minimum gender at the point $x_0~.$ Re-arranging once more the local basis, we can assume for convenience that the above form is simply $\omega^1$ and can therefore write, at the point $x_0~,$ $(d\omega^1)^{\textbf{h}+1}\equiv 0\hspace{3 mm}mod~\{\omega^1,~\cdots~,\omega^r\}~,$ $\textbf{h}$ being the minimum value for which the above congruence holds. A simple calculation will then show that the value $\textbf{h}=1$ is the only possible value for which the rank of the characteristic system is  larger than $r-2$ and consequently the rank of the derived system is lesser than $r-2~,$ precisely equal to $r-3$. 

\vspace{2 mm}
\noindent
However, we should mention that, in many applications, Pfaffian systems with non-trivial characteristics do appear most naturally and, furthermore, the factoring of a Pfaffian system, modulo its characteristics, is never an easy task since it requires an integration process that only reduced locally to a system of ordinary differential equations when the characteristic leaves are $1-$dimensional.

\vspace{2 mm}
Cartan also discusses those Pfaffian systems for which half of the point-wise dimension of their annihilators $\Sigma$ is not less than the value of their characters. Since, once more, the discussion of this setup involves non-vanishing characteristics, we shall leave it aside. There is however a special case when the equality of the above mentioned values holds and that merits to be commented. Let us then assume that the previously mentioned system $\mathcal{P}$ also verifies the following properties:

\vspace{3 mm}
(a) The character of $\mathcal{P}$ is equal to $2~,$
 
\vspace{2 mm}
(b) $\mathcal{P}$ has null characteristics and 
 
\vspace{2 mm}
(c) All the maximal integral contact elements that contain a specific $1-$dimensional integral element $E_1$ also contain an integral element $E_2\supset E_1$ whose dimension is not less than $2~.$

\vspace{3 mm}
\noindent
Under these conditions, the equality $\frac{n-r}{2}=\rho_k$ holds. Translated in terms of maximal integral manifolds, the condition (c) means that if two such manifolds have in common a line they also have in common at least a surface. Cartan calls such systems \textit{systatical (systatiques)}.

\section{Singular Pfaffian Systems}
For generic Pfaffian systems, the successive characters assume their maximum values hence, in particular, for systems with character equal to 2 the values are:
\begin{equation*}
s_1~=~s_2~=~\cdots~=~s_{\rho_k-1}~=~2~.    
\end{equation*}
As for $s_{\rho_k}~,$ it is equal to zero, one or two according to the remainder of the division of $n-r$ by $3~.$ This being so, $s_{\rho_k}=0$ when $n-r=3\rho_k-2~,$ $s_{\rho_k}=1$ when $n-r=3\rho_k-1~,$ and $s_{\rho_k}=2$ when $n-r=3\rho_k~.$ 

\vspace{2 mm}
\noindent
We shall say, together with Cartan, that the Pfaffian system $\mathcal{P},$ with character equal to $2~,$ is \textit{singular} when $s_{\rho_k-1}\leq 1~.$ Still quoting Cartan, for a singular system having character two and null characteristics the maximum dimension of the linear integral contact elements \textit{i.e.}, the value of $\rho_k$ is at least equal to 3 and, consequently, the point-wise dimension of the distribution $\Sigma$ is at least equal to 6 . Another rather surprising property, consequence of the above assertions, is the following:

\vspace{2 mm}
\newtheorem{bosta}[PropositionCounter]{Proposition}
\begin{bosta}
All the systatical Pfaffian systems are singular as soon as $\rho_k\geq 3~.$
\end{bosta}

\vspace{2 mm}
Cartan's \textit{Mémoire} then terminates with an attempt to extend some of the above properties to non-systatical systems and culminates in the following statement:

\vspace{2 mm}
\noindent
\textit{A singular Pfaffian system with null Cauchy characteristics and character equal to 2 can be integrated by means of a system of ordinary differential equations. In case the system admits non-trivial characteristics, its characteristic system can also be integrated by means of ordinary differential equations whenever the derived system of the given system is non-integrable.}

\vspace{2 mm}
\noindent
Curiously enough, the second order \textit{Monge characteristics} play an important role in the proofs.

\vspace{3 mm}
We finally arrive at what really matters namely, the integration of the above considered Pfaffian systems via the reductions provided by Jordan-Hölder sequences. The choice of the examples that follow is far from being casual or fortuitous since they are all backed up by relevant geometrical problems. In the previous discussion, we have seen to what extent the derived systems are relevant. Consequently, we begin by integrating the derived system and afterwards locate, among the manifolds thus obtained, those that are maximal for the given system. We are thus faced with constructing Jordan-Hölder resolutions for the derived system and each specific situation will exhibit its particular techniques. It might however be of some advantage to try initially integrating one of the systems that make part of the sequence of iterated derived systems and, in this case, the chain (43) will already provide some useful terms. In particular, it is possible to detect the appearance of ordinary differential equations, these showing up when the quotients are $1-$dimensional. \nocite{Kumpera1999}

\section{Examples}
Allowing initially $dim~M=5$ and $rank~\mathcal{P}=3~,$ the system generated by
\begin{equation*}
\{\omega^1=dx^1+x^4dx^5,~\omega^2=dx^2,~\omega^3=dx^3\}
\end{equation*}
is \textit{not} integrable since
\begin{equation*}
d\omega^1=dx^4\wedge dx^5\not\equiv 0\hspace{3 mm}mod~\{\omega^1,\omega^2,\omega^3\}
\end{equation*}
and $\mathcal{P}_1$ is generated by $\{\omega^2,\omega^3\}$ hence is, of course, integrable. The character of $\mathcal{P}$ is equal to 1 since, at any point $x_0~,$ the maximal integral manifold is the line defined by 
$x^1=c^1,~x^2=c^2,~x^3=c^3,~x^5=c^5$ for some constants $c^i$ and is a line parallel to the $x^4-$axis. It should be observed that the annihilator $\Sigma_{x_0}$ of $\mathcal{P}_{x_0}$ is $2-$dimensional.

\vspace{2 mm}
\noindent
Let us now replace $\omega^3$ by $\tilde{\omega}^3=dx^3+x^5dx^1~.$ Then
\begin{equation*}
d\tilde{\omega}^3=dx^5\wedge dx^1=dx^5\wedge\omega^1~,
\end{equation*}
hence $\mathcal{P}_1$ is again generated by $\{\omega^2,\omega^3\}$ but is \textit{not} integrable since $d\tilde{\omega}^3\not\equiv 0\hspace{3 mm}mod~\{\omega^2,\omega^3\}~.$ This derived system, whose annihilator is point-wise $3-$dimensional, just admits the $2-$dimensional integral manifolds defined by the equations $x^1=c^1,~x^2=c^2,~x^3=c^3~.$

\vspace{2 mm}
\noindent
We next give a glance at those Pfaffian systems with character equal to 2 and it will just suffice to add one more dimension to our space. Let us therefore consider the system $\mathcal{P}$ generated, on a $6-$space, by
\begin{equation*}
\{\omega^1=dx^1+x^4dx^5,~\omega^2=dx^2+x^5dx^6,~\omega^3=dx^3\}~.    
\end{equation*}
Then the annihilator $\Sigma_{x_0}$ of $\mathcal{P}$ is $3-$dimensional, the derived system $\mathcal{P}_1$ is generated by $\{dx^3\}$ hence is integrable and its rank is two units less than the rank of $\mathcal{P}~.$ Consequently, $\mathcal{P}$ is not integrable. On the other hand, the generators vanish on the affine line defined by the equations
\begin{equation*}
x^1=c^1,~x^2=c^2,~x^3=c^3,~x^5=c^5,~x^6=c^6
\end{equation*}
and, moreover, this line is a maximal integral manifold. The character of $\mathcal{P}$ is therefore equal to $2~.$

\vspace{2 mm}
\noindent
Let us now replace $\omega^3$ by $\tilde{\omega}^3=dx^3+x^5dx^1~.$ Then $\mathcal{P}$ continues to be of character 2 but the derived system $\mathcal{P}_1~,$ generated by $\tilde{\omega}^3~,$ is no longer integrable since $d\tilde{\omega}^3=dx^4\wedge dx^1\not\equiv 0\hspace{3 mm}mod~\{\tilde{\omega}^3\}$ and only admits the $4-$dimensional maximal integral manifolds ($4-$dimensional affine spaces) defined by the equations
\begin{equation*}
x^1=c^1,~x^3=c^3~.
\end{equation*}
The maximal integral element of $\mathcal{P}$ are only $1-$dimensional and its maximal integral manifolds are affine lines.

\vspace{2 mm}
We continue our little game by constructing now a Pfaffian system of character 2 whose derived system has rank three units less than the rank of the given system \textit{i.e.}, we shall construct a singular Pfaffian system with character equal to $2~.$ For this, we simply consider the system $\mathcal{P}$ generated, on a $6-$dimensional space, by
\begin{equation*}
\{\omega^1=dx^1+x^4dx^5,~\omega^2=dx^2+x^5dx^6,~\omega^3=dx^3+x^6dx^4\}~.
\end{equation*}
The maximal integral manifolds of $\mathcal{P}$ are affine lines and a direct calculation shows that $\mathcal{P}_1=0~.$ The rank of the derived system is three units less than its own rank.

\vspace{2 mm}
\noindent
As for the \textit{gender} of a Pfaffian system, the first example considered above illustrates the content of the Lemma 4 since the local section $\omega^1+x^2\omega^3$ has its gender equal to $2~.$

\vspace{2 mm}
Let us finally terminate our discussion by exhibiting three \textit{Jordan-Hölder resolutions}. To begin with, we continue to consider the first example in $5-$space and observe that the pseudo-group $\Gamma$ of all the local automorphisms of $\mathcal{P}$ leaves also invariant the derived system $\mathcal{P}_1$ hence, inasmuch, its annihilator $\Sigma_1$ whose point-wise dimension is equal to 3. Needless to say that the integral manifolds of $\Sigma_1$ are all parallel affine three dimensional sub-spaces in $5-$space. Moreover, $\Gamma$ also leaves invariant the annihilator $\Sigma\subset\Sigma_1~.$ We denote by $\Gamma_1$ the sub-pseudo-group of all the elements $\varphi~,$ belonging to $\Gamma~,$ that \textit{maintain invariant} each connected integral manifold of $\Sigma_1$ \textit{i.e.}, $\varphi$ transforms the points of such an integral manifold into points of the same manifold. Equivalently, this simply means that the above mentioned integral manifolds are the intransitivity classes of $\Gamma_1~.$ We can now factor locally, the given $5-$space, modulo these intransitivity classes and obtain, as quotient, an open set \textit{U} in numerical $2-$space. Furthermore, the Pfaffian system $\mathcal{P}$ being invariant under $\Gamma_1~,$ factors to the quotient and yields a quotient Pfaffian system $\overline{P}$ that is generated by $\{dx^1,dx^2\}$ hence its maximal integral leaves are the points of \textit{U}. Fixing a point $y_0\in U$ and denoting by $Y_0$ the inverse image $q^{-1}(y_0),$ where \textit{q} is the quotient projection modulo the above intransitivity classes, we can now restrict the system $\mathcal{P}$ to the sub-manifold $Y_0~,$ this restricted system $\mathcal{Q}$ being generated by the $1-$form $dx^1+x^4dx^5.$ Since $dim~Y_0=3$ and since the generating form of the system is of maximum Darboux class equal to $3~,$ the restricted system is a Darboux system of rank 1 and only admits $1-$dimensional integral manifolds that, in fact, are the affine lines in $5-$space that constitute the maximal integral manifolds of $\mathcal{P}~.$ The integration of this system is consequently achieved via the Jordan-Hölder method since the quotient pseudo-group $\Gamma/\Gamma_1$ is equal to the set of all the local transformations in \textit{U}, hence is simple.

\vspace{2 mm}
We shall continue with the same Pfaffian system but apply the Jordan-Hölder methodology in a different manner. The character of the system being equal to 1 and, consequently, the maximal integral manifolds being $1-$dimensional curves that actually constitute a foliation, we consider the sub-pseudo-group $\Gamma_1$ as being composed by all those local transformations that keep invariant the connected maximal integral curves, in much the same way as considered previously for the integrals of the first derived system. As before, let \textit{U} be a local quotient of the $5-$space modulo the above integral curves. Then the dimension of \textit{U} is equal to 4 and we can, inasmuch, factor locally the system $\mathcal{P}$ to the space \textit{U}. Since the annihilator, at each point, of the system $\mathcal{P}$ contains the tangent space to the integral curve at that point, we infer that the quotient system $\mathcal{Q}$ is still of rank 3 and is generated by $\{dx^2,dx^3\}$ together with the quotient of $\omega^1.$ However, since the dimension of \textit{U} is equal to 4 and the quotient system $\mathcal{Q}$ has rank equal to $3~,$ we infer that it is necessarily integrable, its maximal integrals being again curves. Proceeding as prescribed in the Jordan-Hölder scheme, we fix a certain, though arbitrary, integral curve of $\mathcal{Q}$ and consider the inverse image of this curve (more precisely, of the image of this curve) modulo the established quotient projection described above. In other terms, we consider the surface \textit{Y} obtained as the union of all the integral lines of $\mathcal{P}$ that project onto points of the selected (though fixed) integral curve of $\mathcal{Q}~.$ Knowing \textit{a priori} that the integral curves are parallel affine lines, we shall thus obtain a \textit{ruled surface}. Let us now restrict $\mathcal{P}$ to the above surface and denote this restricted Pfaffian system by $\overline{\mathcal{P}}.$ Since it must vanish on the tangent spaces to the integral curves, its rank must be equal to 1 and its integral curves are precisely the integral curves of $\mathcal{P}$ contained in the above surface. The Jordan-Hölder method is therefore a most outstanding integration procedure since it enables us to obtain the maximal integral manifolds of the initially given system by a very precise and constructive technique. It should be noted, inasmuch, that the whole process can of course be performed without knowing \textit{a priori} the nature of the maximal integral manifolds that we are looking for (curves in our present context), this being most often the case and, in fact, the sole purpose of an integration method. Moreover, for the system discussed here, the usual \textit{quadrature} technique (Calculus I), will enable us to find a first integral for the integral curves.  

\vspace{2 mm}
Let us now consider the last example given above, where the three forms $\omega^i~,$ generating the system, are of Darboux class equal to $3~.$ We proceed as above and consider the sub-pseudogroup $\Gamma_1$ composed by all those elements of $\Gamma$ that preserve the maximal integral manifolds of $\mathcal{P},$ here again $1-$dimensional curves since the character of $\mathcal{P}$ is equal to $2~.$ The local quotient space \textit{U} is, in this case, $5-$dimensional and the resulting quotient system $\mathcal{Q}$ continues to be of rank 3 , the point-wise dimension of its annihilator being thereafter equal to $2~.$ However, this quotient system acquires character equal to $1~.$\footnote{somehow, the reduction in the dimension squeezes the system into a smaller ambient space.} Taking a $1-$dimensional integral curve $\gamma$ of $\mathcal{Q}~,$ considering the inverse image $Y=q^{-1}(im~\gamma)$ defined in the initially given $6-$space and restricting $\mathcal{P}$ to the surface \textit{Y}, we finally obtain a $rank~1$ Pfaffian system whose integral curves are precisely the integral curves of $\mathcal{P}$ contained in \textit{Y} and for which a first integral can be obtained by a quadrature.

\vspace{2 mm}
We should at present confess to the reader that the above discussion is actually \textit{a big fake}. In fact, the aim of the Jordan-Hölder procedure is, of course, the obtainment of maximal integral manifolds of the given system whereas our procedure was to exhibit this technique, on examples, with the help of the \textit{a priori} knowledge of the solutions. Nevertheless, we believe that the reader is able to devise much more sophisticated examples (or look for them in Cartan's \textit{Mémoires}) where the true aim will be, of course, the determination of maximal integral manifolds not known in advance. One last word is due. Élie Cartan never claimed that the Jordan-Hölder procedure was an easy matter. It can in fact become extremely involved, reason for which Cartan \textit{never} mentioned neither Jordan nor Hölder.

\vspace{2 mm}
A very interesting type of non-integrable Pfaffian systems is given by the \textit{Flag Systems} considered in various occasions by Élie Cartan (\textit{e.g.}, \cite{Cartan1901}, \cite{Cartan1910}, \cite{Cartan1911}). We claim that a flag system $\mathcal{F}$ with null characteristics has always its character equal to $n-r-1~,$ where $n=dim~M$ and $r=rank~\mathcal{F}~.$ In fact, proceeding as indicated in \cite{Kumpera2014} and replacing the descending chain of \textit{derived} systems by the more appropriate \textit{structured multi-fibration} (6), in the above reference, where each successive derived system is replaced by an equivalent system that now has also null characteristics, we see that the maximal dimensional integral manifolds, at each stage, are the integral curves contained in the inverse image of the integral curves of the Darboux system, these last curves being also the maximal dimensional integral manifolds. Our claim then follows since the point-wise dimension of the annihilator distribution $\mathcal{F}^{\perp}$ is precisely $n-r~.$ More generally, we can consider the \textit{Multi-Flag Systems} or, still better, the \textit{Truncated Multi-Flag Systems} that have much relevance in the study of under-determined ordinary differential equations in what concerns the \textit{Monge property} (\cite{Kumpera2002}). We leave the details to the reader and, in particular, the calculation of the characters for such systems. It should be noted that the difference between flag and multi-flag systems can be recognized by looking at their local equivalence pseudo-groups that are non-isomorphic. Whereas for flag systems the Monge condition is a necessary and sufficient condition for the validity of the Monge property (\cite{Cartan1911}), for multi-flag systems the condition is only sufficient (\cite{Kumpera2002}). 

\vspace{2 mm} 
Last but not least, let us give a glance at a rare jewel left to us by Élie Cartan namely, his article on Galois theory (\cite{Cartan1938/39}), where he bases his argumentation on the classical Picard-Vessiot theory for linear differential systems and shows how to integrate such systems via the Jordan-Hölder method.

\section{Equations on Jet bundles.}
On the base manifold \textit{M}, we consider the $k-$th order jet space $J_kTM$ of all the $k-$jets of local sections of the tangent bundle
\begin{equation*}
\pi:TM~\longrightarrow~M~,
\end{equation*}
\textit{i.e.}, $k-$jets of local vector fields on \textit{M}. Then $\alpha_k:J_kTM~\longrightarrow~M~,$ $\alpha_k$ denoting the \textit{source} map, becomes again a vector bundle, its linear operations being simply the extensions, to $k-$jets, of the linear operations in \textit{TM}. A linear partial differential equation of order \textit{k} is, by definition, a vector sub-bundle $\mathcal{R}$ of $J_kTM.$ The dual (covariant) version reads as follows. We take the dual bundle $T^*M~,$ consider the jet space $J_kT^*M$ and then put in evidence a vector sub-bundle $\mathcal{S}.$ It should be remarked that $J_kT^*M\equiv(J_kTM)^*.$ Given a linear differential equation $\mathcal{R}\subset J_kTM,$ its annihilator becomes a vector sub-bundle of $J_kT^*M$ \textit{i.e.}, $\mathcal{R}^{\perp}=\mathcal{S}\subset J_kT^*M$ and, inasmuch, we can consider the annihilator, in $J_kTM,$ of any sub-bundle $\mathcal{S}\subset J_kT^*M.$ We denote $n=dim~M$ and $n_k=dim~J_kTM.$ On account of the previous remarks, it also follows that $n_k=dim~J_kT^*M.$ Taking a local coordinate system $(U,(x^i))$ in the base space \textit{M}, we derive corresponding coordinate systems for $J_kTM$ and $J_kT^*M,$ defined on the respective inverse images of the open set \textit{U}, and any linear partial differential equation will be determined, locally by a system of linear equations with variable coefficients. If $rank~\mathcal{R}=\rho$ denotes the point-wise dimension of the fibres, then $rank~\mathcal{S}=n_k-\rho$ and the number of independent linear equations defining locally the system is also equal to $n_k-\rho~.$

\vspace{2 mm}
As is usual, we indicate by $\mathcal{C}_k$ the canonical contact system defined on $J_kTM$ and, just to mark the distinction, by $\mathcal{C}_k^*$ the canonical contact system defined on $J_kT^*M.$ We recall that a local section in the $k-$jets is holonomic \textit{i.e.}, results from a section in the base tangent or co-tangent bundles if and only if it is annihilated by the contact Pfaffian syatem and, furthermore, a local transformation in the jet space is the prolongation of a base space transformation (in $TM$ or $T^*M$) if and only if it preserves the corresponding contact structure. In the present case, we shall restrict our attention to the pseudo-group of all the local transformations that are vector bundle morphisms. This is the starting point of the Picard-Vessiot theory and the set of all the vector bundle automorphisms of the equation is called its \textit{Galois group}. Most often, this pseudo-group operates non-transitively on the equation though the $(k+1)-$st order groupoid of all the $(k+1)-$jets that preserve the equation, via the \textit{semi-holonomic} action, is always transitive on account of linearity. This action is defined as follows. Given a jet $j_{k+1}\varphi(v)$ on $TM$ or $T^*M,$ we consider the flow 
\begin{equation*}
j_k\varphi:w~\longmapsto~j_k\varphi(w)
\end{equation*}
and subsequently take the jet $j_1(j_k\varphi)(X),~\alpha_k(X)=v~.$ Restricting the contact structures to the corresponding equations, we obtain the Pfaffian systems $\overline{\mathcal{C}}_k$ and $\overline{\mathcal{C}}_k^*$ associated to the linear equations. These systems are not integrable but their $n-$th dimensional integral manifolds transversal to the $\alpha_k-$fibres are precidely the images of the $k-$jets of solutions. On account of the transitivity at the groupoid level, we can apply the Jordan-Hölder integration method as outlined in \cite{Kumpera2016}. When the linear pseudo-group operates transitively on the equations, the integration process can be achieved by just calling out the simple linear Lie groups. We finally observe that Picard (\cite{Picard1896}), Vessiot (\cite{Vessiot1904},\cite{Vessiot1912}) and also Drach (\cite{Drach1898}) only considered first order linear systems. Many years earlier, Sophus Lie studied linear differential equations invariant under linear groups with constant coefficients (\cite{Lie1876})  As for Cartan's \textit{Mémoire}, he shows us much, in \cite{Cartan1938/39}, via examples.

\bibliographystyle{plain}
\bibliography{references}

\end{document}